\documentclass[10pt,a4paper]{article}

\usepackage{authblk}
\usepackage{graphicx}\usepackage[latin1]{inputenc}
\usepackage[english]{babel}
\usepackage{amssymb,amsmath}
\usepackage{a4wide,bm}
\usepackage{enumitem}
\usepackage{tikz}
\usepackage{soul}
\usepackage{pifont}
\newcommand{\cmark}{\ding{51}}%
\newcommand{\xmark}{\ding{55}}%
\usepackage{hyperref}

\newcommand{\parder}[2]{{\frac{\partial #1}{\partial #2}}}

\def\epsilon{\varepsilon}
\def\im{\mathrm{i}}
\def\realpart{\mathop{\rm Re}\nolimits}
\def\imaginary{\mathop{\rm Im}\nolimits}

\def\Jscr{\mathcal{J}}
\def\Lscr{\mathcal{L}}

\def\Oscr{\mathcal{O}}

\def\Rscr{\mathcal{R}}
\def\Uscr{\mathcal{U}}
\def\Vscr{\mathcal{V}}

\def\interi{\mathbb{Z}}
\def\reali{{\mathbb{R}}}
\def\complessi{{\mathbb{C}}}
\def\toro{\mathbb{T}}

\newtheorem{theorem}{Theorem}[section]
\newtheorem{theorem*}{Theorem}

\newtheorem{remark}{Remark}[section]

\title{Continuation of spatially localized periodic solutions \\in
  discrete NLS lattices via normal forms}

\author[2,3]{Veronica Danesi}
\author[1,3]{Marco Sansottera}
\author[1,3]{Simone Paleari}
\author[1,3,*]{Tiziano Penati}

\affil[1]{\small Department of Mathematics ``F.Enriques'',
  University of Milan, via Saldini 50, 20133, Milan, Italy}

\affil[2]{Department of Mathematics, University of Rome ``Tor Vergata'', via della Ricerca Scientifica 1, 00133, Rome, Italy}

\affil[3]{GNFM (Gruppo Nazionale di Fisica Matematica) -- Indam
  (Istituto Nazionale di Alta Matematica ``F. Severi''), Roma, Italy}

\affil[*]{Corresponding author: tiziano.penati@unimi.it}

\begin{document}

\maketitle

\begin{abstract}
  We consider the problem of the continuation with respect to a small
  parameter $\epsilon$ of spatially localised and time periodic solutions in
  1-dimensional dNLS lattices, where $\epsilon$ represents the strength of
  the interaction among the sites on the lattice. Specifically, we
  consider different dNLS models and apply a recently developed normal
  form algorithm in order to investigate the continuation and the
  linear stability of degenerate localised periodic orbits on lower and
  full dimensional invariant resonant tori.  We recover results already
  existing in the literature and provide new insightful ones, both for
  discrete solitons and for invariant subtori.
\end{abstract}


\section{Introduction}\label{s:0}

The discrete nonlinear Schr{\"o}dinger (dNLS) equation is a
paradigmatic physical model in many different areas of Physics, such
as condensed matter, photonic crystals and waveguides.  Indeed, it is
a widely investigated nonlinear lattice model (see, e.g.,
\cite{HenT99, EilJ03, Kev_book09, PalP19, Mal20}) thanks to the
possibility to mathematically combine rigorous different analytic
approaches (e.g., perturbative, variational, spectral, etc.) and to
accurately explore its dynamical features with reliable numerical
methods.

The aim of this work is to investigate the existence of spatially
localised and time periodic solutions in dNLS models, i.e.,
\begin{equation} 
  \label{e.dNLS.eqs}
  \im\dot\psi_j = \psi_j -\epsilon\left(\Lscr\psi\right)_j +
  \gamma\psi_j|\psi_j|^2 \ ,\quad j\in\Jscr\ ,
\end{equation}
where $\Jscr$ is a suitable finite set of indices, with $|\Jscr|=n$,
$\epsilon\in\reali$ is a small parameter (since we focus on the
so-called anticontinuum limit), $\psi_j$ are complex functions and
the linear operator $\Lscr$ reads
\begin{equation}
  \Lscr\psi = \sum_{l=1}^d\kappa_l (\Delta_l\psi)\ , \quad
  (\Delta_l\psi)_j := \psi_{j+l}-2\psi_j+\psi_{j-l}\ ,\qquad 
  j\in\Jscr\ ,
\end{equation}
where $\kappa_l$ are real parameters describing the
$l$-nearest-neighbours coupling.  We also consider periodic boundary
conditions as $\Jscr$ is taken finite\footnote{One can also consider
  boundary conditions vanishing at infinity as
  $\psi\in\ell^2(\complessi)$ in the case of infinite $\Jscr$.  This
  case is not properly covered by the normal form technique here
  proposed: the formal algorithm applies, but the analytic estimates
  needs to be extended.}.

Since we are mainly interested in periodic solutions $\psi(t)$
of~\eqref{e.dNLS.eqs} which are spatially localised on a subset of the
lattice, we excite only $m$ sites (with $m<n=|\Jscr|$) and introduce
the subset
\begin{equation}
  \label{e.S}
  S:=\{j_1,\ldots,j_m\}\subset\Jscr\ ,
\end{equation}
where we stress that the indexes in $S$ do not have to be necessarily
consecutive. Hence, we are including also configurations
where the localisation of the amplitude (hence of the energy), is
clustered, with \emph{holes} separating the different clusters along the
lattice.

For $\epsilon=0$, the unperturbed excited oscillators
$\{\psi_j^{(0)}\}_{j\in\Jscr}$ are set in \emph{complete
  resonance} in order to have a periodic flow on a resonant torus
$\toro^m$. The typical choice exploited in the literature (see, e.g.,
\cite{Kap01,KapK01,Kev09,PelKF05,PelKF05b,QX07,Pan11,PenSPKK18,ParKS20})
is the $1:\text{\textendash}:1$ resonance, obtained by choosing a common
amplitude $R$, or a common frequency $\omega$, for all the
$\{\psi_j^{(0)}\}_{j\in\Jscr}$; in this way the solution takes the
form of the so-called \emph{rotating frame ansatz}
\begin{equation}
  \label{e.ansatz}
  \psi^{(0)}(t) = e^{-\im\omega t}\phi^{(0)}\ ,\qquad
  t\in[0,T=2\pi/\omega]\ ,
\end{equation}
where the unperturbed spatial profile $\phi^{(0)}$ and the frequency $\omega$ 
read
\begin{equation}
  \label{e.torus}
  \phi_j^{(0)} = \begin{cases} R e^{\im \theta_j} \, , &j\in S\ , \\
    0 \,, &j\in \Jscr\setminus S\ ,
  \end{cases}
  \qquad \omega(R) = 1+\gamma R^2\ .
\end{equation}

\smallskip

The system of equations~\eqref{e.dNLS.eqs} can be written in Hamiltonian form as
$$
\im\dot\psi_j =
\parder{H}{\overline\psi_j}\ ,\qquad\hbox{with}\quad
H=H_0+\epsilon H_1\ ,
$$
where
\begin{equation}
  \label{e.KdNLS} 
  H_0 = \sum_{j\in\Jscr}|\psi_j|^2 +
  \frac\gamma2\sum_{j\in\Jscr}|\psi_j|^4\ , \qquad
  H_1 = \sum_{l=1}^d\kappa_l\sum_{j\in\Jscr}|\psi_{j+l}-\psi_j|^2 \ .
\end{equation}

The approach we here propose entirely works at the Hamiltonian level:
we exploit the nearly integrable structure of the problem and perform
a sequence of near to the identity canonical transformation that puts
the Hamiltonian in normal form up to a certain order in the small
parameter $\epsilon$.  Precisely, the normal form at order $r$ reads
$H^{(r)}= K^{(r)} + \epsilon^{r+1} \Rscr^{(r+1)}$ where $K^{(r)}$ is
already in normal form, while $\Rscr^{(r+1)}$ is the remainder.  The
existence and linear stability of the so-called \emph{discrete
  solitons} (consecutive sites) or \emph{multi-pulse discrete
  solitons} (nonconsecutive sites) is investigated with a resonant
normal form algorithm, recently developed for fully resonant maximal
and lower dimensional invariant tori, see \cite{PenSD18,SanDPP20}.
Discrete solitons correspond to time periodic solutions which, at
$\epsilon=0$, belong to a fully resonant low-dimensional torus $\toro^m$
of the unperturbed Hamiltonian $H_0$; as $\epsilon>0$, generically only a
finite number of the periodic orbits survive to the breaking of the
resonant torus, and they turn out to be spatially localised on the few
variables defining $\toro^m$.  The relative equilibria $x^*$ of
$K^{(r)}$ provide accurate approximations of the periodic orbits we
are looking at. Moreover, the linearization of $K^{(r)}$ allows to
investigate the approximate linear stability of the periodic
orbits. Eventually, the effective linear stability can be derived from
Theorem~2.3 of~\cite{SanDPP20}, which exploits classical results on
perturbations of the resolvent.

The present work focuses specifically on periodic solutions which are
at leading order degenerate. Such a situation is natural in the study
of spatially localised time periodic solutions, when nonconsecutive
sites/oscillators are excited in the unperturbed model (see, e.g.,
\cite{PenSPKK18,PenKSKP19,PenDP20}), thus leading generically to
$d$-parameter families of critical points of the time average $\langle
H_1 \rangle$; hence continuation cannot be easily obtained via
implicit function theorem. Sometimes degeneracy occurs also in
topologically isolated solutions, due to the semidefinite nature of
the critical points.

Our normal form approach also allows to enrich the
investigation of time periodic localised structures to 2-dimensional 
subtori foliated by periodic orbits: this
naturally occurs in~\eqref{e.dNLS.eqs}, when the resonances among the
excited frequencies of the unperturbed dynamics differ from the
$1:\text{\textendash}:1$, because different amplitudes have been chosen for the
selected sites. 

Indeed, considering the $1:\text{\textendash}:1$ resonance the Hamiltonian field $X_H$ is
parallel to the generator $X_P$ of the symmetry $X_H(\psi^{(0)}(t))=\omega 
X_P(\psi^{(0)}(t))$ for $t\in[0,T]$ with $P:=\sum_{j\in\Jscr}|\psi_j|^2$.  
The ansatz~\eqref{e.ansatz} takes advantage of this peculiarity, as 
$\psi^{(0)}(t)$ coincides with the orbit of the symmetry $e^{\im\varphi}$ 
($\varphi\in\reali$) acting on the surface of constant energy. Instead, for a 
generic resonance the symmetry vector field $X_P$ is transversal to
$X_H$, hence a given periodic orbit is transported by the action of
the symmetry group and one has to deal with a 2-dimensional 
resonant subtorus.

Another benefit of our approach is that it allows to investigate the long-time 
stability of such periodic solutions: indeed this preliminary normal form 
construction, at a suitable order $r$ depending on the \emph{degeneracy 
degree} of the orbit, gives the linearization around the
approximate solution the right shape for a subsequent stability
analysis, e.g., with perturbation methods like Birkhoff normal
forms (see also \cite{Gio12,Bru20} for related studies).

Moreover, the constructive algorithm allows to approximate the
periodic orbit with arbitrary precision, though some additional
computational effort is unavoidable in order to increase the desired
order of accuracy. This is shown in Section \ref{s:3}, where we have
numerically explored as a case study a multi-pulse discrete soliton in
the standard dNLS model. The numerical simulations of the normal form
at orders $r=2$ and $r=3$ support what is theoretically predicted in
terms of accuracy of the approximated solutions and of their linear
stability.

\medskip

The goal of the present work is to illustrate with significant examples how to 
investigate the dynamics of discrete solitons via normal form in the presence 
of degeneracy.  In order to keep the presentation as simple as possible we 
consider models where an explicit normal form up to order $r=3$ (implemented 
with Mathematica) is sufficient.
We report hereafter a summary of the results obtained for each model 
considered, that is schematically represented in a plot depicting the 
configuration of the discrete soliton to be continued together with the 
geometry of the near-neighbour interactions.  The detailed study is deferred to 
section~\ref{s:2}.

\bigskip

\noindent
\textbf{Two-sites multi-pulse discrete solitons for dNLS $S=\{-1,1\}$}

\begin{center}
  \begin{tikzpicture}[scale = 0.5]
    \begin{scope} 
      \foreach \name/\posX/\posY in {A/0/0,B/1/0,C/2/0,D/3/0,E/4/0,F/5/0,G/6/2,H/7/0,I/8/2,L/9/0,M/10/0,N/11/0, O/12/0,P/13/0,Q/14/0}
      \node (\name) at (\posX,\posY) {};
      \draw [dotted, black!50] (-1,0) -- (14,0);
      \draw [dotted, black!50] (-1,2) -- (14,2);
      \draw [dotted, black!50] (6,2) -- (6,0);
      \draw [dotted, black!50] (8,2) -- (8,0);
      \draw (B) circle (2mm) \foreach \n in {C,D,E,F,G,H,I,L,M,N,O,P} {(\n) circle (2mm)};
      \fill [black] \foreach \n in {G,I} {(\n) circle (2.0mm)};
      \fill [white] \foreach \n in {B,C,D,E,F,H,L,M,N,O,P} {(\n) circle (1.9mm)};
      \draw [red] (A) -- (B) -- (C) -- (D) -- (E) -- (F) -- (G)-- (H) -- (I) -- 
      (L) -- (M) -- (N) -- (O) -- (P) -- (Q) ;
      \node[label=below:{\footnotesize-6}] at (B) {\phantom{$\times$}};
      \node[label=below:{\footnotesize-5}] at (C) {\phantom{$\times$}};
      \node[label=below:{\footnotesize-4}] at (D) {\phantom{$\times$}};
      \node[label=below:{\footnotesize-3}] at (E) {\phantom{$\times$}};
      \node[label=below:{\footnotesize-2}] at (F) {\phantom{$\times$}};
      \node[label=below:{\footnotesize 0}] at (H) {\phantom{$\times$}};
      \node[label=below:{\footnotesize 2}] at (L) {\phantom{$\times$}};
      \node[label=below:{\footnotesize 3}] at (M) {\phantom{$\times$}};
      \node[label=below:{\footnotesize 4}] at (N) {\phantom{$\times$}};
      \node[label=below:{\footnotesize 5}] at (O) {\phantom{$\times$}};
      \node[label=below:{\footnotesize 6}] at (P) {\phantom{$\times$}};
      \node[label=below:{\footnotesize-1}] at (6,0) {\phantom{$\times$}};
      \node[label=below:{\footnotesize1}] at (8,0) {\phantom{$\times$}};
      \node[label=left:{\footnotesize $0$}] at (-1,0) {};
      \node[label=left:{\footnotesize $I^*$}] at (-1,2) {};
    \end{scope}
  \end{tikzpicture}
\end{center}

\begin{center}
  \begin{tabular}{c|c|c|c}
    Normalisation order & Candidate for continuation & Continuation & 
    Stability \\
    & $(q_2)$ & & \\
    \hline
    $r=1$ & $\vartheta\in\toro$ & ? & ?\\
    \hline
    $r=2$ & $0$ & \cmark & unstable\\
    & $\pi$ & \cmark & stable
  \end{tabular}
\end{center}

\bigskip
\noindent
\textbf{Three-sites multi-pulse discrete solitons for dNLS  $S=\{-2,-1,1\}$}

\begin{center}
  \begin{tikzpicture}[scale = 0.5]
    \begin{scope} 
      \foreach \name/\posX/\posY in {A/0/0,B/1/0,C/2/0,D/3/0,E/4/0,F/5/2,G/6/2,H/7/0,I/8/2,L/9/0,M/10/0,N/11/0,O/12/0,P/13/0,Q/14/0}
      \node (\name) at (\posX,\posY) {};
      \draw [dotted, black!50] (-1,0) -- (14,0);
      \draw [dotted, black!50] (-1,2) -- (14,2);
      \draw [dotted, black!50] (5,2) -- (5,0);
      \draw [dotted, black!50] (6,2) -- (6,0);
      \draw [dotted, black!50] (8,2) -- (8,0);
      \draw (B) circle (2mm) \foreach \n in {C,D,E,F,G,H,I,L,M,N,O,P} {(\n) circle (2mm)};
      \fill [black] \foreach \n in {F,G,I} {(\n) circle (2.0mm)};
      \fill [white] \foreach \n in {B,C,D,E,H,L,M,N,O,P} {(\n) circle (1.9mm)};
      \draw [red] (A) -- (B) -- (C) -- (D) -- (E) -- (F) -- (G)-- (H) -- (I) -- (L) -- (M) -- (N) -- (O) -- (P) -- (Q) ;
      \node[label=below:{\footnotesize-6}] at (B) {\phantom{$\times$}};
      \node[label=below:{\footnotesize-5}] at (C) {\phantom{$\times$}};
      \node[label=below:{\footnotesize-4}] at (D) {\phantom{$\times$}};
      \node[label=below:{\footnotesize-3}] at (E) {\phantom{$\times$}};
      \node[label=below:{\footnotesize 0}] at (H) {\phantom{$\times$}};
      \node[label=below:{\footnotesize 2}] at (L) {\phantom{$\times$}};
      \node[label=below:{\footnotesize 3}] at (M) {\phantom{$\times$}};
      \node[label=below:{\footnotesize 4}] at (N) {\phantom{$\times$}};
      \node[label=below:{\footnotesize 5}] at (O) {\phantom{$\times$}};
      \node[label=below:{\footnotesize 6}] at (P) {\phantom{$\times$}};
      \node[label=below:{\footnotesize-2}] at (5,0) {\phantom{$\times$}};
      \node[label=below:{\footnotesize-1}] at (6,0) {\phantom{$\times$}};
      \node[label=below:{\footnotesize1}] at (8,0) {\phantom{$\times$}};
      \node[label=left:{\footnotesize $0$}] at (-1,0) {};
      \node[label=left:{\footnotesize $I^*$}] at (-1,2) {};
    \end{scope}
  \end{tikzpicture}
\end{center}

\begin{center}
  \begin{tabular}{c|c|c|c}
    Normalisation order & Candidate for continuation & Continuation & Stability \\
    & $(q_2,q_3)$ & & \\
    \hline
    $r=1$ & $(0,\ \vartheta)\in\toro$ & ? & ?\\
    & $(\pi,\ \vartheta)\in\toro$ & ? & ?\\
    \hline
    $r=2$ & $(0,\ 0)$ & \cmark & unstable\\
    & $(0,\ \pi)$ & \cmark & stable\\
    & $(\pi,\ 0)$ & \cmark & unstable\\
    & $(\pi,\ \pi)$ & \cmark & unstable
  \end{tabular}
\end{center}

\bigskip
\noindent
\textbf{Four-sites vortex-like structures in a Zigzag model}

\vspace{-5pt}
\begin{center}
  \begin{tikzpicture}[scale = 0.5]
    \begin{scope} 
      \foreach \name/\posX/\posY in {C/3/0,D/4.5/0,E/6/0,F/7.5/0,G/9/1.,H/10.5/1.,I/12/0,L/13.5/0,M/15/0,N/16.5/0,
        T/3.5/1.5,U/5./1.5,V/6.5/1.5,W/8./1.5,Y/9.5/2.5,X/11/2.5,Z/12.5/1.5,J/14/1.5,K/15.5/1.5,a/17./1.5}
      \node (\name) at (\posX,\posY) {};
      \draw (D) circle (2mm) \foreach \n in {D,E,F,G,H,I,L,M,U,V,W,Y,X,Z,J,K} {(\n) circle (2mm)};
      \fill [black] \foreach \n in {G,H,Y,X} {(\n) circle (2.0mm)};
      \draw [red]  (C) -- (D) -- (E) -- (F) -- (G)-- (H) -- (I) -- (L) -- (M) -- (N)  ;
      \draw [red]  (T) -- (U) -- (V) -- (W) -- (Y)-- (X) -- (Z) -- (J) -- (K) -- (a)  ;
      \draw [red]   (D) --(U) (E) -- (V) (F) -- (W) (G)-- (Y) (H) -- (X) (I) -- (Z) (L) -- (J) (M) -- (K) ;
      \draw [red]    (U) -- (E)  (V) -- (F) (W) -- (G) (Y) -- (H) (X) -- (I) (Z) -- (L) (J) -- (M)   ;
      \node[label=below:{\footnotesize-6}] at (D) {\phantom{$\times$}};
      \node[label=below:{\footnotesize-4}] at (E) {\phantom{$\times$}};
      \node[label=below:{\footnotesize-2}] at (F) {\phantom{$\times$}};
      \node[label=below:{\footnotesize 0}] at (G) {\phantom{$\times$}};
      \node[label=below:{\footnotesize 2}] at (H) {\phantom{$\times$}};
      \node[label=below:{\footnotesize 4}] at (I) {\phantom{$\times$}};
      \node[label=below:{\footnotesize 6}] at (L) {\phantom{$\times$}};
      \node[label=below:{\footnotesize 8}] at (M) {\phantom{$\times$}};
      \node[label=above:{\footnotesize-5}] at (U) {\phantom{$\times$}};
      \node[label=above:{\footnotesize-3}] at (V) {\phantom{$\times$}};
      \node[label=above:{\footnotesize-1}] at (W) {\phantom{$\times$}};
      \node[label=above:{\footnotesize 1}] at (Y) {\phantom{$\times$}};
      \node[label=above:{\footnotesize 3}] at (X) {\phantom{$\times$}};
      \node[label=above:{\footnotesize 5}] at (Z) {\phantom{$\times$}};
      \node[label=above:{\footnotesize 7}] at (J) {\phantom{$\times$}};
      \node[label=above:{\footnotesize 9}] at (K) {\phantom{$\times$}};
    \end{scope}
  \end{tikzpicture}
\end{center}

\begin{center}
  \begin{tabular}{c|c|c|c}
    Normalisation order & Candidate for continuation & Continuation & Stability \\
    & $(q_2,q_3, q_4)$ & & \\
    \hline
    $r=1$ & $(0,\ 0,\ 0)$ & \cmark & ?\\
    & $(0,\ 0,\ \pi)$ & \cmark & unstable\\
    & $(\pi,\ 0,\ 0)$ & \cmark & unstable\\
    & $(\pi,\ 0,\ \pi)$ & ? & ?\\
    & $(\vartheta,\ \pi,\ \vartheta-\pi)\in\toro$ & ? & ?\\
    & $(\vartheta,\ \pi,\ -\pi)\in\toro$ & ? & ?\\
    \hline
    $r=2$ & $(0,\ 0,\ 0)$ & \cmark & stable\\
    & $(\pi,\ 0,\ \pi)$ & \cmark & unstable\\
    & $(0,\ \pi,\ 0)$ & \cmark & unstable\\
    & $(0,\ \pi,\ \pi)$ & \cmark & unstable\\
    & $(\pi,\ \pi,\ 0)$ & \cmark & unstable\\
    & $(\pi,\ \pi,\ \pi)$ & \cmark & unstable
  \end{tabular}
\end{center}

\bigskip
\noindent
\textbf{Four-sites vortex solutions in a railway dNLS model}

\vspace{-5pt}
\begin{center}
  \begin{tikzpicture}[scale = 0.5]
    \begin{scope} 
      \foreach \name/\posX/\posY in {C/3/0,D/4.5/0,E/6/0,F/7.5/0,G/9/1.,H/10.5/1.,I/12/0,L/13.5/0,M/15/0,N/16.5/0,
        T/3.5/1.5,U/5/1.5,V/6.5/1.5,W/8/1.5,Y/9.5/2.5,X/11/2.5,Z/12.5/1.5,J/14/1.5,K/15.5/1.5,a/17/1.5}
      \node (\name) at (\posX,\posY) {};
      \draw (D) circle (2mm) \foreach \n in {D,E,F,G,H,I,L,M,U,V,W,Y,X,Z,J,K} {(\n) circle (2mm)};
      \fill [black] \foreach \n in {G,H,Y,X} {(\n) circle (2.0mm)};
      \draw [red] (C) -- (D) -- (E) -- (F) -- (G)-- (H) -- (I) -- (L) -- (M) -- (N)  ;
      \draw [red] (T) -- (U) -- (V) -- (W) -- (Y)-- (X) -- (Z) -- (J) -- (K) -- (a)  ;
      \draw [red] (D) -- (U) (E) -- (V) (F) -- (W) (G)-- (Y) (H) -- (X) (I) -- (Z) (L) -- (J) (M) -- (K) ;
      \node[label=below:{\footnotesize-6}] at (D) {\phantom{$\times$}};
      \node[label=below:{\footnotesize-4}] at (E) {\phantom{$\times$}};
      \node[label=below:{\footnotesize-2}] at (F) {\phantom{$\times$}};
      \node[label=below:{\footnotesize 0}] at (G) {\phantom{$\times$}};
      \node[label=below:{\footnotesize 2}] at (H) {\phantom{$\times$}};
      \node[label=below:{\footnotesize 4}] at (I) {\phantom{$\times$}};
      \node[label=below:{\footnotesize 6}] at (L) {\phantom{$\times$}};
      \node[label=below:{\footnotesize 8}] at (M) {\phantom{$\times$}};
      \node[label=above:{\footnotesize-7}] at (U) {\phantom{$\times$}};
      \node[label=above:{\footnotesize-5}] at (V) {\phantom{$\times$}};
      \node[label=above:{\footnotesize-3}] at (W) {\phantom{$\times$}};
      \node[label=above:{\footnotesize-1}] at (Y) {\phantom{$\times$}};
      \node[label=above:{\footnotesize 1}] at (X) {\phantom{$\times$}};
      \node[label=above:{\footnotesize 3}] at (Z) {\phantom{$\times$}};
      \node[label=above:{\footnotesize 5}] at (J) {\phantom{$\times$}};
      \node[label=above:{\footnotesize 7}] at (K) {\phantom{$\times$}};
    \end{scope}
  \end{tikzpicture}
\end{center}

\begin{center}
  \begin{tabular}{c|c|c|c}
    Normalisation order & Candidate for continuation & Continuation & Stability \\
    & $(q_2,q_3, q_4)$ & & \\
    \hline
    $r=1$ & $(0,\ 0,\ 0)$ & \cmark & ?\\
    & $(\pi,\ 0,\ \pi)$ & \cmark & ?\\
    & $(\vartheta,\ \pi,\ -\vartheta)\in\toro$ & ? & ?\\
    & $(\vartheta,\ \pi,\ \vartheta+\pi)\in\toro$ & ? & ?\\
    & $(\vartheta,\ -2\vartheta,\ \vartheta+\pi)\in\toro$ & ? & ?\\
    \hline
    $r=2$ & $(0,\ 0,\ 0)$ & \cmark & ?\\
    & $(\pi,\ 0,\ \pi)$ & \cmark & ?\\
    & $(0,\ 0,\ \pi)$ & \cmark & unstable\\
    & $(0,\ \pi,\ \pi)$ & \cmark & unstable\\
    & $(\pi,\ 0,\ 0)$ & \cmark & unstable\\
    & $(\pi,\ \pi,\ 0)$ & \cmark & unstable\\
    & $(\pi/2,\ \pi,\ 3\pi/2)$ & ? & ?\\
    & $(3\pi/2,\ \pi,\ \pi/2)$ & ? & ?\\
    & $(\vartheta,\ -2\vartheta,\ \vartheta+\pi)\in\toro$ & ? & ?\\
    \hline
    $r=3$ & $(0,\ 0,\ 0)$ & \cmark & stable\\
    & $(\pi,\ 0,\ \pi)$ & \cmark & unstable\\
    & $(\pi/2,\ \pi,\ 3\pi/2)$ & \xmark & \\
    & $(3\pi/2,\ \pi,\ \pi/2)$ & \xmark & \\
    & $(\pi,\ \pi,\ \pi)$ & \cmark & unstable\\
    & $(0,\ \pi,\ 0)$ & \cmark & unstable
  \end{tabular}
\end{center}

\bigskip\pagebreak
\noindent
\textbf{Discrete soliton in dNLS models with purely nonlinear interaction (1:1:1 res.)}

\begin{center}
  \begin{tikzpicture}[scale = 0.5]
    \begin{scope} 
      \foreach \name/\posX/\posY in {A/0/0,B/1/0,C/2/0,D/3/0,E/4/0,F/5/0,G/6/2,H/7/2,I/8/2,L/9/0,M/10/0,N/11/0,
        O/12/0,P/13/0,Q/14/0}
      \node (\name) at (\posX,\posY) {};
      \draw [dotted, black!50] (-1,0) -- (14,0);
      \draw [dotted, black!50] (-1,2) -- (14,2);
      \draw [dotted, black!50] (6,2) -- (6,0);
      \draw [dotted, black!50] (7,2) -- (7,0);
      \draw [dotted, black!50] (8,2) -- (8,0);
      \draw (B) circle (2mm) \foreach \n in {C,D,E,F,G,H,I,L,M,N,O,P} {(\n) circle (2mm)};
      \fill [black] \foreach \n in {G,H,I} {(\n) circle (2.0mm)};
      \fill [white] \foreach \n in {B,C,D,E,F,L,M,N,O,P} {(\n) circle (1.9mm)};
      \draw [red] (A) -- (B) -- (C) -- (D) -- (E) -- (F) -- (G)-- (H) -- (I) -- (L) -- (M) -- (N) -- (O) -- (P) -- (Q) ;
      \node[label=below:{\footnotesize-6}] at (B) {\phantom{$\times$}};
      \node[label=below:{\footnotesize-5}] at (C) {\phantom{$\times$}};
      \node[label=below:{\footnotesize-4}] at (D) {\phantom{$\times$}};
      \node[label=below:{\footnotesize-3}] at (E) {\phantom{$\times$}};
      \node[label=below:{\footnotesize-2}] at (F) {\phantom{$\times$}};
      \node[label=below:{\footnotesize 2}] at (L) {\phantom{$\times$}};
      \node[label=below:{\footnotesize 3}] at (M) {\phantom{$\times$}};
      \node[label=below:{\footnotesize 4}] at (N) {\phantom{$\times$}};
      \node[label=below:{\footnotesize 5}] at (O) {\phantom{$\times$}};
      \node[label=below:{\footnotesize 6}] at (P) {\phantom{$\times$}};
      \node[label=below:{\footnotesize-1}] at (6,0) {\phantom{$\times$}};
      \node[label=below:{\footnotesize 0}] at (7,0) {\phantom{$\times$}};
      \node[label=below:{\footnotesize 1}] at (8,0) {\phantom{$\times$}};
      \node[label=left:{\footnotesize $0$}] at (-1,0) {};
      \node[label=left:{\footnotesize $I^*$}] at (-1,2) {};
    \end{scope}
  \end{tikzpicture}
\end{center}

\begin{center}
  \begin{tabular}{c|c|c|c}
    Normalisation order & Candidate for continuation & Continuation & Stability \\
    & $(q_2,q_3)$ & & \\
    \hline
    $r=1$ & $(0,\ 0)$ & ? & ?\\
    & $(0,\ \pi)$ & ? & ?\\
    & $(\pi,\ 0)$ & ? & ?\\
    & $(\pi,\ \pi)$ & \cmark & unstable\\
    \hline
    $r=2$ & $(0,\ 0)$ & ? & ?\\
    & $(0,\ \pi)$ & ? & ?\\
    & $(\pi,\ 0)$ & ? & ?\\
    \hline
    $r=3$ & $(0,\ 0)$ & \cmark & stable\\
    & $(0,\ \pi)$ & \cmark & unstable\\
    & $(\pi,\ 0)$ & \cmark & unstable
  \end{tabular}
\end{center}

\bigskip
\noindent
\textbf{Discrete soliton in dNLS models with purely nonlinear interaction (2:1:1 res.)}

\begin{center}
  \begin{tikzpicture}[scale = 0.5]
    \begin{scope} 
      \foreach \name/\posX/\posY in {A/0/0,B/1/0,C/2/0,D/3/0,E/4/0,F/5/0,G/6/1,H/7/2.3,I/8/2.3,L/9/0,M/10/0,N/11/0,
        O/12/0,P/13/0,Q/14/0}
      \node (\name) at (\posX,\posY) {};
      \draw [dotted, black!50] (-1,0) -- (14,0);
      \draw [dotted, black!50] (-1,1) -- (14,1);
      \draw [dotted, black!50] (-1,2.3) -- (14,2.3);
      \draw [dotted, black!50] (6,1) -- (6,0);
      \draw [dotted, black!50] (7,2) -- (7,0);
      \draw [dotted, black!50] (8,2) -- (8,0);
      \draw (B) circle (2mm) \foreach \n in {C,D,E,F,G,H,I,L,M,N,O,P} {(\n) circle (2mm)};
      \fill [black] \foreach \n in {G,H,I} {(\n) circle (2.0mm)};
      \fill [white] \foreach \n in {B,C,D,E,F,L,M,N,O,P} {(\n) circle (1.9mm)};
      \draw [red] (A) -- (B) -- (C) -- (D) -- (E) -- (F) -- (G)-- (H) -- (I) -- (L) -- (M) -- (N) -- (O) -- (P) -- (Q) ;
      \node[label=below:{\footnotesize-6}] at (B) {\phantom{$\times$}};
      \node[label=below:{\footnotesize-5}] at (C) {\phantom{$\times$}};
      \node[label=below:{\footnotesize-4}] at (D) {\phantom{$\times$}};
      \node[label=below:{\footnotesize-3}] at (E) {\phantom{$\times$}};
      \node[label=below:{\footnotesize-2}] at (F) {\phantom{$\times$}};
      \node[label=below:{\footnotesize 2}] at (L) {\phantom{$\times$}};
      \node[label=below:{\footnotesize 3}] at (M) {\phantom{$\times$}};
      \node[label=below:{\footnotesize 4}] at (N) {\phantom{$\times$}};
      \node[label=below:{\footnotesize 5}] at (O) {\phantom{$\times$}};
      \node[label=below:{\footnotesize 6}] at (P) {\phantom{$\times$}};
      \node[label=below:{\footnotesize-1}] at (6,0) {\phantom{$\times$}};
      \node[label=below:{\footnotesize 0}] at (7,0) {\phantom{$\times$}};
      \node[label=below:{\footnotesize 1}] at (8,0) {\phantom{$\times$}};
      \node[label=left:{\footnotesize $0$}] at (-1,0) {};
      \node[label=left:{\footnotesize $I^*$}] at (-1,1) {};
      \node[label=left:{\footnotesize $1+2I^*$}] at (-1,2.3) {};	
    \end{scope}
  \end{tikzpicture}
\end{center}

\begin{center}
  \begin{tabular}{c|c|c|c}
    Normalisation order & Candidate for continuation & Continuation & Stability \\
    & $(q_2,q_3)$ & & \\
    \hline
    $r=1$ & $(\vartheta,\ 0)\in\toro$ & \cmark & stable\\
    & $(\vartheta,\ \pi)\in\toro$ & \cmark & unstable
  \end{tabular}
\end{center}

\bigskip
\noindent
\textbf{Discrete soliton in dNLS models with purely nonlinear interaction (2:1:2 res.)}

\begin{center}
  \begin{tikzpicture}[scale = 0.5]
    \begin{scope} 
      \foreach \name/\posX/\posY in {A/0/0,B/1/0,C/2/0,D/3/0,E/4/0,F/5/0,G/6/1,H/7/2.3,I/8/1,L/9/0,M/10/0,N/11/0,
        O/12/0,P/13/0,Q/14/0}
      \node (\name) at (\posX,\posY) {};
      \draw [dotted, black!50] (-1,0) -- (14,0);
      \draw [dotted, black!50] (-1,1) -- (14,1);
      \draw [dotted, black!50] (-1,2.3) -- (14,2.3);
      \draw [dotted, black!50] (6,1) -- (6,0);
      \draw [dotted, black!50] (7,2) -- (7,0);
      \draw [dotted, black!50] (8,1) -- (8,0);
      \draw (B) circle (2mm) \foreach \n in {C,D,E,F,G,H,I,L,M,N,O,P} {(\n) 
        circle 
        (2mm)};
      \fill [black] \foreach \n in {G,H,I} {(\n) circle (2.0mm)};
      \fill [white] \foreach \n in {B,C,D,E,F,L,M,N,O,P} {(\n) circle (1.9mm)};
      \draw [red] (A) -- (B) -- (C) -- (D) -- (E) -- (F) -- (G)-- (H) -- (I) -- (L) -- (M) -- (N) -- (O) -- (P) -- (Q) ;
      \node[label=below:{\footnotesize-6}] at (B) {\phantom{$\times$}};
      \node[label=below:{\footnotesize-5}] at (C) {\phantom{$\times$}};
      \node[label=below:{\footnotesize-4}] at (D) {\phantom{$\times$}};
      \node[label=below:{\footnotesize-3}] at (E) {\phantom{$\times$}};
      \node[label=below:{\footnotesize-2}] at (F) {\phantom{$\times$}};
      \node[label=below:{\footnotesize 2}] at (L) {\phantom{$\times$}};
      \node[label=below:{\footnotesize 3}] at (M) {\phantom{$\times$}};
      \node[label=below:{\footnotesize 4}] at (N) {\phantom{$\times$}};
      \node[label=below:{\footnotesize 5}] at (O) {\phantom{$\times$}};
      \node[label=below:{\footnotesize 6}] at (P) {\phantom{$\times$}};
      \node[label=below:{\footnotesize-1}] at (6,0) {\phantom{$\times$}};
      \node[label=below:{\footnotesize 0}] at (7,0) {\phantom{$\times$}};
      \node[label=below:{\footnotesize 1}] at (8,0) {\phantom{$\times$}};
      \node[label=left:{\footnotesize $0$}] at (-1,0) {};
      \node[label=left:{\footnotesize $I^*$}] at (-1,1) {};
      \node[label=left:{\footnotesize $1+2I^*$}] at (-1,2.3) {};	
    \end{scope}
  \end{tikzpicture}
\end{center}

\begin{center}
  \begin{tabular}{c|c|c|c}
    Normalisation order & Candidate for continuation & Continuation & Stability \\
    & $(q_2,q_3)$ & & \\
    \hline
    $r=1$ & $(\vartheta_1,\ \vartheta_2)\in\toro^2$ & ? & ?\\
    \hline
    $r=2$ & $(\vartheta,\ 0)\in\toro$ & \cmark & unstable\\
    & $(\vartheta,\ \pi)\in\toro$ & \cmark & stable
  \end{tabular}
\end{center}

The paper is structured as follows. In Section~\ref{s:1} we review the
normal form scheme introduced in~\cite{SanDPP20} together with the
main Theorems about continuation of periodic orbit and its linear
stability.  In Section~\ref{s:2} we describe in detail all the above
mentioned examples.  The results have been obtained by
implementing\footnote{The actual code can be found at \url{https://github.com/marcosansottera/periodic_orbits_NF}.}
the normal form algorithm in Mathematica. In Section~\ref{s:3} we
explore numerically the approximate solutions in the case study of
multi-pulse discrete soliton in the standard dNLS
chain. Section~\ref{s:4} is devoted to some final comments.

\section{Theoretical framework}
\label{s:1}

We here briefly recall the normal form scheme presented 
in~\cite{SanDPP20}, so as to make the paper quite self-contained.  We refer
to the quoted paper for all the details.  The main feature that we 
want to stress is that the normal form algorithm is completely 
constructive and can be effectively implemented in a computer algebra system.  
Thus, in a specific application, one can easily check all the 
assumptions in Theorems~\ref{teo:forma-normale-r} and~\ref{t.lin.stab.2}. This 
is what we highlight in section~\ref{s:2} through the examples presented 
before.

\subsection{Preliminary transformations}
The real Hamiltonian~\eqref{e.KdNLS} is written as a function of the complex 
amplitudes $\psi_j$.  Introducing the real canonical variables
\begin{equation}
\label{e.real.coord}
x_j = \frac\im{\sqrt2}(\psi_j-\bar\psi_j)=-\sqrt2\imaginary(\psi_j)\ ,\quad\quad
y_j=\frac1{\sqrt2}(\bar\psi_j+\psi_j=\sqrt2\realpart(\psi_j)\ ,
\end{equation}
the Hamiltonian reads again $H=H_0+\epsilon
H_1$ with
\begin{equation}
  \label{e.ex.dnls}
  \begin{aligned}
    H_{0} &= \sum_{j\in\Jscr} \frac{1}{2} (x_j^2 + y_j^2) +
    \frac{\gamma}{8} (x_j^2 + y_j^2)^2 \\
    H_1 &= \sum_{l=1}^d\kappa_l\sum_{j\in\Jscr} (x_j^2 + y_j^2)
    -\sum_{l=1}^d\kappa_l\sum_{j\in\Jscr} (x_{j+l}x_j+ y_{j+l}y_j)\ .
  \end{aligned}
\end{equation}
Since $S$ corresponds to the set of the excited sites, we 
introduce the following variables
\begin{align*}
x_j &= \sqrt{2I_j}\cos\theta_j\ ,& y_j &= -\sqrt{2I_j}\sin\theta_j\ ,& j&\in S\ ,\\
x_j &= \frac{1}{\sqrt{2}}(\xi_j + \im \eta_j)\ ,&  y_j &= \frac{\im}{\sqrt{2}}(\xi_j - \im \eta_j) \ ,& 
j&\in \Jscr\setminus S\ .\\
\end{align*}
Thus the Hamiltonian~\eqref{e.ex.dnls} takes the form
\begin{equation}
H(I,\theta,\xi,\eta,\epsilon) = h_0(I)+g_0 (\xi,\eta) +\epsilon
    H_1(I,\theta,\xi,\eta;\epsilon)\ ,
\label{frm:H-modello}
\end{equation}
with
$$
h_0(I) = \sum_{j\in S}\left(I_j+\frac\gamma2 I_j^2
    \right)\ ,\qquad g_0(\xi,\eta) = \sum_{j\in\Jscr\setminus
      S}\left(\im \xi_j\eta_j
      -\frac\gamma2\xi_j^2\eta_j^2\right)\ ,
$$
and $H_1$ that is given by the expression in~\eqref{e.ex.dnls} written in the 
new variables.

\subsection{Normal form algorithm ($1:\text{\textendash}:1$ resonance)}

According to the geometrical interpretation given in the Introduction,
all the unperturbed periodic orbits foliate a $m$-dimensional torus
$\toro^m$ of the phase space: the torus corresponds to ${I}_j=I^*=R^2$
for $j\in S$ and $\xi_j=\eta_j=0$ for the remaining $j\in
\Jscr\setminus S$. Any orbit on such a torus is uniquely identified by
a point in the quotient space $\toro^{m-1}=\toro^m /\toro$; such a
point can be well represented by introducing a set of $m-1$ new
\emph{phase shifts} angles
\begin{equation}
\label{e.phidef}
q_j := \theta_{j_{l+1}}-\theta_{j_l}\  ,\qquad\qquad
l=1,\ldots,m\ .
\end{equation}
The definition of these new angles is related to the
$1:\text{\textendash}:1$ resonance among the unperturbed oscillators, thus we
will refer to them as \emph{resonant variables}.

In order to reveal the structure of the dynamics around the
unperturbed low-dimensional torus, we locally expand the Hamiltonian
in a neighbourhood of it.  Specifically, we expand $H$ in power series
of $J_j=I_j-I^*$ and introduce the resonant angles $\hat q=(q_1,q)$,
with $q_1=\theta_1$ and and $q$ as in~\eqref{e.phidef}, and we
complete canonically the transformation with the corresponding actions
$\hat p=(p_1,p)$; in
particular it follows that $p_1=\sum_{l\in S} J_l$. We finally split
the Hamiltonian~\eqref{frm:H-modello} as
\begin{equation}
  \label{e.H.nf.0}
\begin{aligned}
H^{(0)} &=\omega p_1 + \sum_{j\in\Jscr\setminus
  S}\im\xi_j\eta_j + f_4^{(0,0)} + \sum_{\ell>4}
f_\ell^{(0,0)}\\ &\qquad+f_0^{(0,1)}+ f_1^{(0,1)}+ f_2^{(0,1)}+ f_3^{(0,1)}
+ f_4^{(0,1)} + \sum_{\ell>4} f_\ell^{(0,1)} + \Oscr(\epsilon^2)\ ,
\end{aligned}
\end{equation}
where $\omega(I^*)=1+\gamma I^*$ is the frequency of any periodic
orbit on the unperturbed torus $\hat p=0$ and $f^{(r,s)}_\ell$ is a
polynomial of degree $l$ in ${\hat{p}}$ and degree $m$ in $(\xi,\eta)$
satisfying $\ell=2l+m$ and with coefficients depending on the angles
$\hat{q}$.  The index $r$ identifies the order of normalisation
($r=0$ corresponding to the original Hamiltonian), while $s$ keeps track of the
order in $\epsilon$. 

The standard approach to continue the periodic orbit surviving the
breaking of the unperturbed lower dimensional torus $I_j=I^*$ consists
in averaging the leading term of the perturbation, namely
$f_0^{(0,1)}(q_1,q)$, with respect to the fast angle $q_1$ and to look
for critical points of the averaged function on the torus
$\toro^{m-1}$. The explicit form of $f_0^{(0,1)}$ depends on the
choice of the set $S$ and of the coupling $H_1$, but it always
consists of trigonometric terms of the form $\cos(k\cdot q)$; hence,
solutions of $\nabla_{q} f_0^{(0,1)}=0$ always include
$q_l\in\{0,\pi\}$, but additional solutions, the so-called
\emph{phase-shift} solutions, might appear.  If the critical points
are not degenerate, continuation easily follows from an implicit
function theorem argument.  Instead, for the degenerate ones, like
$d$-parameter families with $d\geq 1$, it is necessary to take into
account higher order terms in the perturbation.

To this end in \cite{SanDPP20} we implement a normal form construction
for elliptic low dimensional and completely resonant tori that is
reminiscent of the Kolmogorov algorithm, see
also~\cite{SanLocGio-2010,GioLocSan-2014,SanLib2019}. Shortly, we perform a
sequence of canonical transformations which remove $f_1^{(0,1)}$ and
the part of $f_3^{(0,1)}$ which is linear in the transversal variables
$(\xi,\eta)$, and we average over the fast angle $q_1$ also the terms
$f_2^{(0,1)}$ and the part of $f_4^{(0,1)}$ which is quadratic in
$\hat p$ (see \cite{CarL21} for a strictly related construction
applied to the FPU model). First and second order nonresonance
conditions between $\omega$ and the linear frequencies $\Omega_j$
   \begin{align} 
     k_1\omega\pm \Omega_j &\neq 0\ ,\quad  k_1\in\interi \ ,
    \label{melnikov1}
    \\
     k_1\omega\pm \Omega_l \pm \Omega_k  &\neq 0\ ,\quad  k_1\in\interi\setminus \lbrace 0\rbrace \ ,
    \label{melnikov2}
   \end{align}
are needed to ensure the existence of such transformations; these are
the so-called first and second Melnikov conditions. In the dNLS
model~\eqref{frm:H-modello} here considered we have
$\Omega_l=\Omega_k=1$ for all $k,l\in\Jscr\setminus S$,
hence~\eqref{melnikov2} is turned into its simplified form
$k_1\omega\pm 2\neq 0$.  In addition, we perform a translation of the
actions $\hat{p}$ so as to keep fixed the linear frequency $\omega$;
here anharmonicity of $h_0(I)$ is relevant, which corresponds to the
so-called twist condition for $f_4^{(0,0)}(\hat p)$: there exists
$m>0$ such that
\begin{equation}
  m \sum_{i=1}^n |v_i| \leq \sum_{i=1}^n |\sum_{j=1}^n C_{ij}
  v_j|\ ,\quad\forall v\in\reali^{n}\quad\hbox{where}\quad
  C_{ij}=\frac{\partial^2 f_4^{(0,0)}}{\partial \hat{p}_i \partial
    \hat{p}_j}
  \label{frm:twist}\ .
\end{equation}
In this way the Hamiltonian is brought in normal form at order $r=1$.
Iterating $r$-times the same procedure, we get the
Hamiltonian in normal form at order $r\geq 2$, $H^{(r)} = K^{(r)} + \Oscr(\epsilon^{r+1})$, with
$$
K^{(r)} = \omega p_1 + \sum_{j\in\Jscr\setminus S}{\im\xi_j\eta_j} +
f_4^{(r,0)} + \sum_{\ell>4} f_\ell^{(r,0)}+ Z_0^{(r)} + Z_2^{(r)} + Z_3^{(r)} + Z_4^{(r)} +
\sum_{s=1}^r\sum_{\ell>4} f_\ell^{(r,s)}\ ,
$$
where
\begin{displaymath}
  Z_\ell^{(r)} = \sum_{s=1}^r f_\ell^{(r,s)}\ , \qquad \ell=0,2,3,4\ .
\end{displaymath}
A key ingredient in our construction is that the translation which
keeps fixed the frequency $\omega$ of the periodic orbit depends on a
parameter vector $q^*\in\toro^{m-1}$; in particular, the translation
of $\hat p$ is such that $Z_2^{(r)}(q^*,\hat p,0;q^*)= 0$.

At leading order, periodic orbits of the form
\begin{equation}
\label{e.appr.sol}
\dot q_1 = \omega\ ,\qquad q=q^*\ ,\qquad {\hat p}=\xi=\eta=0\ ,
\end{equation}
correspond to relative equilibria of the truncated normal form
$K^{(r)}$, provided $q^*$ satisfies
\begin{equation}
\label{e.ex.qstar}
{\nabla_q Z^{(r)}_0}(q;q^*)\Big|_{q=q^*}=0\ .
\end{equation}
Then, continuation of the approximate periodic
orbit~\eqref{e.appr.sol} could follow by means of a fixed point
method, once suitable spectral conditions are verified.  More
precisely, we introduce the smooth map
$\Upsilon(x):\Uscr(x^*)\subset\reali^{2n-1}\to
\Vscr(x^*)\subset\reali^{2n-1}$ as
\begin{equation}\label{frm:Ups}
\Upsilon(x(0);\epsilon,q_{1}(0))=
  \begin{pmatrix}
  q_1(T)-q_1(0)-\omega T\\
  q(T)-q(0)\\
  p(T)-p(0)\\
  \xi(T)-\xi(0) \\
  \eta(T)-\eta(0)
  \end{pmatrix} \ ,
\end{equation}
parameterised by the initial phase $q_1(0)$ and $\epsilon$, with $T$
the period of the periodic orbit; the map $\Upsilon$ is basically the
variation over the period $T$ of the Hamiltonian flow (a part from the
coordinate $p_1$). The main result (proved in \cite{SanDPP20})
used in the examples is the following

\begin{theorem}
\label{teo:forma-normale-r}
  Consider the map $\Upsilon$ defined in~\eqref{frm:Ups} in a
  neighbourhood of the lower dimensional torus $\hat{p}=0$, $\xi=\eta=0$ and let
  $x^*(\epsilon)=({q^*}(\epsilon),0,0,0)$, with ${q^*}(\epsilon)$
  satisfying~\eqref{e.ex.qstar}.  Assume that
  \begin{equation}
    \label{e.small.Ups}
    \left|\Upsilon(x^*(\epsilon);\epsilon,q_1(0))\right| \leq C_1
    \epsilon^{r+1}\ ,
  \end{equation}
  where $C_1$ is a positive constant depending on $\Uscr$ and
  $r$. Assume also that
  $M(\epsilon):=\Upsilon'(x^*(\epsilon);\epsilon,q_1(0))$ is invertible
  and there exists $\alpha>0$ with $2\alpha<r+1$ such that
  \begin{equation}
    \label{e.small.eig}
    |\lambda|\gtrsim |\epsilon|^\alpha\ ,\qquad\hbox{for all}\quad\lambda\in
    \Sigma\bigl(M(\epsilon)\bigr)\ .
  \end{equation}
  Then, there exist $C_0>0$ and $\epsilon^*>0$ such that for any
  $0\leq|\epsilon|<\epsilon^*$ there exists a unique $x^*_{\rm
    p.o.}(\epsilon)=(q^*_{\rm p.o.}(\epsilon),\hat{p}_{\rm
    p.o.}(\epsilon),\xi_{\rm p.o.}(\epsilon),\eta_{\rm
    p.o.}(\epsilon))\in\Uscr$ which solves
  \begin{equation}
    \label{e.exist.appr}
    \Upsilon(x^*_{\rm p.o.};\epsilon,q_1(0))=0
    \ ,\qquad\hbox{with}\quad |x^*_{\rm p.o.}-x^*| \leq
    C_0\epsilon^{r+1-\alpha}\ .
\end{equation}
Moreover, the approximate linear stability of $x^*_{\rm p.o.}$ is
encoded in the linear stability of the relative equilibrium $x^*$,
hence in the Floquet multipliers of the matrix $\exp\bigl(L(\epsilon)T\bigr)$,
where $L:=DX_{K^{(r)}}(x^*)$ has the block diagonal form
\begin{equation}
L(\epsilon)=\left( \begin{matrix}
{L_{11}}(\epsilon) & O\\
O & L_{22}(\epsilon) \\
\end{matrix}\right)\ .
\label{frm:L11L22}
\end{equation}
\end{theorem}

The validity of~\eqref{e.small.eig} is checked by exploiting the
definition of $L$. Indeed, as explained in \cite{SanDPP20}, the matrix
$M(\epsilon)$ can be obtained by removing the row and the column
corresponding to $p_1$ and $q_1$ from $\Phi-\textrm{Id}$, where $\Phi$
is the monodromy matrix, which is well approximated by
$\exp(LT)$. Hence the scaling in $\epsilon$ of the smallest
eigenvalues of $M(\epsilon)$ can be extract from $L$. At the same
time, $L$ provide the approximate linear stability of the periodic
orbit. Indeed, thanks to the block-diagonal structure, its spectrum
splits into two different components: $\Sigma(L_{22})\subset\im\reali$
since the quadratic part $K_2^{(r)}(\xi,\eta)$ is positive definite
for $\epsilon$ small enough (by continuity at $\epsilon=0$), while
$\Sigma(L_{11})$ is generically made of $m-1$ pairs of eigenvalues
$\lambda_j(\epsilon)\to 0$ as $\epsilon\to 0$ (and of a couple of zero
eigenvalues). Hence approximate linear stability depends only on the
\emph{internal} Floquet-exponents $\Sigma(L_{11})$.  The effective
linear stability of the periodic orbit can be derived from the
approximate spectrum if the approximate Floquet multipliers are well
distinct and the perturbation is sufficiently small, see Theorem~2.3
in~\cite{SanDPP20}, that we report below for completeness.
\begin{theorem}
\label{t.lin.stab.2}
Assume that $L_{11}(\epsilon)$ has $2m-2$ distinct nonzero eigenvalues
and let $\tilde c>0$ and $\beta<r+1-\alpha$, with $2\alpha<r+1$ as in
Theorem~\ref{teo:forma-normale-r}, be such that
\begin{equation}
\label{e.dist.eig}
|\lambda_j-\lambda_k|>\tilde c\epsilon^{\beta}\ ,\qquad\hbox{for
  all}\quad\lambda_{j},\lambda_{k}\in\Sigma(L_{11}(\epsilon))\setminus\{0\}\ .
\end{equation}
Then there exists $\epsilon^*>0$ such that if $|\epsilon|<\epsilon^*$
and $\mu=e^{\lambda T}\in\Sigma(\exp(L_{11}(\epsilon) T))$, there
exists one Floquet multiplier $\nu$ of $x^*_{p.o.}$ inside the complex
disk $D_\epsilon(\mu)=\{z\in\complessi\ :\ |z-\mu|<c
  \epsilon^{r-\alpha+1}\}$, with $c>0$ a suitable constant independent
of $\mu$.
\end{theorem}

\section{Applications to dNLS models}
\label{s:2}

In the present section we show how the results already illustrated in
the Introduction have been obtained exploiting the normal form
construction with the aid of a computer algebra system. We will mainly
restrict to the so-called \emph{focusing} case in~\eqref{e.KdNLS},
i.e., with $\gamma=1$, which implies $f_4^{(0,0)}$ positive definite
and the actions $I_j$ of the torus are globally defined on the phase
space. We choose values of $I^*$ such that~\eqref{melnikov1}
and~\eqref{melnikov2} are satisfied up to the required normalisation
order. Furthermore, we will usually take $\Jscr=\{-N,\ldots,N\}$ with
$N=\leq 10$ and $S$ involving at most $4$ sites: this choice allows to
explore meaningful configurations with at most $r=3$ normal form
steps, thus keeping the presentation more compact and easy to follow.

Actually, the results here presented have been obtained via an
implementation of the normal form algorithm in Mathematica that can be
found at \url{https://github.com/marcosansottera/periodic_orbits_NF}.

\subsection{Multi-pulse solutions in the standard dNLS model}

We start with the well known standard dNLS model, where only
$\kappa_1=1$ in~\eqref{e.KdNLS}, namely only nearest-neighbours
interactions are active.  We are going to consider two different kind
of sets $S$, both dealing with problem of degeneracy due to
nonconsecutive excited sites. In the first case we take only two
nonconsecutive sites $S=\{-l,l\}$, with $l\geq 1$, the larger is the
distance $2l$ among the sites, the greater is the number of normal
form steps needed to remove the degeneracy, i.e., $r=2l$. In the
second case we take 3 sites, giving an asymmetric configuration
$S=\{-2,-1,1\}$.  This is the easiest asymmetric example which
exhibits degeneracy, due to the lack of the interaction at order
$\epsilon$ between the second and the fourth site. In agreement with
the existing literature (see, e.g., \cite{HenT99,PelKF05,Kev_book09}),
it will be shown that only standard in/out-of-phase solutions do
exist. Linear stability analysis provides a scaling of the Floquet
exponents coherent with the literature and Theorem~\ref{t.lin.stab.2}
can be always applied in our examples. In addition, the normal form
remarkably shows the effect of switching from focusing to defocusing
dNLS, obtained by changing the sign of $\gamma$: nondegenerate saddle
and centre eigenspaces exchange their stability, while degenerate ones
keep unchanged their stability whenever the order of degeneracy is
even, as with $S=\{-2,-1,1\}$.

\bigskip
\noindent
\textbf{Two-sites multi-pulse discrete solitons for dNLS $S=\{-1,1\}$}

\bigskip
\begin{center}
  \begin{tikzpicture}[scale = 0.8]
    \begin{scope} 
      \foreach \name/\posX/\posY in {A/0/0,B/1/0,C/2/0,D/3/0,E/4/0,F/5/0,G/6/2,H/7/0,I/8/2,L/9/0,M/10/0,N/11/0,O/12/0,P/13/0,Q/14/0}
      \node (\name) at (\posX,\posY) {};
      \draw [dotted, black!50] (-1,0) -- (14,0);
      \draw [dotted, black!50] (-1,2) -- (14,2);
      \draw [dotted, black!50] (6,2) -- (6,0);
      \draw [dotted, black!50] (8,2) -- (8,0);
      \draw (B) circle (2mm) \foreach \n in {C,D,E,F,G,H,I,L,M,N,O,P} {(\n) circle (2mm)};
      \fill [black] \foreach \n in {G,I} {(\n) circle (2.0mm)};
      \fill [white] \foreach \n in {B,C,D,E,F,H,L,M,N,O,P} {(\n) circle (1.9mm)};
      \draw [red] (A) -- (B) -- (C) -- (D) -- (E) -- (F) -- (G)-- (H) -- (I) -- (L) -- (M) -- (N) -- (O) -- (P) -- (Q) ;
      \node[label=below:{\footnotesize-6}] at (B) {\phantom{$\times$}};
      \node[label=below:{\footnotesize-5}] at (C) {\phantom{$\times$}};
      \node[label=below:{\footnotesize-4}] at (D) {\phantom{$\times$}};
      \node[label=below:{\footnotesize-3}] at (E) {\phantom{$\times$}};
      \node[label=below:{\footnotesize-2}] at (F) {\phantom{$\times$}};
      \node[label=below:{\footnotesize 0}] at (H) {\phantom{$\times$}};
      \node[label=below:{\footnotesize 2}] at (L) {\phantom{$\times$}};
      \node[label=below:{\footnotesize 3}] at (M) {\phantom{$\times$}};
      \node[label=below:{\footnotesize 4}] at (N) {\phantom{$\times$}};
      \node[label=below:{\footnotesize 5}] at (O) {\phantom{$\times$}};
      \node[label=below:{\footnotesize 6}] at (P) {\phantom{$\times$}};
      \node[label=below:{\footnotesize-1}] at (6,0) {\phantom{$\times$}};
      \node[label=below:{\footnotesize 1}] at (8,0) {\phantom{$\times$}};
      \node[label=left:{\footnotesize $0$}] at (-1,0) {};
      \node[label=left:{\footnotesize $I^*$}] at (-1,2) {};
    \end{scope}
  \end{tikzpicture}
\end{center}
In the first case, the perturbation $H_1$, given by the nearest
neighbours interactions, reads
\begin{displaymath}
H_1=2\sum_{j\in S}I_j + 2\sum_{j\in\Jscr\setminus S}\im\xi_j\eta_j-\sum_{j\in\Jscr}(x_{j+1}x_j+y_{j+1}y_j)
\end{displaymath}
where the products $x_{j+1}x_j+y_{j+1}y_j$ are of the following
types
$$
\begin{aligned}
x_{j+1}x_j+y_{j+1}y_j &= \im(\xi_{j+1}\eta_j+\xi_j\eta_{j+1})\ ,& \hbox{if } &j\hbox{ and }j+1\not\in S\ , \\
x_{j+1}x_j+y_{j+1}y_j &= \sqrt{I_j}\left(\cos(\theta_j)(\xi_{j+1}+\im\eta_{j+1})  -\im\sin(\theta_j)(\xi_{j+1}-\im\eta_{j+1})\right)\ , & \hbox{if } &j\in S\ , \\
x_{j+1}x_j+y_{j+1}y_j &= \sqrt{I_{j+1}}\left(\cos(\theta_{j+1})(\xi_{j}+\im\eta_{j})  -\im\sin(\theta_{j+1})(\xi_{j}-\im\eta_{j})\right)\ , & \hbox{if } &j+1\in S\ , 
\end{aligned}
$$
thus no term of the form $\cos(\theta_l-\theta_{-l})$ appears at
order $\Oscr(\epsilon)$. Expanding $H_0$ and $H_1$ in Taylor
series of the actions around $I^*$, forgetting constant terms and
introducing the resonant angles $\hat{q}=(q_1,q)$ and their conjugated
actions $\hat{p}=(p_1,p)$, i.e.,
$$
 \left\lbrace \begin{aligned}
& q_1 = \theta_{-l} \ ,\\
& q_2 = \theta_{l} - \theta_{-l} \ ,
\end{aligned}
\right.
\qquad
 \left\lbrace \begin{aligned}
& p_1 = J_l + J_{-l} \ ,\\
& p_2 = J_l \ ,
\end{aligned}
 \right. 
$$ 
we can rewrite the initial Hamiltonian as
$$
\begin{aligned}
H^{(0)} & = \omega p_1 + \sum_{j\in\Jscr\setminus S}{\im\xi_j\eta_j} +
f_4^{(0,0)}(\hat p, \xi,\eta) + f_1^{(0,1)}(\hat q,\xi,\eta) \\ &\qquad +
f_2^{(0,1)}(\hat p,\xi,\eta) + f_3^{(0,1)}(\hat q,\hat
p,\xi,\eta) + \sum_{\ell\geq
  5}f_\ell^{(0,1)}(\hat q,\hat p,\xi,\eta) \ ,
\end{aligned} 
$$ with $\omega=1+I^*$. Notice that $f_0^{(0,1)}$ and
$f_4^{(0,1)}$ are missing and that $f_2^{(0,1)}$ does not depend on
$\hat q$: this is due to the lack of coupling terms
$x_jx_{j+1}+y_jy_{j+1}$ with both $j$ and $j+1$ belonging to $S$.

For $l=1$, $S=\{-1,1\}$ and the normal form
at order one gives $Z_0^{(1)}=f_0^{(1,1)}(q)\equiv 0$, hence any $q_2\in\toro$ is a
critical point and the problem is trivially degenerate.

The normal form at order two gives $Z_0^{(2)}=f_0^{(2,2)}=2\epsilon^2 \cos(q_2)$, thus 
\begin{displaymath}
\nabla_qZ_0^{(2)} = -2\epsilon^2\sin(q_2)=0\ ,
\end{displaymath}
provides only the standard solutions $q_2\in\{0,\pi\}$. In order
to conclude the existence of these two in/out-of-phase configurations,
we need to check condition~\eqref{e.small.eig} with $\alpha<\frac32$,
explicit symbolic computations performed with Mathematica gives
$\alpha=1$. The stability analysis shows that $q_2=0$ (the so-called
Page mode) is the unstable configuration, while $q_2=\pi$ (the
so-called Twist mode) is the stable one, with approximated Floquet
exponents
$$
\lambda_0= \pm \left(2\epsilon-\frac{\epsilon^3}{(I^*)^2}\right)\ ,\quad
\lambda_\pi= \pm \im\left(2\epsilon+\frac{\epsilon^3}{(I^*)^2}\right)\ .
$$
Theorem~\ref{t.lin.stab.2} applies with
$\beta=1<2=r+1-\alpha$, hence Floquet multipliers are $\epsilon^2$-close
to the approximate ones $e^{\lambda T}$ (where $T$ is the period).

For $l>1$, the procedure for the continuation is clearly the same: it
turns out that degeneracy persists up to order $r=2l-1$, namely
$f_0^{(s,1)}\equiv 0$ for $s\leq 2l-1$. At order $r=2l$ one has
$Z_0^{(r)}=f_0^{(r,r)}$ such that
\begin{displaymath}
\nabla_q Z_0^{(r)}=c(I^*)\epsilon^{r}\sin(q_2)=0\ ,
\end{displaymath}
with $c(I^*)$ a constant depending on $I^*$, which again provides only
standard solutions $q_2=\{0,\pi\}$. Existence of these two
in/out-of-phase configurations is ensured by~\eqref{e.small.eig} with $\alpha=r/2<(r+1)/2$.  Stable and
unstable configurations are expected to be respectively $q_2=\pi$ and
$q_2=0$, with approximate Floquet exponents of order
$\Oscr(\epsilon^{l})$.

\bigskip\pagebreak
\noindent
\textbf{Three-sites multi-pulse discrete solitons for dNLS  $S=\{-2,-1,1\}$}

\begin{center}
  \begin{tikzpicture}[scale = 0.8]
    \begin{scope} 
      \foreach \name/\posX/\posY in {A/0/0,B/1/0,C/2/0,D/3/0,E/4/0,F/5/2,G/6/2,H/7/0,I/8/2,L/9/0,M/10/0,N/11/0,O/12/0,P/13/0,Q/14/0}
      \node (\name) at (\posX,\posY) {};
      \draw [dotted, black!50] (-1,0) -- (14,0);
      \draw [dotted, black!50] (-1,2) -- (14,2);
      \draw [dotted, black!50] (5,2) -- (5,0);
      \draw [dotted, black!50] (6,2) -- (6,0);
      \draw [dotted, black!50] (8,2) -- (8,0);
      \draw (B) circle (2mm) \foreach \n in {C,D,E,F,G,H,I,L,M,N,O,P} {(\n) circle (2mm)};
      \fill [black] \foreach \n in {F,G,I} {(\n) circle (2.0mm)};
      \fill [white] \foreach \n in {B,C,D,E,H,L,M,N,O,P} {(\n) circle (1.9mm)};
      \draw [red] (A) -- (B) -- (C) -- (D) -- (E) -- (F) -- (G)-- (H) -- (I) -- (L) -- (M) -- (N) -- (O) -- (P) -- (Q) ;
      \node[label=below:{\footnotesize-6}] at (B) {\phantom{$\times$}};
      \node[label=below:{\footnotesize-5}] at (C) {\phantom{$\times$}};
      \node[label=below:{\footnotesize-4}] at (D) {\phantom{$\times$}};
      \node[label=below:{\footnotesize-3}] at (E) {\phantom{$\times$}};
      \node[label=below:{\footnotesize 0}] at (H) {\phantom{$\times$}};
      \node[label=below:{\footnotesize 2}] at (L) {\phantom{$\times$}};
      \node[label=below:{\footnotesize 3}] at (M) {\phantom{$\times$}};
      \node[label=below:{\footnotesize 4}] at (N) {\phantom{$\times$}};
      \node[label=below:{\footnotesize 5}] at (O) {\phantom{$\times$}};
      \node[label=below:{\footnotesize 6}] at (P) {\phantom{$\times$}};
      \node[label=below:{\footnotesize-2}] at (5,0) {\phantom{$\times$}};
      \node[label=below:{\footnotesize-1}] at (6,0) {\phantom{$\times$}};
      \node[label=below:{\footnotesize 1}] at (8,0) {\phantom{$\times$}};
      \node[label=left:{\footnotesize $0$}] at (-1,0) {};
      \node[label=left:{\footnotesize $I^*$}] at (-1,2) {};
    \end{scope}
  \end{tikzpicture}
\end{center}
The perturbation $H_1$, given by the nearest neighbours interactions,
reads
\begin{displaymath}
H_1=2\sum_{j\in S}I_j + 2\sum_{j\in\Jscr\setminus S}\im\xi_j\eta_j-\sum_{j\in\Jscr}(x_{j+1}x_j+y_{j+1}y_j)
\end{displaymath}
where the products $x_{j+1}x_j+y_{j+1}y_j$ are of the following
types
$$
\begin{aligned}
x_{j+1}x_j+y_{j+1}y_j &=
\im(\xi_{j+1}\eta_j+\xi_j\eta_{j+1})\ ,& \hbox{if } &j\hbox{ and }j+1\not\in S\ , \\
x_{j+1}x_j+y_{j+1}y_j &=
\sqrt{I_j} \bigl(\cos(\theta_{j})(\xi_{j+1}+\im\eta_{j+1})
  -\im\sin(\theta_j)(\xi_{j+1}-\im\eta_{j+1})\bigr)\ , & \hbox{if } &j\in S\,,\ j+1\not\in S\ , \\
x_{j+1}x_j+y_{j+1}y_j &=
\sqrt{I_{j+1}}\bigl(\cos(\theta_{j+1})(\xi_{j}+\im\eta_{j})
  -\im\sin(\theta_{j+1})(\xi_{j}-\im\eta_{j})\bigr)\ , & \hbox{if } &j\not\in S\,,\ j+1\in S\ , \\
x_{-2}x_{-1}+y_{-2}y_{-1} &=2\sqrt{I_{-2}I_{-1}}\cos(\theta_{-2}-\theta_{-1})
\end{aligned}
$$
Expanding $H_0$ and $H_1$ in Taylor series of the actions
around $I^*$, forgetting constant terms and introducing the resonant
angles $\hat{q}=(q_1,q)$ and their conjugated actions
$\hat{p}=(p_1,p)$
$$
 \left\lbrace \begin{aligned}
& q_1 = \theta_{-2} \ ,\\
& q_2 = \theta_{-1} - \theta_{-2} \ ,\\ 
& q_3 = \theta_1 - \theta_{-1} \ ,\\
\end{aligned}
\right. 
\qquad
 \left\lbrace \begin{aligned}
& p_1 = J_{-2} + J_{-1} + J_1 \ ,\\
& p_2 = J_{-1} + J_1 \ ,\\
& p_3= J_1 \ ,\\
\end{aligned}
\right.
$$
we can rewrite the initial Hamiltonian as
$$
\begin{aligned}
H^{(0)} &= \omega p_1 + \sum_{j\in\Jscr\setminus S}{\im\xi_j\eta_j} +
f_4^{(0,0)}(\hat p, \xi,\eta) + f_0^{(0,1)}(q_2) + f_1^{(0,1)}(\hat
q,\xi,\eta) \\ &\qquad+ f_2^{(0,1)}(q,\hat p,\xi,\eta)+f_3^{(0,1)}(\hat
q,\hat p,\xi,\eta)+ f_4^{(0,1)}(q,\hat p)+ \sum_{\ell\geq
  5}f_\ell^{(0,1)}(\hat q,\hat p,\xi,\eta) \ ,
\end{aligned} 
$$ with $\omega = (1+I^*)$. The normal form at order one gives
$Z_0^{(1)}=f_0^{(1,1)}=-2\epsilon I^*\cos(q_2)$, thus the critical points $q^*$ of
 on the torus $\toro^2$ are two disjoint one-parameter families $P_1(\vartheta)=(0,\vartheta)$ and
 $P_2(\vartheta)=(\pi,\vartheta)$, where $\vartheta=q_3$.

The normal form at order two gives 
$Z_0^{(2)}=Z_0^{(1)} + f_0^{(2,2)}$ with
$$
f_0^{(2,2)} = \epsilon^2 \bigl(4\cos(q_2) + 2\cos(q_3) -\frac12\cos(2q_2)\bigr)\ ,
$$
thus the critical points are the four in/out-of-phase
solutions $(q_2,q_3)\in\{(0,0),(0,\pi),(\pi,0),(\pi,\pi)\}$. In order
to conclude the existence of these configurations, we need to check
condition~\eqref{e.small.eig} with $\alpha<\frac32$, explicit symbolic computations with
Mathematica gives $\alpha=1$. Linear stability analysis provides $(0,\pi)$ as the only
stable configurations with approximate Floquet exponents
\begin{displaymath}
  \lambda_{(0,\pi)}=
\begin{cases}
\pm\im\left(2\sqrt{I^*\epsilon}+\dfrac{\epsilon^{3/2}}{4\sqrt{I^*}}+\Oscr(\epsilon^{3/2})\right) \ ,\\
\pm\im\sqrt3\left(\epsilon -\dfrac{\epsilon^2}{8I^*} + \Oscr(\epsilon^{2})\right) \ ,
\end{cases}
\end{displaymath}
while the other three configurations are all unstable, with
\begin{displaymath}
\lambda_{(0,0)}=
\begin{cases}
\pm2\im\sqrt{I^*\epsilon}+\Oscr(\sqrt\epsilon)\ ,\\
\pm\sqrt3\epsilon+\Oscr(\epsilon)\ ,
\end{cases}
\ 
\lambda_{(\pi,0)}=
\begin{cases}
\pm2\sqrt{I^*\epsilon}+\Oscr(\sqrt\epsilon)\ ,\\
\pm\sqrt3\epsilon+\Oscr(\epsilon)\ ,
\end{cases}
\ 
\lambda_{(\pi,\pi)}=
\begin{cases}
\pm2\sqrt{I^*\epsilon}+\Oscr(\sqrt\epsilon)\ ,\\
\pm\im\sqrt3\epsilon+\Oscr(\epsilon)\ .
\end{cases}\ .
\end{displaymath}
Approximate linear stability corresponds to effective linear
stability, since Theorem~\ref{t.lin.stab.2} applies with
$\beta=1<2=r+1-\alpha$, hence Floquet multipliers are located
$\epsilon^2$-close to the approximate ones and fulfil the usual
symmetries of the spectrum of a symplectic matrix.

\begin{remark}
It is interesting to investigate what happens to the Floquet exponents once the sign of
the nonlinear coefficient $\gamma$ is changed. It turns out, as
already stressed in the literature, that eigenvalues of order
$\Oscr(\sqrt\epsilon)$ switch from real to imaginary and vice versa, hence
stable and unstable eigenspaces are exchanged. However, eigenvalues of
order $\Oscr(\epsilon)$ keep their nature. This is the effect of a
cancellation of $\gamma$ in front of the equation $-2\epsilon^2\sin(q_3)$,
as already stressed in the ``seagull'' example
in \cite{SanDPP20}. Hence the new stable configuration would be in
this case $(\pi,\pi)$.
\end{remark}

\bigskip
\subsection{Four-sites vortex-like solutions in a Zigzag dNLS cell}
\begin{center}
  \begin{tikzpicture}[scale = 0.8]
    \begin{scope} 
      \foreach \name/\posX/\posY in {C/3/0,D/4.5/0,E/6/0,F/7.5/0,G/9/1.,H/10.5/1.,I/12/0,L/13.5/0,M/15/0,N/16.5/0,
        T/3.5/1.5,U/5./1.5,V/6.5/1.5,W/8./1.5,Y/9.5/2.5,X/11./2.5,Z/12.5/1.5,J/14/1.5,K/15.5/1.5,a/17./1.5}
      \node (\name) at (\posX,\posY) {};
      \draw (D) circle (2mm) \foreach \n in {D,E,F,G,H,I,L,M,U,V,W,Y,X,Z,J,K} {(\n) circle (2mm)};
      \fill [black] \foreach \n in {G,H,Y,X} {(\n) circle (2.0mm)};
      \draw [red]  (C) -- (D) -- (E) -- (F) -- (G)-- (H) -- (I) -- (L) -- (M) -- (N)  ;
      \draw [red]  (T) -- (U) -- (V) -- (W) -- (Y)-- (X) -- (Z) -- (J) -- (K) -- (a)  ;
      \draw [red]   (D) --(U) (E) -- (V) (F) -- (W) (G)-- (Y) (H) -- (X) (I) -- (Z) (L) -- (J) (M) -- (K) ;
      \draw [red]    (U) -- (E)  (V) -- (F) (W) -- (G) (Y) -- (H) (X) -- (I) (Z) -- (L) (J) -- (M)   ;
      \node[label=below:{\footnotesize-6}] at (D) {\phantom{$\times$}};
      \node[label=below:{\footnotesize-4}] at (E) {\phantom{$\times$}};
      \node[label=below:{\footnotesize-2}] at (F) {\phantom{$\times$}};
      \node[label=below:{\footnotesize 0}] at (G) {\phantom{$\times$}};
      \node[label=below:{\footnotesize 2}] at (H) {\phantom{$\times$}};
      \node[label=below:{\footnotesize 4}] at (I) {\phantom{$\times$}};
      \node[label=below:{\footnotesize 6}] at (L) {\phantom{$\times$}};
      \node[label=below:{\footnotesize 8}] at (M) {\phantom{$\times$}};
      \node[label=above:{\footnotesize-5}] at (U) {\phantom{$\times$}};
      \node[label=above:{\footnotesize-3}] at (V) {\phantom{$\times$}};
      \node[label=above:{\footnotesize-1}] at (W) {\phantom{$\times$}};
      \node[label=above:{\footnotesize 1}] at (Y) {\phantom{$\times$}};
      \node[label=above:{\footnotesize 3}] at (X) {\phantom{$\times$}};
      \node[label=above:{\footnotesize 5}] at (Z) {\phantom{$\times$}};
      \node[label=above:{\footnotesize 7}] at (J) {\phantom{$\times$}};
      \node[label=above:{\footnotesize 9}] at (K) {\phantom{$\times$}};
    \end{scope}
  \end{tikzpicture}
\end{center}
Let us consider the Hamiltonian system~\eqref{e.KdNLS} with
$\kappa_1=\kappa_2=1$, namely the so-called Zigzag model. This is a
particular case of two coupled one-dimensional dNLS models, where the
Zigzag coupling provides a one-dimensional Hamiltonian systems. We
want to investigate the continuation of vortex-like localised
structures given by four consecutive excited sites; hence the low
dimensional resonant torus is $I_l=I^*$ for $l\in S=\{0,1,2,3\}$
and $\xi_l=\eta_l=0$ for $l\in \Jscr\setminus S$. These configurations
have been the object of investigation of \cite{PenKSKP19} where
nonexistence of four-sites solutions with phase differences $q_l$
different from $\{0,\pi\}$ have been obtained with a Lyapunov-Schmidt
reduction.  We here show how to recover the same results via normal
form and we correct a minor statement on the nondegeneracy of the
isolated configurations.

Here the perturbation, involving nearest and
next-to-nearest neighbours interactions, reads
\begin{displaymath}
H_1=4\sum_{j\in S}I_j + 4\sum_{j\in\Jscr\setminus
  S}\im\xi_j\eta_j-2\sum_{j\in\Jscr} \bigl((x_{j+1}x_j+y_{j+1}y_j)+(x_{j+2}x_j+y_{j+2}y_j)\bigr)
\end{displaymath}
where, as in the previous examples, the products in the last sum must
be expressed in the $(I,\theta,\xi,\eta)$ variables.  Expanding
$H_0$ and $H_1$ in Taylor series of the actions around $I^*$,
forgetting constant terms and introducing the resonant angles
$\hat{q}=(q_1,q)$ and their conjugated actions $\hat{p}=(p_1,p)$
$$
 \left\lbrace \begin{aligned}
& q_1 = \theta_0 \ ,\\
& q_2 = \theta_1 - \theta_0 \ ,\\ 
& q_3 = \theta_2 - \theta_1 \ ,\\
& q_4 = \theta_3 - \theta_2 \ ,
\end{aligned}
\right. \ ,
\qquad
 \left\lbrace \begin{aligned}
& p_1 = J_0 + J_1 + J_2 + J_3 \ ,\\
& p_2 = J_1 + J_2 + J_3\ ,\\
& p_3= J_2 + J_3 \ ,\\
& p_4 = J_3\ .
\end{aligned}
\right. \ ,
$$
we can rewrite the initial Hamiltonian in the form
$$
\begin{aligned}
H^{(0)} &= \omega p_1 + \sum_{j\in\Jscr\setminus S}{\im\xi_j\eta_j} +
f_4^{(0,0)}(\hat p, \xi,\eta) + f_0^{(0,1)}(q) + f_1^{(0,1)}(\hat
q,\xi,\eta) \\ &\qquad+ f_2^{(0,1)}(q,\hat p,\xi,\eta)+f_3^{(0,1)}(\hat
q,\hat p,\xi,\eta)+ f_4^{(0,1)}(q,\hat p,\xi,\eta)+ \sum_{\ell\geq
  5}f_\ell^{(0,1)}(\hat q,\hat p,\xi,\eta) \ ,
\end{aligned} 
$$ with $\omega=(1+I^*)$.  The normal form at order one gives
$$
Z_0^{(1)} = f_0^{(1,1)}(q) = 2I^*\epsilon\bigl(\cos(q_2)+\cos(q_3)+\cos(q_4)+\cos(q_2+q_3)+\cos(q_3+q_4)\bigr)\ ,
$$ thus there are four isolated solutions $(0,0,0)$, $(0,0,\pi)$,
$(\pi,0,0)$, $(\pi,0,\pi)$, and two one-parameter families
$P_1(\vartheta)=(\vartheta,\pi,\vartheta-\pi)$ and
$P_2(\vartheta)=(\vartheta,\pi,-\vartheta)$.  In order to apply
Theorem~\ref{teo:forma-normale-r} with $r=1$, critical points need to
be not degenerate; by calculating the determinant in correspondence of
the $q^*$-values determined above, we see that nondegeneracy is
fulfilled only in three of the four isolated solutions $(0,0,0)$,
$(0,0,\pi)$, $(\pi,0,0)$, while the fourth isolated configuration
$(\pi,0,\pi)$ and the two families are degenerate. In particular, the
topologically isolated configuration $(\pi,0,\pi)$ is a degenerate
minimizer of $Z_0^{(1)}$, since along the tangent direction
$(\pi+t,-2t,\pi+t)$ it is possible to observe a growth as
$\Oscr(t^4)$; this represents an example of degenerate isolated
configuration.  In order to conclude the continuation of the three
nondegenerate in/out-of-phase configurations to effective periodic
orbits of the system, we need to check condition~\eqref{e.small.eig}
with $\alpha<1$, explicit symbolic computations with Mathematica give
$\alpha=1/2$.

For the degenerate configurations we have to compute the normal form at
order two. The equation for the critical points of
$Z_0^{(2)}$ takes the form
$$
F(q_2,q_3,q_4,\varepsilon)= F_0 (q_2,q_3,q_4)+ \varepsilon F_1 (q_2,q_3,q_4)=0,
$$ with $F:\mathbb{T}^3 \times \mathcal{U}(0) \rightarrow
\mathbb{R}^3$ (we here omit the explicit expression of $F_1$). We
observe that the vectors
$$
\partial_{\vartheta}P_1=\begin{pmatrix}
1 \\ 0 \\ 1
\end{pmatrix}
\qquad \text{and} \qquad
\partial_{\vartheta}P_2=\begin{pmatrix}
1 \\ 0 \\ -1
\end{pmatrix}
$$ generate the Kernel of $D_qF_0(q)\Big|_{q=P_j(\vartheta)}$ with
$j=1,2$.  The necessary condition to continue the degenerate solutions
of $F_0(q)=0$ to solutions of $F(q,\epsilon)=0$ is then
$F_1(P_j(\vartheta))\perp\partial_\vartheta P_j(\vartheta)$, i.e.,
$$
\langle F_1(P_1(\vartheta)), \partial_{\vartheta}P_1 \rangle = 
4 \sin(2\vartheta)\ ,\quad
\langle F_1(P_2(\vartheta)), \partial_{\vartheta}P_2 \rangle = 
-4\sin(\vartheta)\ ,
$$ and we can deduce that the two families $P_j(\vartheta)$ break down
and only four configurations $(0,\pi,\pi)$, $(\pi,\pi,0)$,
$(0,\pi,0)$, $(\pi,\pi,\pi)$ are solutions of $F(q,\epsilon)$.  Hence
the critical points of $Z_0^{(2)}$ are given only by in/out-of phase
configurations.  To prove the existence of these configurations we
need to check condition~\eqref{e.small.eig} with $\alpha<\frac32$,
explicit symbolic computations with Mathematica give $\alpha=1$.

Concerning the linear stability analysis, we summarise below the
results for the different cases.

\smallskip
\noindent
\textbf{Isolated and nondegenerate solutions.} The three
configurations $(0,0,0)\ ,(0,0,\pi)\ ,(\pi,0,0)$ have all approximate
Floquet exponents of order $\Oscr(\sqrt\epsilon)$ and can be computed
directly from the normal form at order one.  Specifically, $(0,0,0)$ is
the unique stable configuration, with
$$
\begin{aligned}
  \lambda_{1,2}(\epsilon)&=\pm 2\sqrt{2}\im\left(\sqrt{I^*}\sqrt\epsilon + \dfrac{\epsilon^{3/2}}{\sqrt {I^*}}+\Oscr(\epsilon^{3/2})\right) \ ,\\
  \lambda_{3,4}(\epsilon)&=\pm\im\left( 2\sqrt{2}\sqrt {I^*}\sqrt\epsilon + \dfrac{5 \epsilon^{3/2}}{\sqrt{2}\sqrt {I^*}}+\Oscr(\epsilon^{3/2}) \right) \ ,\\
  \lambda_{5,6}(\epsilon)&=\pm\im\left( 2\sqrt{I^*}\sqrt\epsilon + \dfrac{\epsilon^{3/2}}{\sqrt{I^*}} +\Oscr(\epsilon^{3/2}) \right)  \ .
\end{aligned}
$$
while the other two configurations are unstable with
\begin{displaymath}
\lambda_{(0,0,\pi)}=
\begin{cases}
\pm2^{5/4}\im\sqrt{I^*}\sqrt\epsilon+\Oscr(\epsilon)\ ,\\
\pm2^{5/4}\sqrt{I^*}\sqrt\epsilon+\Oscr(\epsilon)\ ,\\
\pm2\im\sqrt{I^*}\sqrt\epsilon+\Oscr(\epsilon),
\end{cases}
\qquad
\lambda_{(\pi,0,0)}
\begin{cases}
\pm2^{5/4}\im\sqrt{I^*}\sqrt\epsilon+\Oscr(\epsilon)\ ,\\
\pm2^{5/4}\sqrt{I^*}\sqrt\epsilon+\Oscr(\epsilon)\ ,\\
\pm2\im\sqrt{I^*}\sqrt\epsilon+\Oscr(\epsilon)\ .
\end{cases}\ .
\end{displaymath}
Concerning the effective stability of $(0,0,0)$, since there are two
couples of exponents which coincide at order $\Oscr(\sqrt\epsilon)$, but
different at order $\Oscr(\epsilon^{3/2})$, we have to take $\beta=3/2$ in the
assumption of Theorem~\ref{t.lin.stab.2}. Being $r+1-\alpha =
3-\frac12=\frac52$, the statement ensures existence of two couples of
distinct Floquet multipliers which are $\epsilon^{5/2}$-close to
$e^{\lambda T}$ on the unitary circle, which means linear stability of
the solution. Instead, for the other two unstable configurations we have
$\beta=\frac{1}{2}$ and $r+1-\alpha=\frac32$.

\smallskip
\noindent
\textbf{Isolated and degenerate solution.} The $(\pi,0,\pi)$ solution is
unstable with
\begin{displaymath}
\lambda_{(\pi,0,\pi)}=
\begin{cases}
\pm2^{3/2}\sqrt{I^*}\sqrt\epsilon+\Oscr(\sqrt{\epsilon})\ ,\\
\pm2\sqrt{I^*}\sqrt\epsilon+\Oscr(\sqrt{\epsilon})\ ,\\
\pm2\im\epsilon+\Oscr(\epsilon)\ .
\end{cases}
\end{displaymath}
Again, Theorem~\ref{t.lin.stab.2} applies with $\beta=1$ and
$r+1-\alpha=2$.

\smallskip
\noindent
\textbf{Degenerate solutions of the two families.} The remaining four
  configurations lying on the two families $P_{1,2}$ have two couples
  of Floquet exponents of order $\Oscr(\sqrt\epsilon)$ and one couple of
  exponents of order $\Oscr(\epsilon)$ due to the degenerate direction. In
  all the cases it is possible to verify the applicability of
  Theorem~\ref{t.lin.stab.2} with $\beta=1$ and $r+1-\alpha=2$, since
  the three couples of eigenvalues are all different at leading
  order, precisely
\begin{align*}
\lambda_{(\pi,\pi,\pi)}&=
\begin{cases}
\pm2\im \epsilon+\Oscr(\epsilon)\ ,\\
\pm\im \sqrt{I^*}\sqrt{2(\sqrt5-1)}\sqrt\epsilon+\Oscr(\sqrt\epsilon)\ ,\\
\pm \sqrt{I^*}\sqrt{2(\sqrt5+1)}\sqrt\epsilon+\Oscr(\sqrt\epsilon)\ ,
\end{cases}&
\lambda_{(0,\pi,0)}&=
\begin{cases}
\pm2\epsilon\\
\pm\im \sqrt{I^*}\sqrt{2({\sqrt5-1})}\sqrt\epsilon+\Oscr(\sqrt\epsilon)\ ,\\
\pm \sqrt{I^*}\sqrt{2({\sqrt5+1})}\sqrt\epsilon+\Oscr(\sqrt\epsilon)\ ,
\end{cases}\\
\\[-10pt]
\lambda_{(0,\pi,\pi)} &=
\begin{cases}
\pm2\im\epsilon+\Oscr(\epsilon)\ ,\\
\pm\im 2\sqrt{I^*}\sqrt\epsilon+\Oscr(\sqrt\epsilon)\ ,\\
\pm 2\sqrt{2I^*}\sqrt\epsilon+\Oscr(\sqrt\epsilon)\ ,
\end{cases}&
\lambda_{(\pi,\pi,0)} &=
\begin{cases}
\pm2\im\epsilon+\Oscr(\epsilon)\ ,\\
\pm\im 2\sqrt{I^*}\sqrt\epsilon+\Oscr(\sqrt\epsilon)\ ,\\
\pm 2\sqrt{2I^*}\sqrt\epsilon+\Oscr(\sqrt\epsilon)\ .
\end{cases}
\end{align*}

\bigskip
\subsection{Nonexistence of minimal square-vortexes in a dNLS railway-model}

\begin{center}
  \begin{tikzpicture}[scale = 0.8]
    \begin{scope} 
      \foreach \name/\posX/\posY in {C/3/0,D/4.5/0,E/6/0,F/7.5/0,G/9/1.,H/10.5/1.,I/12/0,L/13.5/0,M/15/0,N/16.5/0,
        T/3.5/1.5,U/5/1.5,V/6.5/1.5,W/8/1.5,Y/9.5/2.5,X/11/2.5,Z/12.5/1.5,J/14/1.5,K/15.5/1.5,a/17/1.5}
      \node (\name) at (\posX,\posY) {};
      \draw (D) circle (2mm) \foreach \n in {D,E,F,G,H,I,L,M,U,V,W,Y,X,Z,J,K} {(\n) circle (2mm)};
      \fill [black] \foreach \n in {G,H,Y,X} {(\n) circle (2.0mm)};
      \draw [red] (C) -- (D) -- (E) -- (F) -- (G)-- (H) -- (I) -- (L) -- (M) -- (N)  ;
      \draw [red] (T) -- (U) -- (V) -- (W) -- (Y)-- (X) -- (Z) -- (J) -- (K) -- (a)  ;
      \draw [red] (D) -- (U) (E) -- (V) (F) -- (W) (G)-- (Y) (H) -- (X) (I) -- (Z) (L) -- (J) (M) -- (K) ;
      \node[label=below:{\footnotesize-6}] at (D) {\phantom{$\times$}};
      \node[label=below:{\footnotesize-4}] at (E) {\phantom{$\times$}};
      \node[label=below:{\footnotesize-2}] at (F) {\phantom{$\times$}};
      \node[label=below:{\footnotesize 0}] at (G) {\phantom{$\times$}};
      \node[label=below:{\footnotesize 2}] at (H) {\phantom{$\times$}};
      \node[label=below:{\footnotesize 4}] at (I) {\phantom{$\times$}};
      \node[label=below:{\footnotesize 6}] at (L) {\phantom{$\times$}};
      \node[label=below:{\footnotesize 8}] at (M) {\phantom{$\times$}};
      \node[label=above:{\footnotesize-7}] at (U) {\phantom{$\times$}};
      \node[label=above:{\footnotesize-5}] at (V) {\phantom{$\times$}};
      \node[label=above:{\footnotesize-3}] at (W) {\phantom{$\times$}};
      \node[label=above:{\footnotesize-1}] at (Y) {\phantom{$\times$}};
      \node[label=above:{\footnotesize 1}] at (X) {\phantom{$\times$}};
      \node[label=above:{\footnotesize 3}] at (Z) {\phantom{$\times$}};
      \node[label=above:{\footnotesize 5}] at (J) {\phantom{$\times$}};
      \node[label=above:{\footnotesize 7}] at (K) {\phantom{$\times$}};
    \end{scope}
  \end{tikzpicture}
\end{center}
We here consider a minor variation of the Hamiltonian
system~\eqref{e.KdNLS}, the so-called \emph{railway-model}.  It
consists of two coupled dNLS models, where only nearest neighbours
interactions are active. The model, labelling the sites of the lattice according to the picture with
$\Jscr=\{-N,\ldots,N+1\}$, is described by the Hamiltonian
\begin{equation}
  \label{e.railway.ham}
\begin{aligned}
  H^{(0)} &= \sum_{j\in\Jscr}\left(\frac12 (x_j^2 +
    y_j^2) + \frac\gamma8 (x_j^2 + y_j^2)^2\right)
  +3\epsilon\sum_{j\in\Jscr}\frac12(x_j^2 + y_j^2)\\
    &\qquad-\epsilon
  \sum_{j\in\Jscr}(x_{j+1}x_{j-1}+ y_{j+1}y_{j-1})-\epsilon
 \sum_{j=-\lfloor\frac{N-1}2\rfloor}^{\lfloor\frac{N+1}2\rfloor} (x_{2j}x_{2j-1}+ y_{2j}y_{2j-1})\ .
\end{aligned}
\end{equation}
We want to investigate the continuation of the minimal vortex
configuration, namely the localised structures given by four
consecutive excited sites, that we here take as $S=\{-1,0,1,2\}$,
with phase differences between the neighbouring ones all equal to
$\pi/2$. The existence of such \emph{rotating} structures has been
shown in proper two-dimensional lattices in \cite{PelKF05b}, by
expanding at very high perturbation orders the Kernel equation
obtained with a Lyapunov-Schmidt reduction. On the other hand, in
\cite{PenSPKK18} similar structures have been proved not to exists in
the one-dimensional dNLS lattice~\eqref{e.ex.dnls} with
$\kappa_1=\kappa_3=1$, which at first orders in the perturbation
parameter $\epsilon$ exhibits the same averaged term $f_0^{(1,1)}(q)$
as the two-dimensional problem, hence the same critical points.  The
present \emph{railway-model} represents a natural hybrid setting
between one and two-dimensional square lattices.  Here, exploiting our
normal form construction, we are going to show the nonexistence of
the minimal vortex, thus enforcing the proper two-dimensional nature
of these kind of localised solutions.

We introduce action-angle variables
$(I,\theta)$ and complex coordinates $(\xi,\eta)$ and we expand $H_0$
and $H_1$ in Taylor series of the actions around $I^*$; 
forgetting constant terms and introducing the resonant angles
$\hat{q}=(q_1,q)$ and their conjugated actions $\hat{p}=(p_1,p)$
$$
 \left\lbrace \begin{aligned}
& q_1 = \theta_{-1} \ ,\\
& q_2 = \theta_0 - \theta_{-1} \ ,\\ 
& q_3 = \theta_1 - \theta_0 \ ,\\
& q_4 = \theta_2 - \theta_1 \ ,
\end{aligned}
\right.
\qquad
 \left\lbrace \begin{aligned}
& p_1 = J_{-1} + J_0 + J_1 + J_2 \ ,\\
& p_2 = J_0 + J_1 + J_2\ ,\\
& p_3= J_1 + J_2 \ ,\\
& p_4 = J_2\ ,
\end{aligned}
\right. \ ,
$$
we can rewrite the initial Hamiltonian as
$$
\begin{aligned}
H^{(0)} &= \omega p_1 + \sum_{j\in\Jscr\setminus S}{\im\xi_j\eta_j} +
f_4^{(0,0)}(\hat p, \xi,\eta) + f_0^{(0,1)}(q_2,q_3,q_4) +
f_1^{(0,1)}(\hat q,\xi,\eta)  \\ &\qquad+ f_2^{(0,1)}(q,\hat
p,\xi,\eta)+f_3^{(0,1)}(\hat q,\hat
p,\xi,\eta)+ f_4^{(0,1)}(q,\hat
p,\xi,\eta)+ \sum_{\ell\geq 5}f_\ell^{(0,1)}(\hat q,\hat
p,\xi,\eta) \ ,
\end{aligned} 
$$ with as usual $\omega=(1+I^*)$. The normal form at order one
$$
Z_0^{(1)}=f_0^{(1,1)}= 2I^*\epsilon\bigl(\cos(q_2)+\cos(q_2+q_3)+\cos(q_3+q_4)+\cos(q_4)\bigr)
$$
gives the two isolated critical points $(0,0,0)$,
$(\pi,0,\pi)$, and three one-parameter families
$$
P_1(\vartheta)=(\vartheta,\pi,-\vartheta)\ , \quad
P_2(\vartheta)=(\vartheta,\pi,\vartheta+\pi)\ , \quad
P_3(\vartheta)=(\vartheta,-2\vartheta,\vartheta+\pi)\ . 
$$ This is in agreement with the existing literature, see also the
example in \cite{PenSD18}. Let us notice that the three families
intersect in the two vortexes configurations
$\pm\left(\frac\pi{2},\pi,-\frac\pi{2}\right)$. These are completely
degenerate configurations, since the Kernel admits three independent
directions $\partial_{\theta}P_j$ on the tangent space to the torus
$\toro^3$; hence $D_q^2
Z_0^{(1)}\left(\frac\pi{2},\pi,-\frac\pi{2}\right)\equiv 0$. It is immediate
to verify that the two isolated configurations are nondegenerate,
hence we can apply Theorem \ref{teo:forma-normale-r} with $r=1$. For
the three degenerate families $P_j$ we have to compute the normal form
at order two. Similarly to the previous example on the Zigzag model,
the equation for the critical points of $Z_0^{(2)}$ takes the form
$$
F(q_2,q_3,q_4,\varepsilon)= F_0 (q_2,q_3,q_4)+ \varepsilon F_1 (q_2,q_3,q_4)=0,
$$
with $F:\mathbb{T}^3 \times \mathcal{U}(0) \rightarrow \mathbb{R}^3$.
In this case there are three vectors
$$
\partial_{\vartheta}P_1=\begin{pmatrix}
1 \\ 0 \\ -1
\end{pmatrix},
\qquad\qquad
\partial_{\vartheta}P_2=\begin{pmatrix}
1 \\ 0 \\ 1
\end{pmatrix}
\qquad\qquad
\partial_{\vartheta}P_3=\begin{pmatrix}
1 \\ -2 \\ 1
\end{pmatrix}
$$ that generate the Kernel of $D_qF_0(q)\Big|_{q=P_j(\vartheta)}$
with $j=1,2,3$.  The necessary condition for the solutions of $F_0$ to
be also solutions of $F$ is
$F_1(P_j(\vartheta))\perp\partial_\vartheta P_j(\vartheta)$; it turns
out that $F_1(P_1(\vartheta))\equiv 0$, hence nothing can be concluded
on $P_1$ (similarly to what already observed also in the dNLS cell in
\cite{PenSD18}), while for the other two families we get (apart from
a prefactor $\epsilon^2$)
$$
\langle F_1(P_2(\vartheta)), \partial_{\vartheta}P_2 \rangle =
 4 \sin(2\vartheta)\ ,\quad
 F_1(P_3(\vartheta)), \partial_{\vartheta}P_3 \rangle =
 4\sin(2\vartheta)\ ,
$$ and we can deduce that the two families $P_{2,3}(\vartheta)$ break
 down and either the four solutions $(0,\pi,\pi)$, $(\pi,\pi,0)$,
 $(0,0,\pi)$, $(\pi,0,0)$ or the two vortexes $\left(
 \frac{\pi}{2},\pi,\frac{3\pi}{2} \right)$ and $\left(
 \frac{3\pi}{2},\pi,\frac{\pi}{2} \right)$ are allowed.  The
 continuation of the four in/out-of-phase configurations to periodic
 orbits is ensured by~\eqref{e.small.eig} with $\alpha<\frac32$, explicit
 symbolic calculations gives $\alpha=1$.  In the two vortexes, instead,
 condition~\eqref{e.small.eig} is not fulfilled, since
 $P_1(\vartheta)$ is still a 1-parameter family of solutions for
 $\nabla_q Z_0^{(2)}=0$. Hence, a third normal form step is needed to
 study the continuation of the configurations in $P_1(\vartheta)$,
 vortexes included. The equation for the critical points of
 $Z_0^{(3)}$ takes the form
$$ F(q_2,q_3,q_4,\varepsilon)= F_0 (q_2,q_3,q_4)+ \varepsilon F_1
(q_2,q_3,q_4) + \varepsilon^2 F_2 (q_2,q_3,q_4)=0\ ,
$$
with
\begin{displaymath}
\langle F_2(P_1(\vartheta)), \partial_{\vartheta}P_1 \rangle
= \frac{4}{I^*}\sin(\vartheta)\ .
\end{displaymath}
The normal form at order three allows to prove the nonexistence of the
two vortex configurations, being $\sin(\vartheta)=0$ only for
$\vartheta=0,\pi$. Instead, the continuation of the last two
in/out-of-phase solutions $(\pi,\pi,\pi)$ and $(0,\pi,0)$ is ensured
by~\eqref{e.small.eig} with $\alpha=2$, explicit computations gives
$\alpha=3/2$.

Concerning the linear stability analysis, we summarise below the results for the different cases.

\smallskip
\noindent
\textbf{Isolated and nondegenerate solutions.} The two configurations
  $(0,0,0)$, $(\pi,0,\pi)$ have all approximate Floquet exponents of
  order $\Oscr(\sqrt\epsilon)$. In particular $(0,0,0)$ is the unique
  stable configuration, with Floquet exponents
$$
\begin{aligned}
\lambda_{1,2}(\epsilon)&=\pm\im\left( 2\sqrt{I^*}\sqrt\epsilon +
\dfrac{6}{(I^*)^{5/2}}\epsilon^{7/2} +\Oscr(\epsilon^{7/2}) \right)
\ ,\\ \lambda_{3,4}(\epsilon)&=\pm\im\left( 2\sqrt{I^*}\sqrt\epsilon
-\dfrac{\epsilon^{5/2}}{(I^*)^{3/2}}+ \Oscr(\epsilon^{5/2}) \right)
\ ,\\ \lambda_{5,6}(\epsilon)&=\pm\im\left( 2\sqrt2\sqrt{I^*}\sqrt\epsilon
+\dfrac{\sqrt{2} \epsilon^{3/2}}{\sqrt{I^*}} +\Oscr(\epsilon^{3/2})\right) \ .
	\end{aligned}
$$
which split only at order $\epsilon^{5/2}$. This leads to
$\beta=5/2$ in the assumption of Theorem~\ref{t.lin.stab.2}. Since
$r+1-\alpha = 4-\frac12=\frac72$, the statement ensures existence of
two couples of distinct Floquet multipliers which are
$\epsilon^{7/2}$-close to $e^{\lambda T}$ on the unitary circle, which
means linear stability of the solution. Also the Floquet exponents of the
unstable configurations coincide at leading order
\begin{displaymath}
\lambda_{(\pi,0,\pi)}=
\begin{cases}
\pm2\sqrt{I^*}\sqrt\epsilon+\Oscr(\sqrt\epsilon)\ ,\\
\pm2\sqrt{I^*}\sqrt\epsilon+\Oscr(\sqrt\epsilon)\ ,\\
\pm2\sqrt{2}\sqrt{I^*}\sqrt\epsilon+\Oscr(\sqrt\epsilon)\ ,
\end{cases}
\end{displaymath}
so that the normal form at order three is necessary to localise their deformation;
this however does not affect the instability of the true periodic
orbit.

\smallskip
\noindent
\textbf{Degenerate solutions of the family $P_1$.} The $(\pi,\pi,\pi)$ and $(0,\pi,0)$ solutions are unstable
and have approximate Floquet exponents
\begin{displaymath}
\lambda_{(0,\pi,0)}=
\begin{cases}
\pm2\im\sqrt{I^*}\sqrt\epsilon+\Oscr(\sqrt\epsilon)\ ,\\
\pm2\sqrt{I^*}\sqrt\epsilon+\Oscr(\sqrt\epsilon)\ ,\\
\pm\frac{2\im}{\sqrt{I^*}}\epsilon^{3/2}+\Oscr(\sqrt\epsilon^{3/2})\ ,
\end{cases}\qquad
\lambda_{(\pi,\pi,\pi)}=
\begin{cases}
\pm2\im\sqrt{I^*}\sqrt\epsilon+\Oscr(\sqrt\epsilon)\ ,\\
\pm2\sqrt{I^*}\sqrt\epsilon+\Oscr(\sqrt\epsilon)\ ,\\
\pm\frac{2}{\sqrt{I^*}}\epsilon^{3/2}+\Oscr(\sqrt\epsilon^{3/2})\ ,
\end{cases}
\end{displaymath}
also in this case a normal form at order three is needed in order to
apply Theorem~\ref{t.lin.stab.2} with $\beta=\frac32$ and
$r+1-\alpha=\frac52$.

\smallskip
\noindent
\textbf{Degenerate solutions of the family $P_{2}$ and $P_3$.} The four 
configurations
lying on the two families $P_{2}$ and $P_3$ all have the same two couples of
Floquet exponents of order $\Oscr(\sqrt\epsilon)$ and the same one couple
of exponents of order $\Oscr(\epsilon)$ related to the degenerate
direction
\begin{displaymath}
\lambda=
\begin{cases}
\pm2^{5/4}\im\sqrt{I^*}\sqrt\epsilon+\Oscr(\sqrt\epsilon)\ ,\\
\pm2^{5/4}\sqrt{I^*}\sqrt\epsilon+\Oscr(\sqrt\epsilon)\ ,\\
\pm2\im\epsilon+\Oscr(\epsilon)\ .
\end{cases}
\end{displaymath}
In these cases, a normal form at order one is enough since it is possible to verify the applicability of
Theorem~\ref{t.lin.stab.2} with $\beta=1$ and $r+1-\alpha=2$. Indeed, all the couples of
eigenvalues are different at leading order.

\begin{remark}
  The continuation of these localised solutions only requires to
  compute the normal form at order one or two, as previously
  explained.  Instead, the study of the stability of these solutions
  require, for some configurations, the computation of the normal form
  at order three in order to apply Theorem \ref{t.lin.stab.2}.  Indeed
  in these cases, the eigenvalues split at order $5/2$.  This
  highlights the power of the normal form approach, which allows to increase
  the accuracy of the approximation beyond the minimal order needed to
  ensure existence of the continuation.
\end{remark}

\subsection{Discrete solitons in the dNLS model with purely nonlinear coupling}

\begin{center}
	\begin{tikzpicture}[scale = 0.8]
	\begin{scope} 
	\foreach \name/\posX/\posY in {A/0/0,B/1/0,C/2/0,D/3/0,E/4/0,F/5/0,G/6/2,H/7/2,I/8/2,L/9/0,M/10/0,N/11/0,O/12/0,P/13/0,Q/14/0}
	\node (\name) at (\posX,\posY) {};
	\draw [dotted, black!50] (-1,0) -- (14,0);
	\draw [dotted, black!50] (-1,2) -- (14,2);
	\draw [dotted, black!50] (6,2) -- (6,0);
    \draw [dotted, black!50] (7,2) -- (7,0);
	\draw [dotted, black!50] (8,2) -- (8,0);
	\draw (B) circle (2mm) \foreach \n in {C,D,E,F,G,H,I,L,M,N,O,P} {(\n) circle (2mm)};
	\fill [black] \foreach \n in {G,H,I} {(\n) circle (2.0mm)};
	\fill [white] \foreach \n in {B,C,D,E,F,L,M,N,O,P} {(\n) circle (1.9mm)};
	\draw [red] (A) -- (B) -- (C) -- (D) -- (E) -- (F) -- (G)-- (H) -- (I) -- (L) -- (M) -- (N) -- (O) -- (P) -- (Q) ;
	\node[label=below:{\footnotesize-6}] at (B) {\phantom{$\times$}};
	\node[label=below:{\footnotesize-5}] at (C)
	{\phantom{$\times$}};
	\node[label=below:{\footnotesize-4}] at (D)
	{\phantom{$\times$}};
	\node[label=below:{\footnotesize-3}] at (E)
	{\phantom{$\times$}};
	\node[label=below:{\footnotesize-2}] at (F)
	{\phantom{$\times$}};
	\node[label=below:{\footnotesize2}] at (L)
	{\phantom{$\times$}};
	\node[label=below:{\footnotesize3}] at (M)
	{\phantom{$\times$}};
	\node[label=below:{\footnotesize4}] at (N)
	{\phantom{$\times$}};
	\node[label=below:{\footnotesize5}] at (O)
	{\phantom{$\times$}};
	\node[label=below:{\footnotesize6}] at (P)
	{\phantom{$\times$}};
	\node[label=below:{\footnotesize-1}] at (6,0) {\phantom{$\times$}};
	\node[label=below:{\footnotesize0}] at (7,0) {\phantom{$\times$}};
	\node[label=below:{\footnotesize1}] at (8,0) {\phantom{$\times$}};
	\node[label=left:{\footnotesize $0$}] at (-1,0) {};
	\node[label=left:{\footnotesize $I^*$}] at (-1,2) {};
	
	\end{scope}
	\end{tikzpicture}
\end{center}
We consider here a dNLS model slightly different from~\eqref{e.KdNLS},
with purely nonlinear coupling and, in its simplest form, only
nearest-neighbours interactions are active. It is well known that in
this model single-site discrete solitons (such as breathers in
Klein-Gordon models) are strongly localised, with tails decaying more
than exponentially fast (see, e.g.,
\cite{Fla94,RosS05}).  

Specifically, we consider a perturbation $H_1$ of the form
\begin{displaymath}
H_1 = \sum_{j\in\Jscr}|\psi_{j+1}-\psi_j|^4\ ,
\end{displaymath}
where $\Jscr=\{-N,\ldots,N\}$.  We want to investigate the continuation
of localised structure given by three consecutive sites, hence
corresponding to the set $S=\{-1,0,1\}$.  The perturbation $H_1$
is given by the quartic nearest neighbours interaction, which in real
coordinates reads
$$
\begin{aligned}
H_1 &=
\frac12\sum_{j\in\Jscr}(x_j^2+y_j^2)^2+\sum_{j\in\Jscr}(x_{j+1}x_j+y_{j+1}y_j)^2-\sum_{j\in\Jscr}(x_j^2+y_j^2)(x_{j+1}x_j+y_{j+1}y_j)\\
&\qquad -\sum_{j\in\Jscr}(x_j^2+y_j^2)(x_{j-1}x_j+y_{j-1}y_j) +\frac12 \sum_{j\in\Jscr}(x_{j+1}^2+y_{j+1}^2)(x_j^2+y_j^2)\ .
\end{aligned}
$$
Expanding $H_0$ and $H_1$ in Taylor series of the actions around
$I^*$, forgetting constant terms and introducing the resonant angles
and their conjugated actions 
$$
\left\lbrace \begin{aligned}
& q_1 = \theta_{-1} \ ,\\
& q_2 = \theta_0 - \theta_{-1} \ ,\\ 
& q_3 = \theta_1 - \theta_{0}\ ,
\end{aligned}
\right. \ ,
\qquad
\left\lbrace \begin{aligned}
& p_1 = J_{-1} + J_0 +J_1 \ , \\
& p_2 = J_0 + J_1 \ ,\\
& p_3 = J_1\ ,
\end{aligned}
\right. 
$$ 
we can rewrite the initial Hamiltonian as
$$
\begin{aligned}
H^{(0)} & = \omega p_1 + \sum_{j\in\Jscr\setminus S}{\im\xi_j\eta_j} +
f_4^{(0,0)}(\hat p, \xi,\eta) + f_0^{(0,1)}(q) + f_1^{(0,1)}(\hat q,\xi,\eta) +
f_2^{(0,1)}(\hat{q},\hat p,\xi,\eta)\\ 
& \quad + f_3^{(0,1)}(\hat q,\hat
p,\xi,\eta)+ f_4^{(0,1)}(\hat q,\hat p,\xi,\eta)+ \sum_{\ell\geq
	5}f_\ell^{(0,1)}(\hat q,\hat p,\xi,\eta) \ ,
\end{aligned} 
$$ where $\omega=1+I^*$. The normal form at order one gives
$$
Z_0^{(1)}=f_0^{(1,1)} =8 \left( I^*\right)^2\epsilon\left(\cos(2q_2)-\cos(q_2)+\cos(2q_3)-\cos(q_3)\right)\ ,
$$
thus the critical points are the four isolated solutions $(0,\pi)$, $(\pi,0)$,
$(0,\pi)$ and $(\pi,\pi)$. However (as already noticed in the isolated configurations of
the Zigzag model) the nondegeneracy condition is fulfilled only in
the last configuration $(\pi,\pi)$.  The remaining ones are all
degenerate extremizer of $Z_0^{(1)}$; indeed along the tangent
direction related to the zero variable(s) it is possible to observe an
asymptotic growth as $\Oscr(t^4)$. These represent further examples of
critical points which are degenerate, although being isolated.

The normal form at order two is not sufficient to remove the
degeneration for the other configurations, as $Z_0^{(2)} = Z_0^{(1)}
(1+\epsilon g(q))$ (the explicit expression of $g(q)$ is not relevant)
thus the critical points have exactly the same asymptotic behaviour near
$q_{2,3}=0$ and continuation is not granted since~\eqref{e.small.eig} is not satisfied.

The normal form at order three allows to prove the existence of the
degenerate configurations since~\eqref{e.small.eig} with $\alpha=3/2$
is satisfied.

The approximate stability analysis easily shows that $(0,0)$ is the
only stable configuration, with both Floquet exponents of order
$\Oscr(\epsilon^{3/2})$, precisely
$$
\begin{aligned}
& \lambda_{1,2}(\epsilon)= \pm 2\im I^*\left( \epsilon^{3/2} +
  \epsilon^{5/2}+\Oscr(\epsilon^{5/2})\right)\ , \\ & \lambda_{3,4}(\epsilon)= \pm 2\im
  I^*\left( \sqrt{3}\epsilon^{3/2} +\dfrac{1}{\sqrt{3}}
  \epsilon^{5/2}+\Oscr(\epsilon^{5/2})\right)\ .
\end{aligned}
$$ Theorem~\ref{t.lin.stab.2} applies after $r=3$ normal form steps
with $\beta=3/2<5/2=r+1-\alpha$, hence Floquet multipliers are
$\Oscr(\epsilon^{5/2})$-close to the approximate ones. The other three
configurations are all unstable, with 
\begin{displaymath}
\lambda_{(0,\pi)}\,,\  \lambda_{(\pi,0)}=
\begin{cases}
\pm7\sqrt6 I^*\im\epsilon^{\frac32}+\Oscr(\epsilon^{\frac32})\ ,\\
\pm4\sqrt2 I^*\sqrt\epsilon+\Oscr(\sqrt\epsilon)\ ,
\end{cases}
\lambda_{(\pi,\pi)}=
\begin{cases}
  \pm4 I^*\sqrt\epsilon+\Oscr(\sqrt\epsilon)\ ,\\
  \pm4\sqrt3 I^*\sqrt\epsilon+\Oscr(\sqrt\epsilon)\ .
\end{cases}
\end{displaymath}

\subsection{Other resonances and persistence of two dimensional tori.}

We consider the standard dNLS model~\eqref{e.ex.dnls} with
$\kappa_1=1$ and $\kappa_l=0$ for any ${2\leq l\leq d}$. At difference
with most of the literature on localised solutions, we now consider a
resonant torus with resonance different from the classical
$1:\text{\textendash}:1$.  In this case, the action of the symmetry group is
transversal to the action of the periodic flow on the unperturbed
torus; hence, any periodic orbit surviving to the continuation is not
isolated, being part of a 2-dimensional torus foliated by periodic
orbits, obtained by the action of the symmetry on one of the continued
periodic orbit. The objects which survive are then 2-dimensional
resonant subtori of the given resonant torus $I_l=I^*\ ,l\in S$ and
$\xi_l=\eta_l=0\ ,l\in \Jscr\setminus S$.

The normal form allows to approximate, at any finite order, the
subtori surviving to the breaking of the original resonant torus.  The
approximated invariant object allows to prove the persistence of the
considered subtorus.  The persistence of nondegenerate tori in
Hamiltonian systems with symmetries (and even in more generic
dynamical systems) is a known subject, see,
e.g.,~\cite{BamV01,BamG02,Bam15}. Unlike the quoted works, our
approach allows to treat both nondegenerate and degenerate subtori.
Let us remark that the continuation is here made at fixed period and
not at fixed values of the independent conserved quantities.

We now show how to construct the leading order approximation of these
subtori in both a nondegenerate and a degenerate case, in the easiest
case of three consecutive excited sites $S=\{-1,0,1\}$, always
assuming $\gamma=1$. Focusing on these examples, we also explain how
to modify the proof of Theorem~\ref{teo:forma-normale-r} in terms of
the map $\Upsilon$, so as to prove the persistence of these family of
localised and time periodic structures in dNLS models.

\subsection{Nondegenerate case.}
\begin{center}
	\begin{tikzpicture}[scale = 0.8]
	\begin{scope} 
	\foreach \name/\posX/\posY in {A/0/0,B/1/0,C/2/0,D/3/0,E/4/0,F/5/0,G/6/1,H/7/2.3,I/8/2.3,L/9/0,M/10/0,N/11/0,O/12/0,P/13/0,Q/14/0}
	\node (\name) at (\posX,\posY) {};
	\draw [dotted, black!50] (-1,0) -- (14,0);
	\draw [dotted, black!50] (-1,1) -- (14,1);
	\draw [dotted, black!50] (-1,2.3) -- (14,2.3);
	\draw [dotted, black!50] (6,1) -- (6,0);
	\draw [dotted, black!50] (7,2) -- (7,0);
	\draw [dotted, black!50] (8,2) -- (8,0);
	\draw (B) circle (2mm) \foreach \n in {C,D,E,F,G,H,I,L,M,N,O,P} {(\n) circle (2mm)};
	\fill [black] \foreach \n in {G,H,I} {(\n) circle (2.0mm)};
	\fill [white] \foreach \n in {B,C,D,E,F,L,M,N,O,P} {(\n) circle (1.9mm)};
	\draw [red] (A) -- (B) -- (C) -- (D) -- (E) -- (F) -- (G)-- (H) -- (I) -- (L) -- (M) -- (N) -- (O) -- (P) -- (Q) ;
	\node[label=below:{\footnotesize-6}] at (B) {\phantom{$\times$}};
	\node[label=below:{\footnotesize-5}] at (C)
	{\phantom{$\times$}};
	\node[label=below:{\footnotesize-4}] at (D)
	{\phantom{$\times$}};
	\node[label=below:{\footnotesize-3}] at (E)
	{\phantom{$\times$}};
	\node[label=below:{\footnotesize-2}] at (F)
	{\phantom{$\times$}};
	\node[label=below:{\footnotesize2}] at (L)
	{\phantom{$\times$}};
	\node[label=below:{\footnotesize3}] at (M)
	{\phantom{$\times$}};
	\node[label=below:{\footnotesize4}] at (N)
	{\phantom{$\times$}};
	\node[label=below:{\footnotesize5}] at (O)
	{\phantom{$\times$}};
	\node[label=below:{\footnotesize6}] at (P)
	{\phantom{$\times$}};
	\node[label=below:{\footnotesize-1}] at (6,0) {\phantom{$\times$}};
	\node[label=below:{\footnotesize0}] at (7,0) {\phantom{$\times$}};
	\node[label=below:{\footnotesize1}] at (8,0) {\phantom{$\times$}};
	\node[label=left:{\footnotesize $0$}] at (-1,0) {};
	\node[label=left:{\footnotesize $I^*$}] at (-1,1) {};
	\node[label=left:{\footnotesize $1+2I^*$}] at (-1,2.3) {};	
	\end{scope}
	\end{tikzpicture}
\end{center}
Consider the set $S=\{-1,0,1\}$ with excited actions equal to $\{I^*,1+2I^*,1+2I^*\}$, so that at
$\epsilon=0$ the flow lies on a resonant torus with frequencies
$\hat\omega = \omega(1,2,2)$, where $\omega=1+I^*$.
After expanding $H_0$ and $H_1$ in Taylor series of the actions
around $I_l^*$, with $l=S$, we introduce the resonant angles
$\hat{q}=(q_1,q)$ and their conjugated actions $\hat{p}=(p_1,p)$ as
follows
$$
 \left\lbrace \begin{aligned}
& q_1 = \theta_{-1} \ ,\\
& q_2 = \theta_0 - 2\theta_{1} \ ,\\ 
& q_3 = \theta_1 - \theta_0 \ ,\\
\end{aligned}
\right.
\qquad
 \left\lbrace \begin{aligned}
& p_1 = J_{-1} + 2J_0 + 2J_1\ , \\
& p_2 = J_0 + J_1\ ,\\
& p_3= J_1 \ .\\
\end{aligned}
\right.
$$
so that we can rewrite the initial Hamiltonian in the form
$$
\begin{aligned}
H^{(0)} &= \omega p_1 + \sum_{j\in\Jscr\setminus S}{\im\xi_j\eta_j} +
f_4^{(0,0)}(\hat p, \xi,\eta) + f_0^{(0,1)}(\hat q) + f_1^{(0,1)}(\hat
q,\xi,\eta) + \\ &\qquad+ f_2^{(0,1)}(\hat q,\hat
p,\xi,\eta)+f_3^{(0,1)}(\hat q,\hat p,\xi,\eta)+ f_4^{(0,1)}(\hat
q,\hat p,\xi,\eta)+ \sum_{\ell\geq 5}f_\ell^{(0,1)}(\hat q,\hat p,\xi,\eta)
\ .
\end{aligned} 
$$
The normal form at order one gives
\begin{displaymath}
Z_0^{(1)} = f_0^{(1,1)}=-2\epsilon(1+2I^*)\cos(q_3)\ ,
\end{displaymath}
whose critical points are only $q_3=0,\,\pi$, which correspond to two
invariant subtori, foliated by periodic orbits $(q_1=\omega
t+q_1(0),q_2=q_2(0))$. Let us stress that the absence of the resonant
angle $q_2$ has not to be interpreted as the effect of a proper
degeneracy, since we expect a finite number of 2-dimensional subtori
to be continued; thus the two subtori are clearly nondegenerate. The
subtorus $q_3=\pi$ is linearly unstable, its approximate Floquet
exponents are $\lambda=
\pm2\sqrt{1+2I^*}\sqrt{\epsilon}+\Oscr(\sqrt{\epsilon})$, while
$q_3=0$ is linearly elliptic, with
$$
\lambda_{1,2}(\epsilon)= \pm\im\left( 2\sqrt{1+2I^*}\sqrt\epsilon+ \dfrac{\epsilon^{3/2}}{\sqrt{1+2I^*}} ++\Oscr(\epsilon^{3/2}) \right) \ .
$$ Theorem~\ref{t.lin.stab.2} can be applied also in this case with
$\beta=1/2<r+1-\alpha=1$, obtaining the effective linear stability of
the torus.

In order to prove the persistence of the obtained subtori, one can
keep both $q_1(0)$ and $q_2(0)$ as parameters in the map $\Upsilon$
introduced in~\eqref{frm:Ups}, and forget the variation of the second
action $p_2$ (since in this case we have two independent constant of
motion); hence $\Upsilon:\reali^{2n-2}\to\reali^{2n-2}$. Coherently
with such a definition of $\Upsilon$, and under the same assumptions
of Theorem~\ref{teo:forma-normale-r} on the spectrum of
$M(\epsilon)=\Upsilon'(x^*)$, existence and approximation of the
considered subtori are derived via the same Newton-Kantorovich
method. In our case existence can be obtained since
$\alpha=\frac12<=1$ for the smallest eigenvalue of $M(\epsilon)$.

\subsection{Degenerate case}

\begin{center}
	\begin{tikzpicture}[scale = 0.8]
	\begin{scope} 
	\foreach \name/\posX/\posY in {A/0/0,B/1/0,C/2/0,D/3/0,E/4/0,F/5/0,G/6/1,H/7/2.3,I/8/1,L/9/0,M/10/0,N/11/0,O/12/0,P/13/0,Q/14/0}
	\node (\name) at (\posX,\posY) {};
	\draw [dotted, black!50] (-1,0) -- (14,0);
	\draw [dotted, black!50] (-1,1) -- (14,1);
	\draw [dotted, black!50] (-1,2.3) -- (14,2.3);
	\draw [dotted, black!50] (6,1) -- (6,0);
	\draw [dotted, black!50] (7,2) -- (7,0);
	\draw [dotted, black!50] (8,1) -- (8,0);
	\draw (B) circle (2mm) \foreach \n in {C,D,E,F,G,H,I,L,M,N,O,P} {(\n) circle (2mm)};
	\fill [black] \foreach \n in {G,H,I} {(\n) circle (2.0mm)};
	\fill [white] \foreach \n in {B,C,D,E,F,L,M,N,O,P} {(\n) circle (1.9mm)};
	\draw [red] (A) -- (B) -- (C) -- (D) -- (E) -- (F) -- (G)-- (H) -- (I) -- (L) -- (M) -- (N) -- (O) -- (P) -- (Q) ;
	\node[label=below:{\footnotesize-6}] at (B) {\phantom{$\times$}};
	\node[label=below:{\footnotesize-5}] at (C)
	{\phantom{$\times$}};
	\node[label=below:{\footnotesize-4}] at (D)
	{\phantom{$\times$}};
	\node[label=below:{\footnotesize-3}] at (E)
	{\phantom{$\times$}};
	\node[label=below:{\footnotesize-2}] at (F)
	{\phantom{$\times$}};
	\node[label=below:{\footnotesize2}] at (L)
	{\phantom{$\times$}};
	\node[label=below:{\footnotesize3}] at (M)
	{\phantom{$\times$}};
	\node[label=below:{\footnotesize4}] at (N)
	{\phantom{$\times$}};
	\node[label=below:{\footnotesize5}] at (O)
	{\phantom{$\times$}};
	\node[label=below:{\footnotesize6}] at (P)
	{\phantom{$\times$}};
	\node[label=below:{\footnotesize-1}] at (6,0) {\phantom{$\times$}};
	\node[label=below:{\footnotesize0}] at (7,0) {\phantom{$\times$}};
	\node[label=below:{\footnotesize1}] at (8,0) {\phantom{$\times$}};
	\node[label=left:{\footnotesize $0$}] at (-1,0) {};
	\node[label=left:{\footnotesize $I^*$}] at (-1,1) {};
	\node[label=left:{\footnotesize $1+2I^*$}] at (-1,2.3) {};	
	\end{scope}
	\end{tikzpicture}
\end{center}
Consider the set $S=\{-1,0,1\}$ with excited actions equal to
$\{I^*,1+2I^*,I^*\}$, so that at $\epsilon=0$ the flow lies on a
resonant torus with frequencies $\hat\omega = \omega(1,2,1)$, where
again $\omega=1+I^*$. After expanding $H_0$ and $H_1$ in Taylor series
of the actions around $I_l^*$, with $l=S$, we introduce the resonant
angles $\hat{q}=(q_1,q)$ and their conjugated actions
$\hat{p}=(p_1,p)$ as
$$
 \left\lbrace \begin{aligned}
& q_1 = \theta_{-1} \ ,\\
& q_2 = \theta_0 - 2\theta_{-1} \ ,\\ 
& q_3 = \theta_1 - \theta_{-1} \ ,\\
\end{aligned}
\right. 
\qquad
 \left\lbrace \begin{aligned}
& p_1 = J_{-1} + 2J_0 + J_1 \ ,\\
& p_2 = J_0\ ,\\
& p_3= J_1 \ ,\\
\end{aligned}
\right.
$$
so that we can rewrite the initial Hamiltonian in the form
$$
\begin{aligned}
H^{(0)} & = \omega p_1 + \sum_{j\in\Jscr\setminus S}{\im\xi_j\eta_j} +
f_4^{(0,0)}(\hat p, \xi,\eta) + f_0^{(0,1)}(\hat{q}) +
f_1^{(0,1)}(\hat q,\xi,\eta) \\ 
&\qquad + f_2^{(0,1)}(\hat q,\hat
p,\xi,\eta)+f_3^{(0,1)}(\hat q,\hat
p,\xi,\eta)+ f_4^{(0,1)}(\hat q,\hat
p,\xi,\eta)+ \sum_{\ell\geq 5}f_\ell^{(0,1)}(\hat q,\hat
p,\xi,\eta) \ .
\end{aligned} 
$$

Unlike the previous example, the normal form at order one gives $
Z_0^{(1)} \equiv 0$, since the two resonant oscillators at sites
$\{-1,1\}$ are not interacting at order $\Oscr(\epsilon)$.
The normal form at order two gives
\begin{displaymath}
Z_0^{(2)} = f_0^{(2,2)}=\frac{2 (I^*)^2}{(1+I^*)^2}\epsilon^2\cos(q_3)\ ,
\end{displaymath}
whose critical points are $q^*_3=0,\,\pi$.  Explicit calculations with
Mathematica, provide the expected value of $\alpha=1<\frac32$, this
allows to prove the continuation of the subtori.  The subtorus
$q^*_3=0$ is linearly unstable, while $q^*_3=\pi$ is linearly stable 
with
$$
\lambda_{1,2}(\epsilon)= \pm \im \left( \frac{2I^*}{1+I^*} \epsilon+
\dfrac{ 2I^*\left( 5 + 22 I^* + 24 (I^*)^2\right)\epsilon^3}{(1 + I^*)^5}+
\Oscr(\epsilon^3)\right) \ .
$$
Theorem~\ref{t.lin.stab.2} can be applied with
$\beta=1<r+1-\alpha=2$, thus proving the effective linear stability.

\section{Numerical simulations: a case study.}
\label{s:3}

In this section, we numerically investigate some aspects of the normal
form construction, previously applied to different dNLS models and
spatial configurations, focusing on a single case study: the
multi-pulse solution with three excited sites $S=\{-2,-1,1\}$ in the
standard dNLS model.  In particular we numerically highlight
\begin{enumerate}[label=(\roman*)]
\item the increase of the approximation accuracy as the order of the
  normal form is increased, by comparing $r=2$ with $r=3$;

\item the linear stability properties of the approximate periodic
  orbit in the stable case $q^*=(0,\pi)$ and in two different unstable
  ones, $q^*=(0,0)$ and $q^*=(\pi,\pi)$, both having only one unstable
  direction but with different orders in $\epsilon$.
\end{enumerate}

\subsection{Accuracy of the approximate periodic orbit}

\begin{figure}[!ht]
\centering
\includegraphics[width=0.4\textwidth]{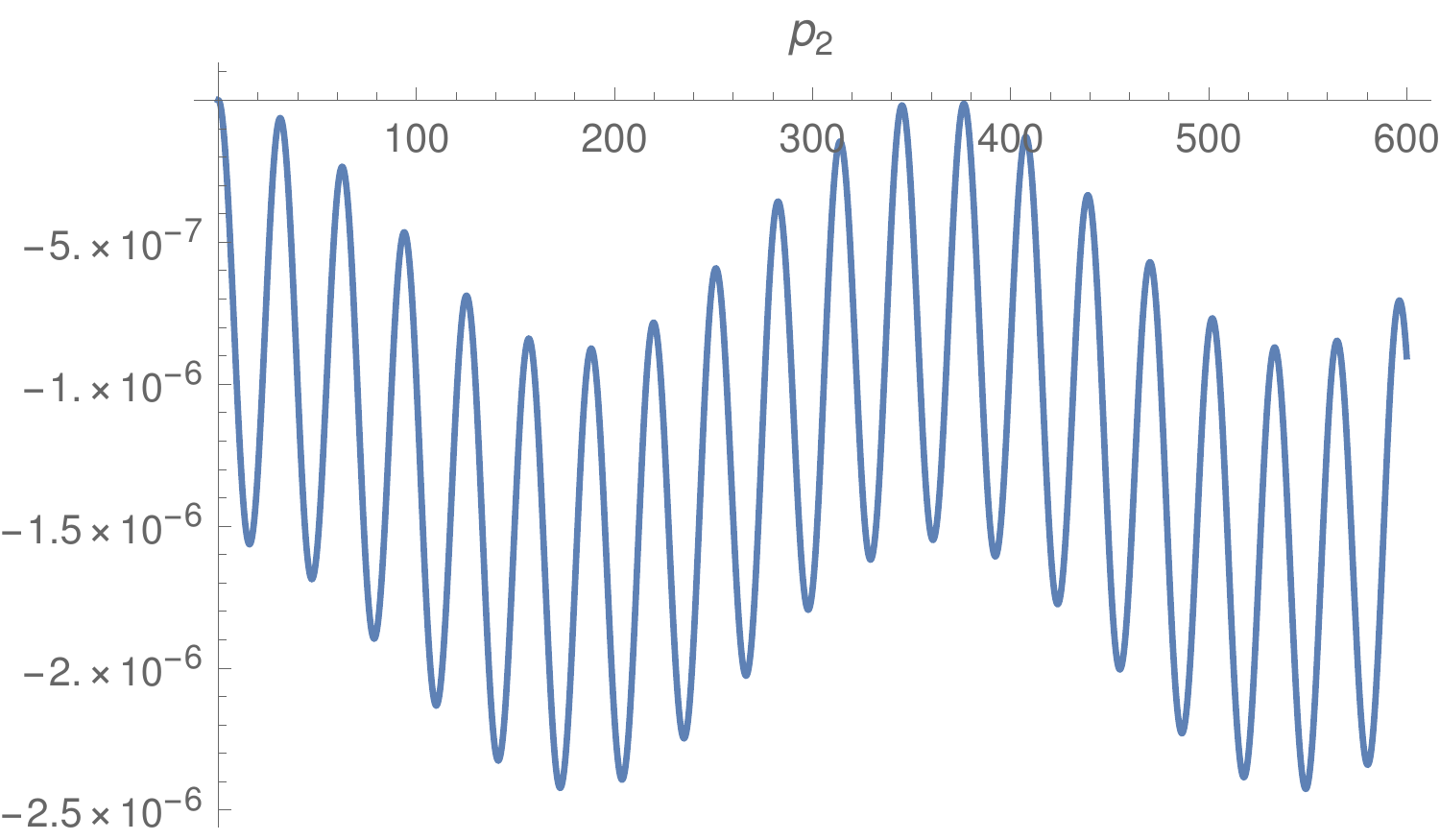}
\includegraphics[width=0.4\textwidth]{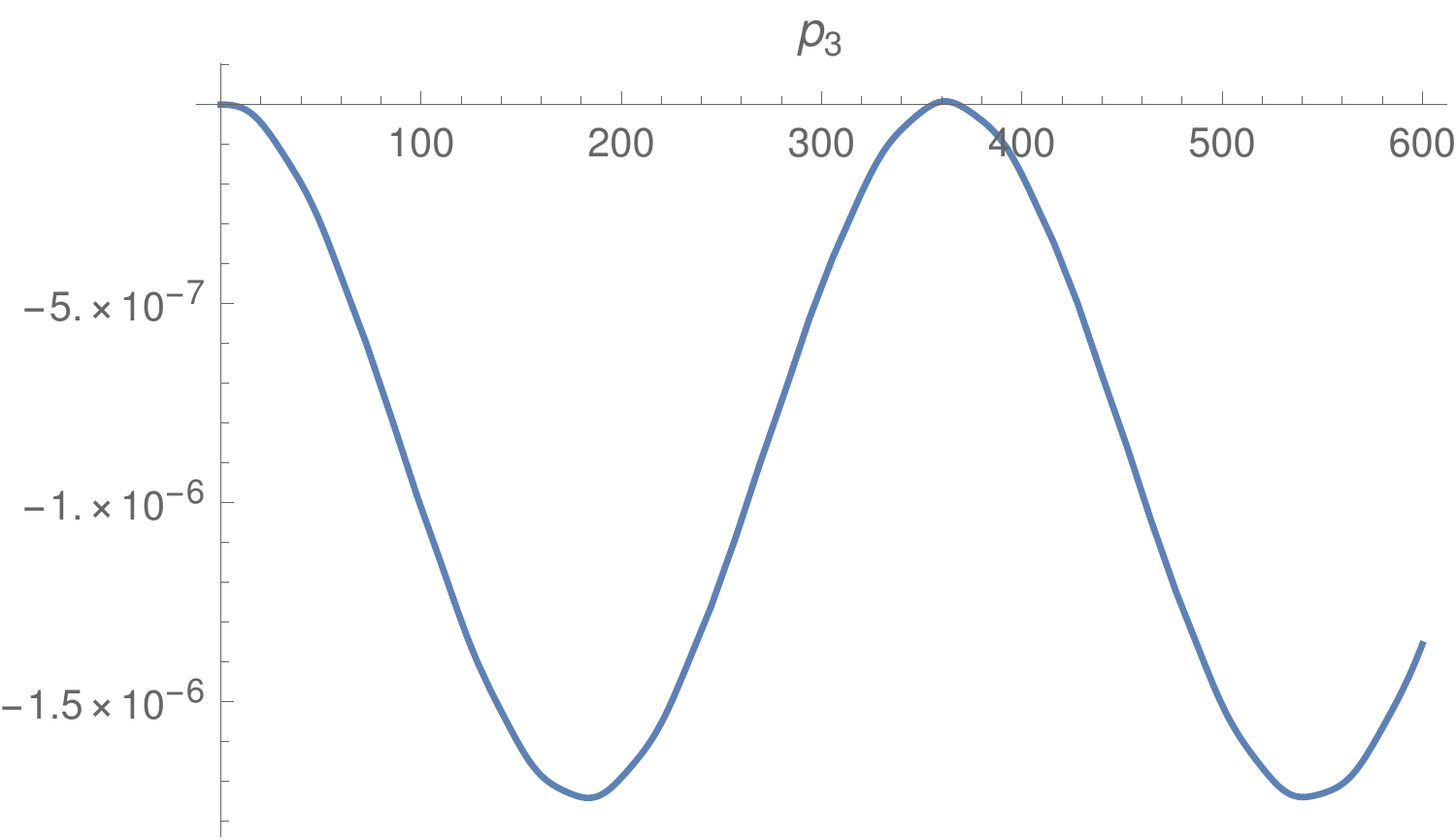}
\includegraphics[width=0.4\textwidth]{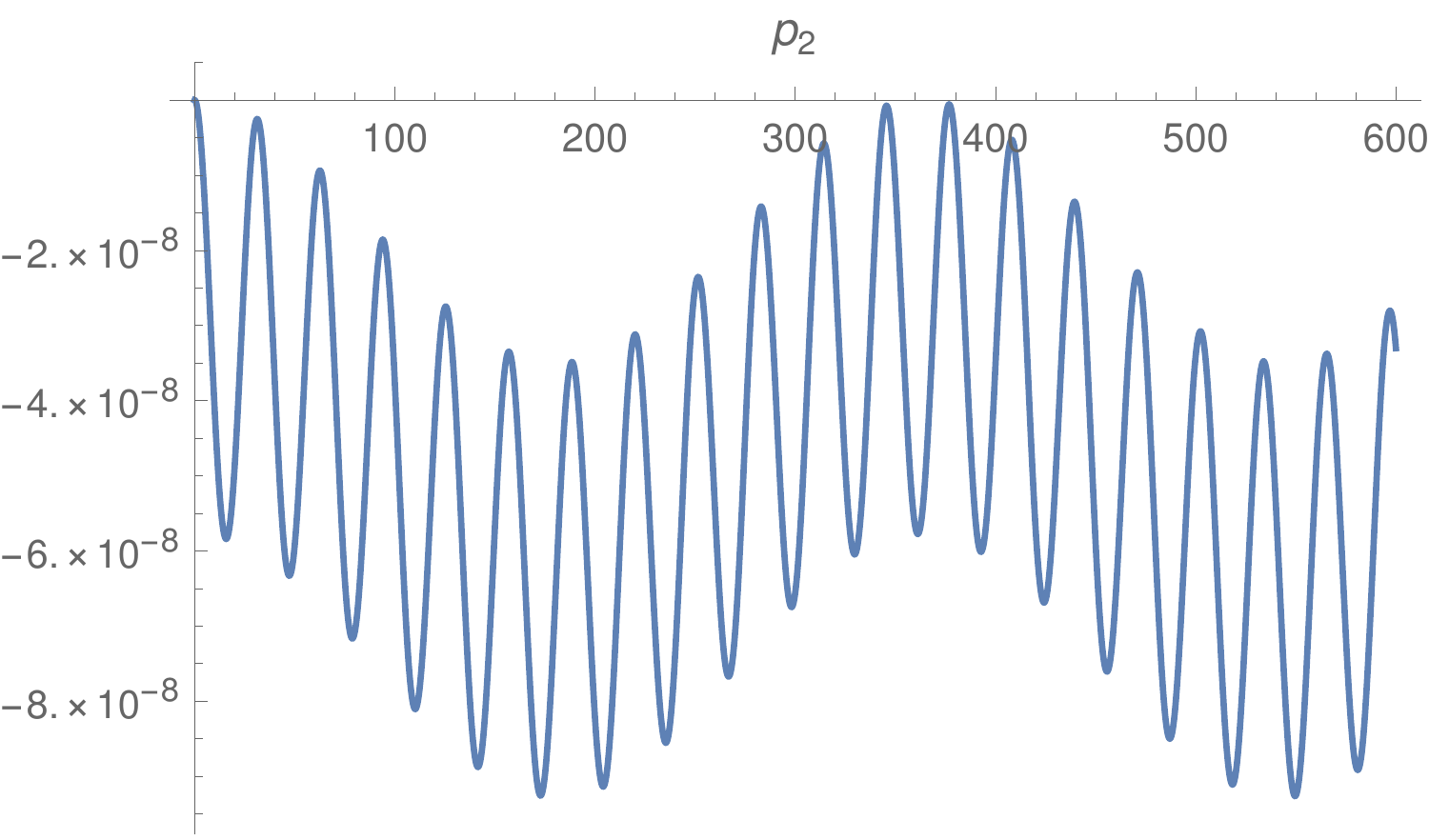}
\includegraphics[width=0.4\textwidth]{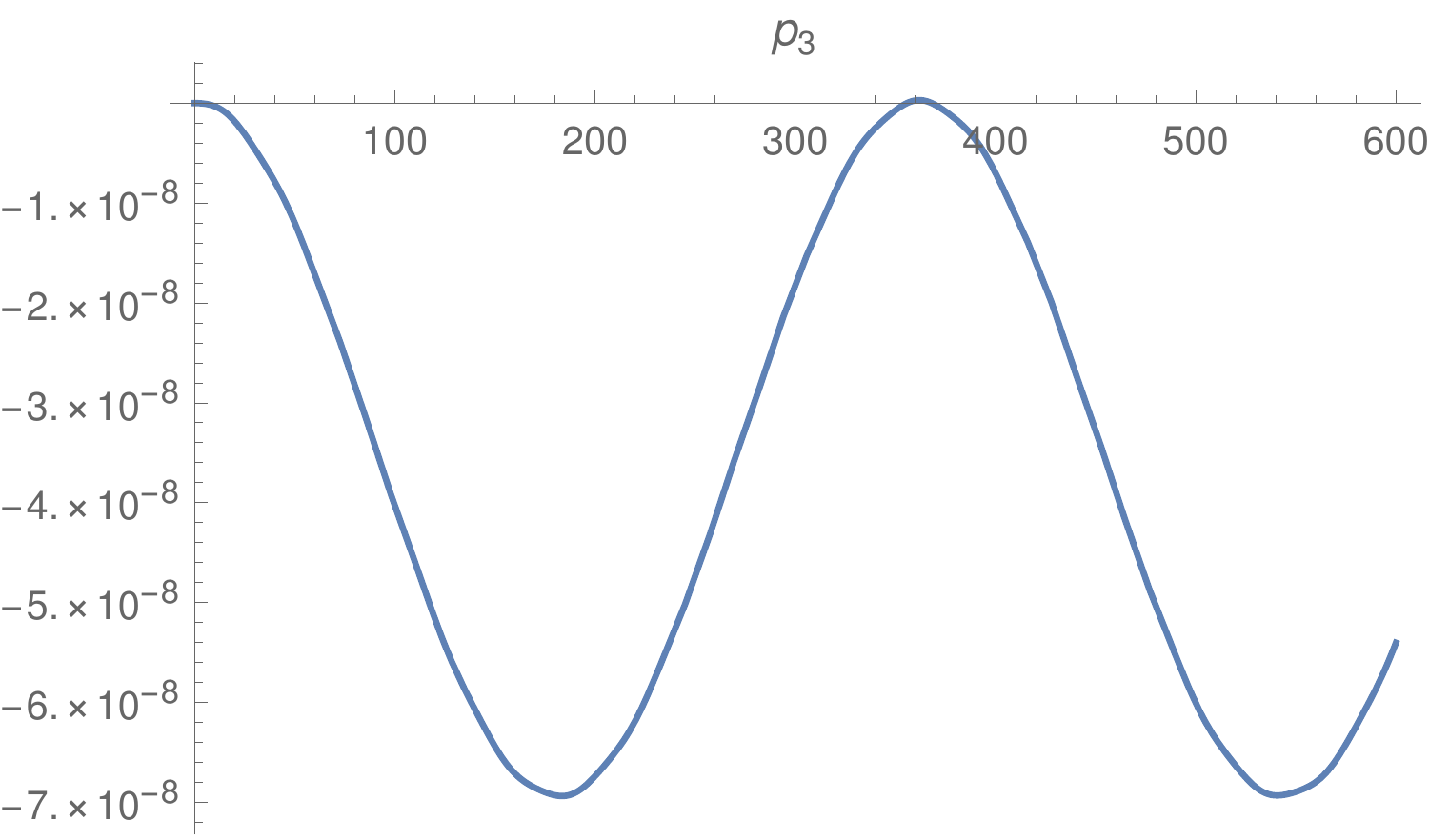}
\caption{Time evolution over the interval $[0,\Oscr(\epsilon^{-1})]$ of the variables
  $p_{2}(t)$ and $p_{3}(t)$. Comparison between $r=2$ (top) and $r=3$ (bottom).}
\label{f.0Pi.comp_p}
\end{figure}

\begin{figure}[!ht]
\centering
\includegraphics[width=0.4\textwidth]{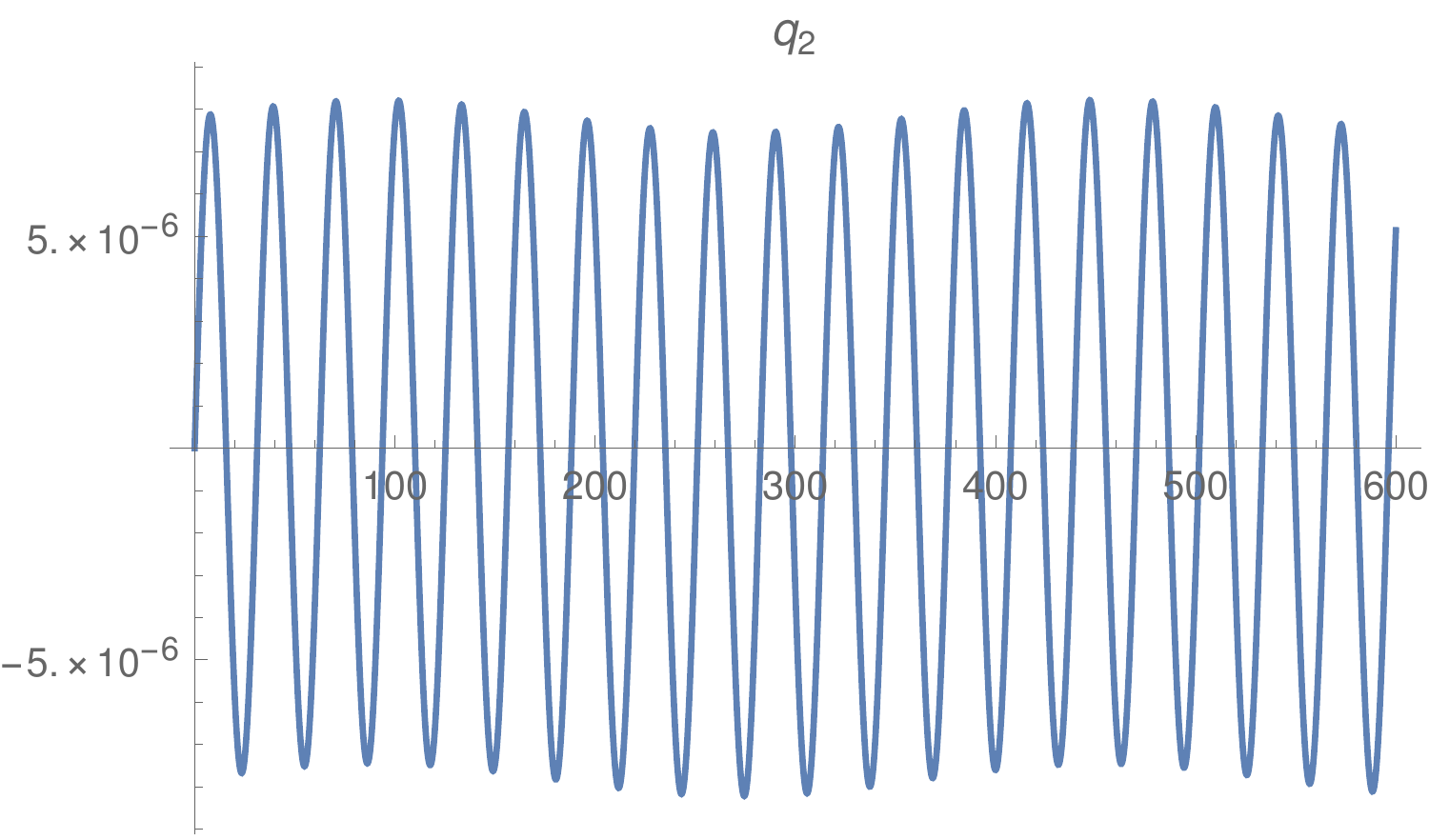}
\includegraphics[width=0.4\textwidth]{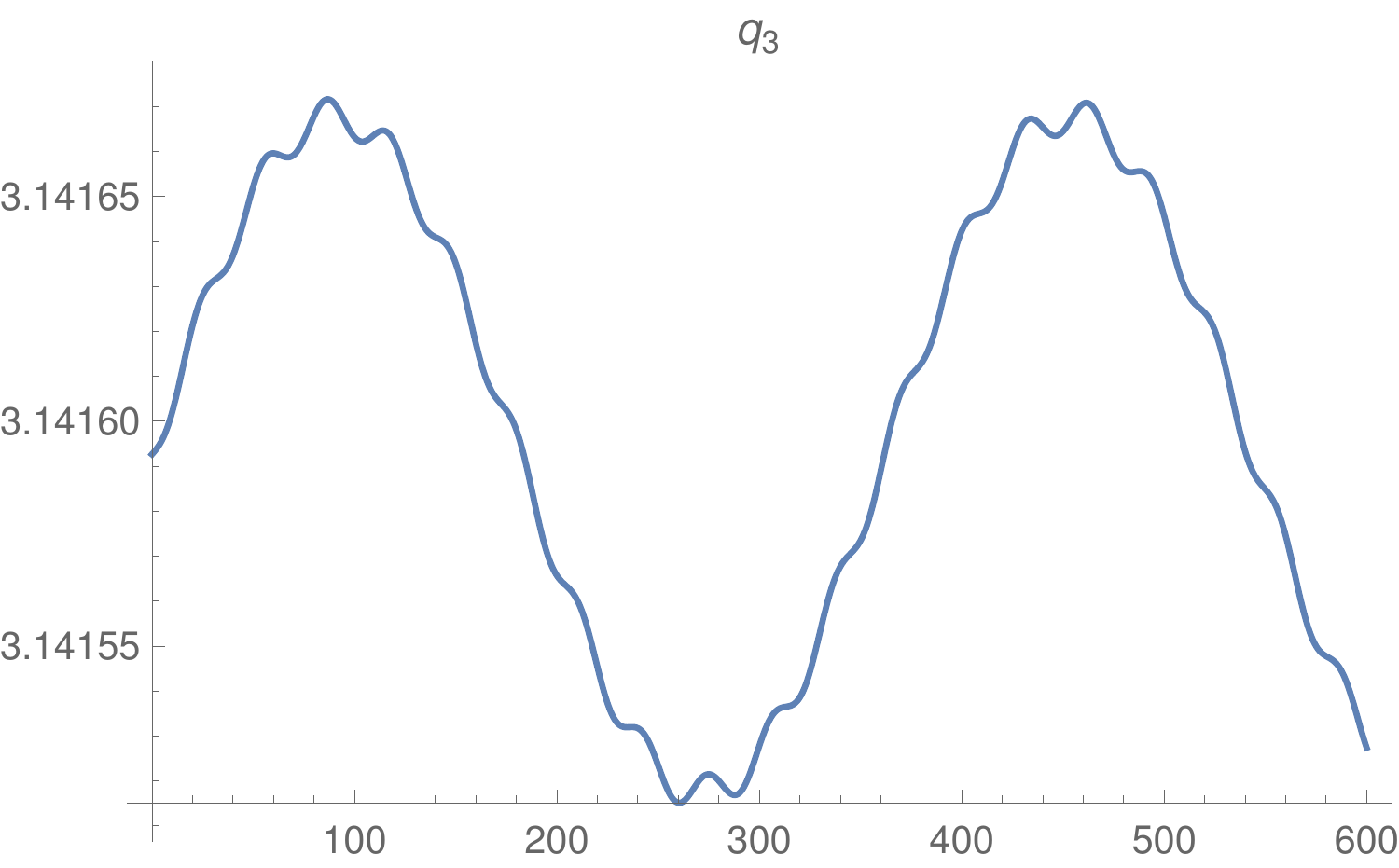}
\includegraphics[width=0.4\textwidth]{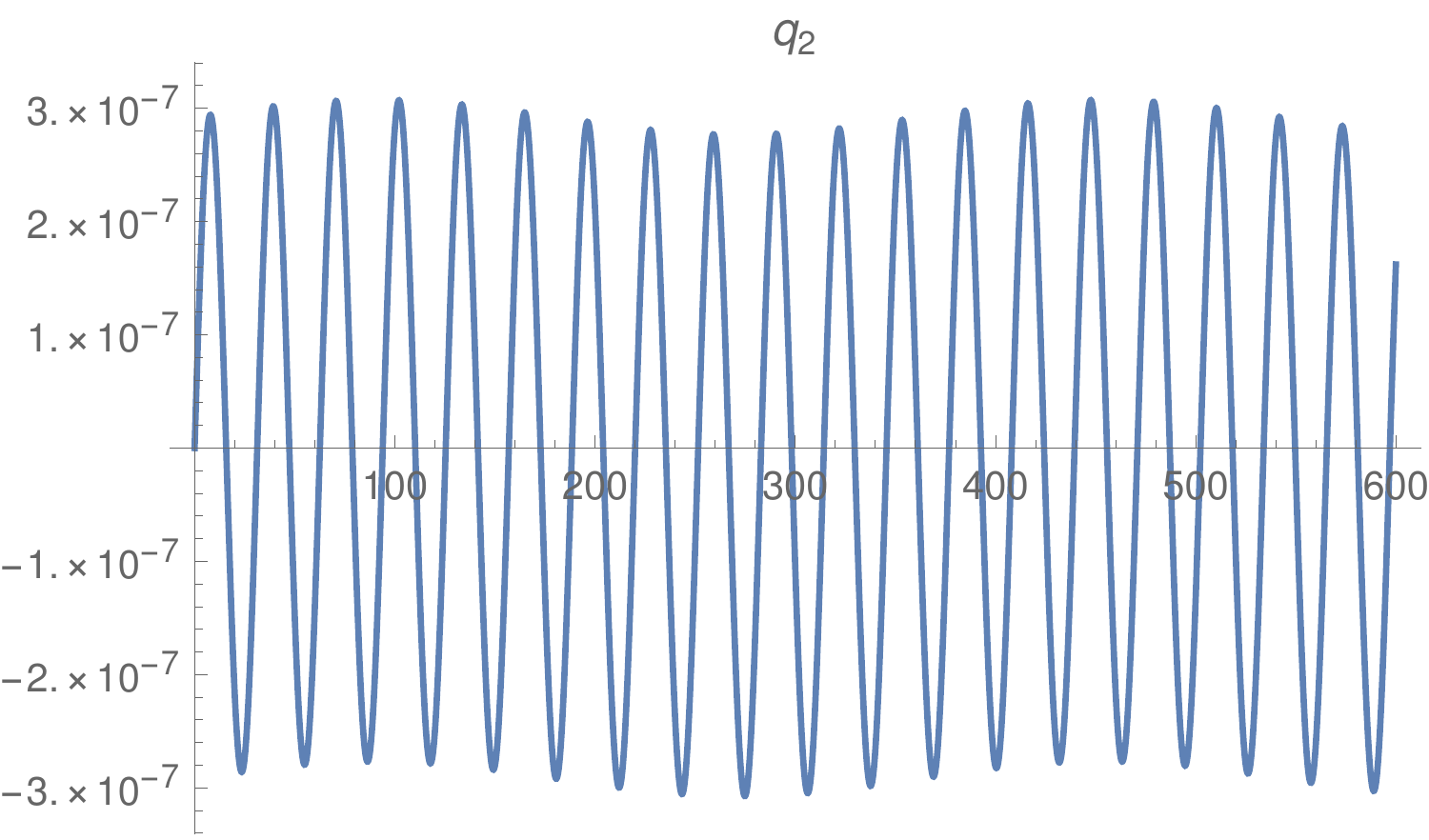}
\includegraphics[width=0.4\textwidth]{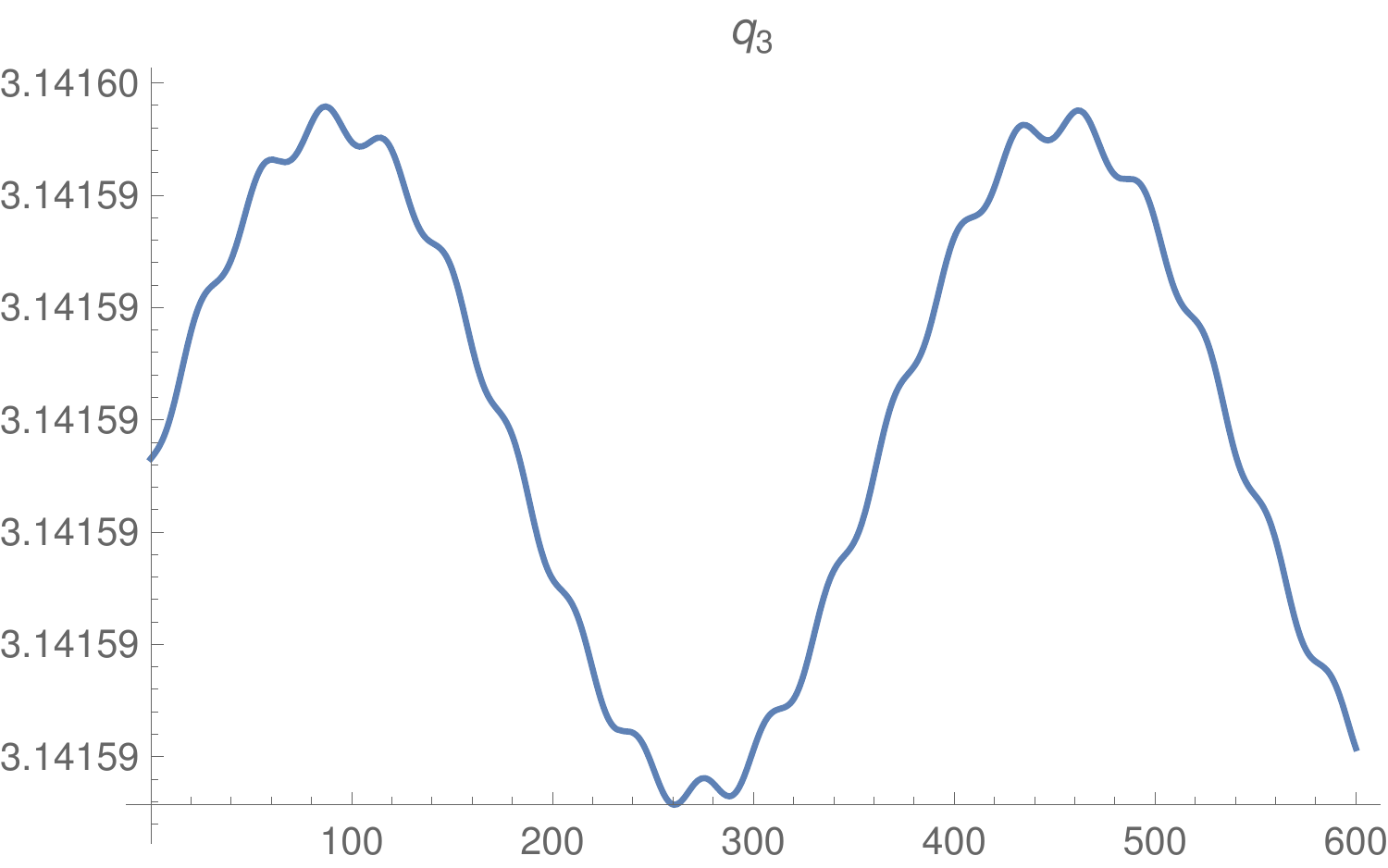}
\caption{Time evolution over the interval $[0,\Oscr(\epsilon^{-1})]$ of the variables
  $q_{2}(t)$ and $q_{3}(t)$. Comparison between $r=2$ (top) and $r=3$ (bottom).}
\label{f.0Pi.comp_q}
\end{figure}

\begin{figure}[!ht]
\centering
\includegraphics[width=0.4\textwidth]{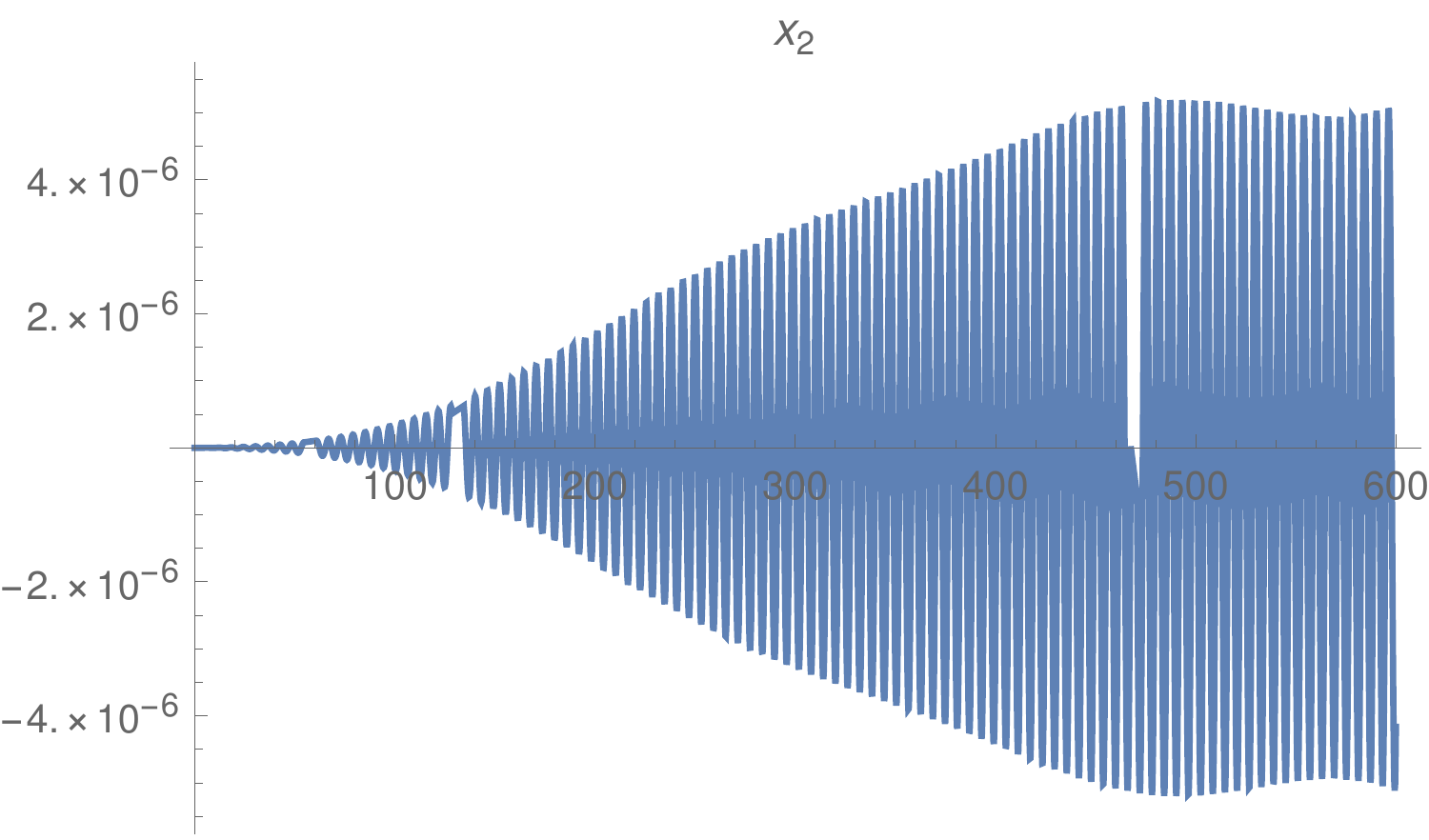}
\includegraphics[width=0.4\textwidth]{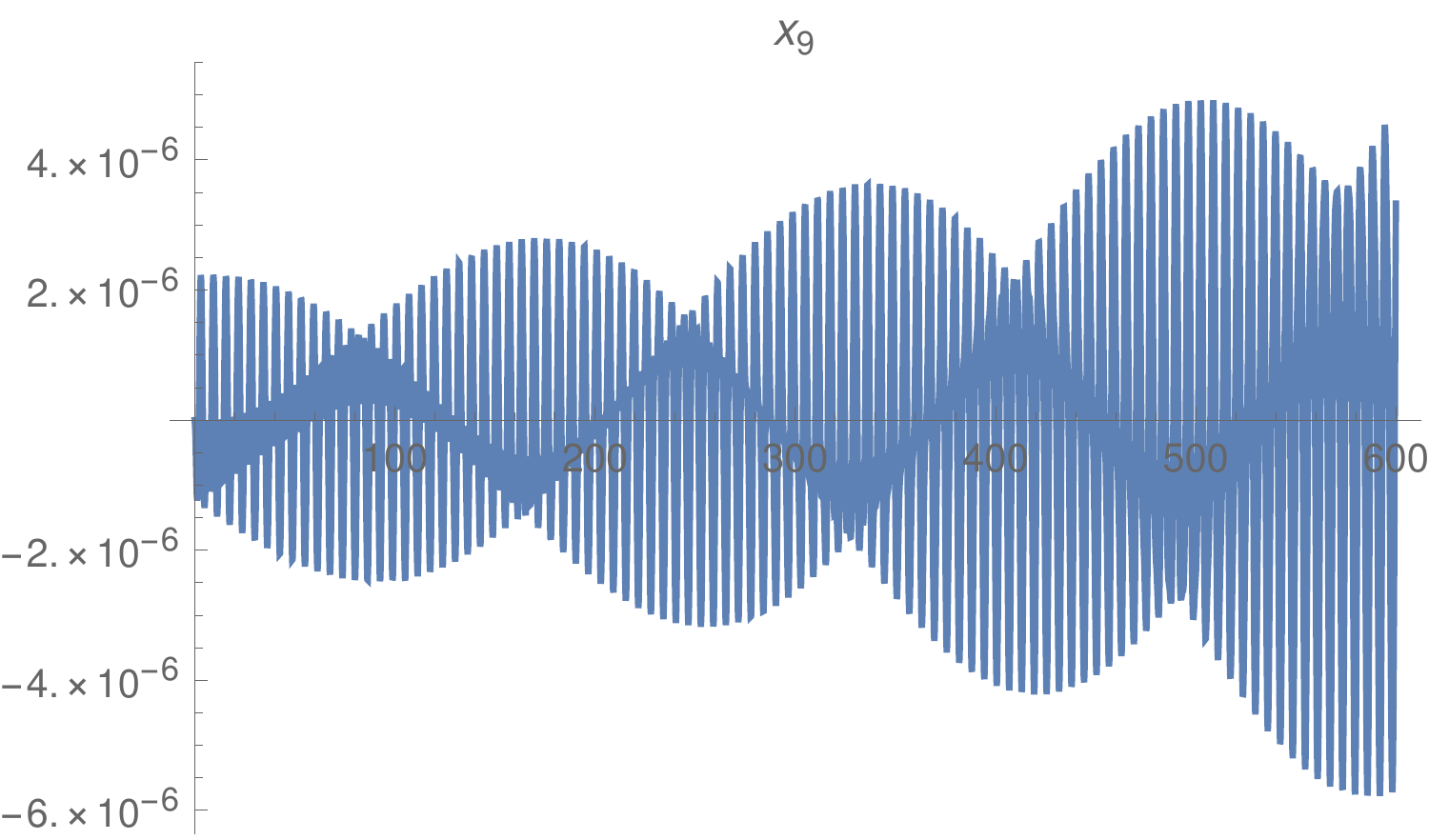}
\includegraphics[width=0.4\textwidth]{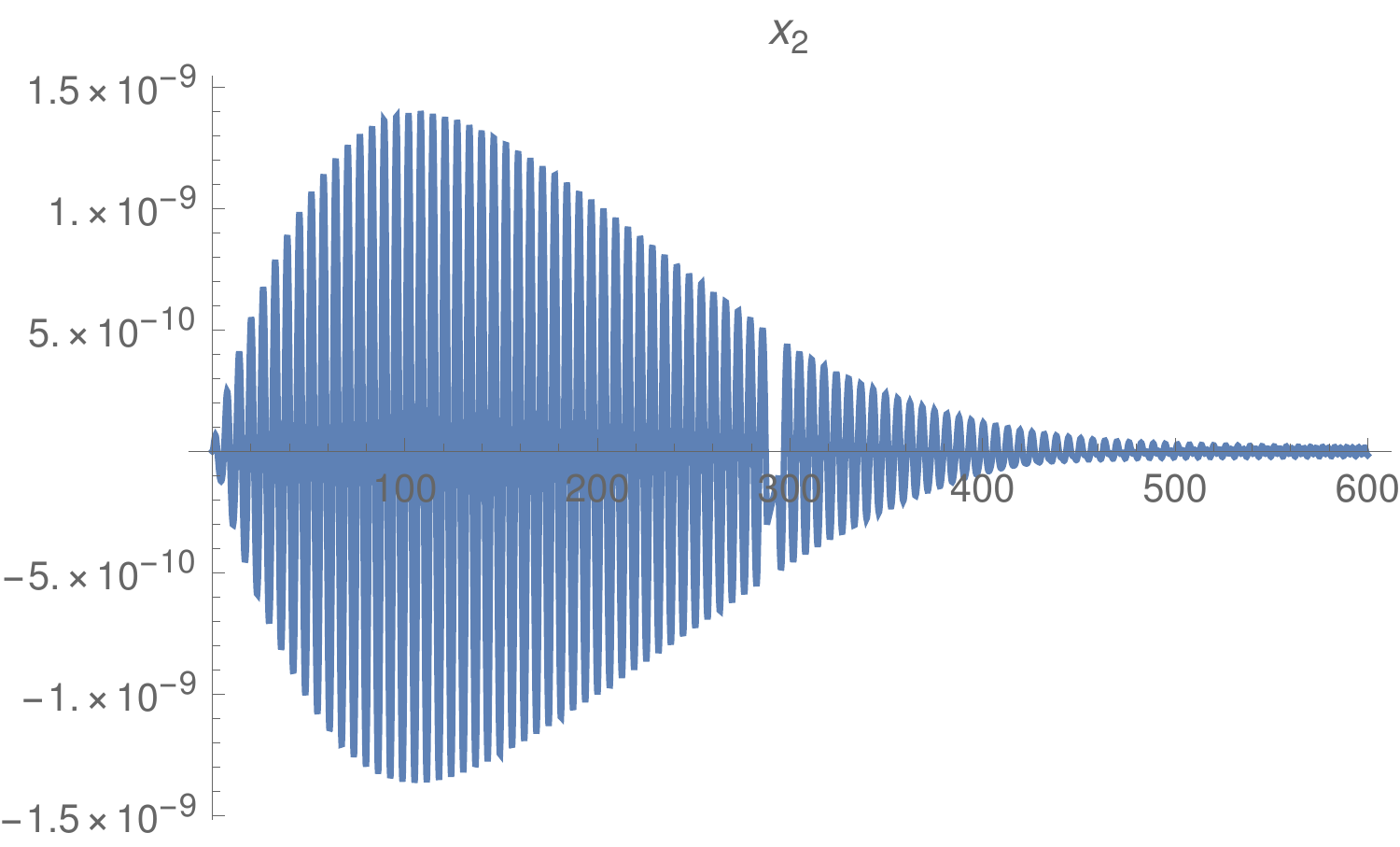}
\includegraphics[width=0.4\textwidth]{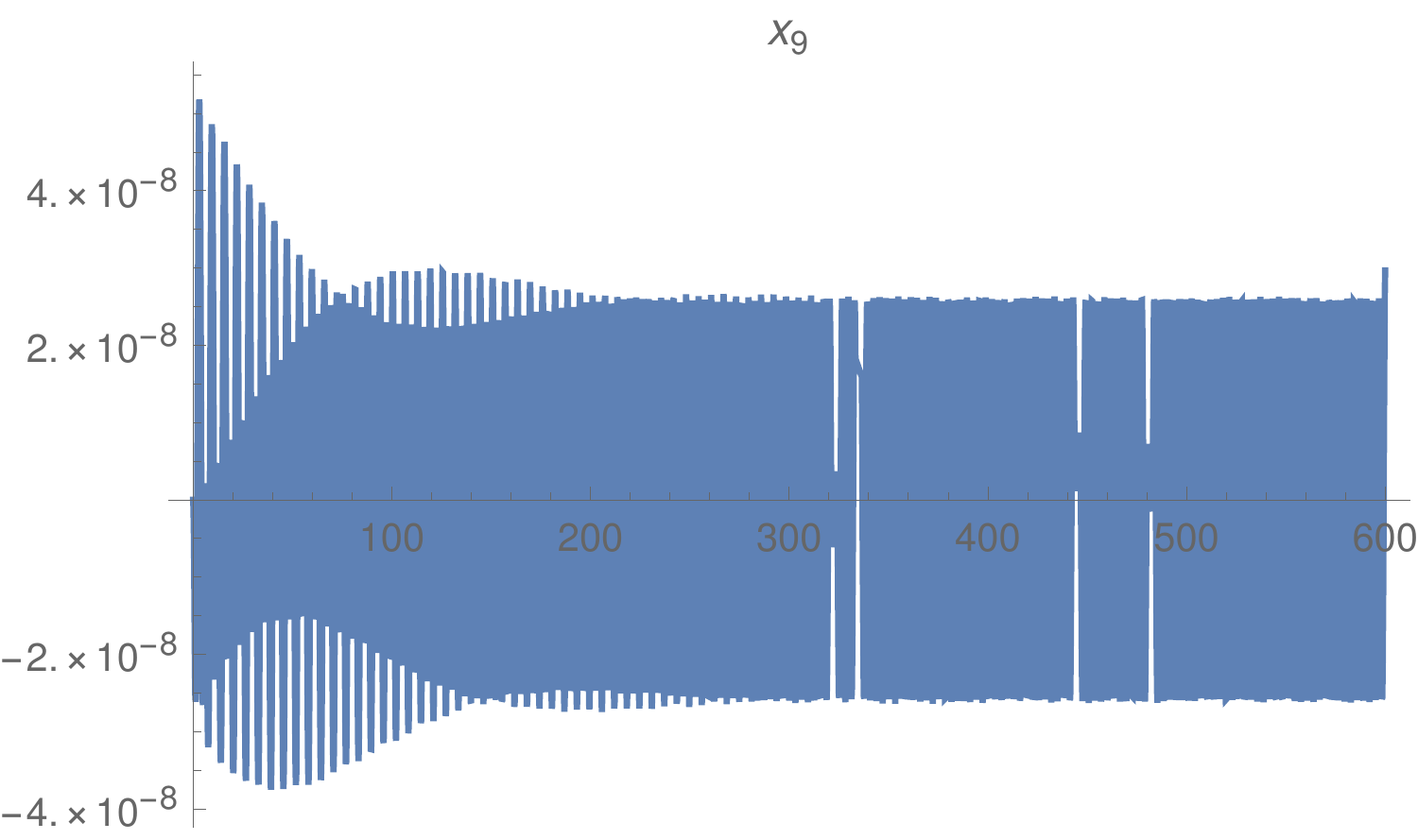}
\caption{Time evolution over the interval $[0,600]$ of the variables
  $x_{2}(t)$ and $x_{9}(t)$. Comparison between $r=2$ (top) and $r=3$ (bottom).}
\label{f.0Pi.comp_x}
\end{figure}

In Figures \ref{f.0Pi.comp_p}, \ref{f.0Pi.comp_q} and
\ref{f.0Pi.comp_x} we have compared the dynamics of some of the
variables of the linearly stable approximate periodic orbit
$q^*=(0,\pi)$ over a significant time interval of order
$\Oscr(\epsilon^{-1})$, according to the slowest frequency
$\omega=\sqrt3\epsilon$. As expected, the increase of the normal form
order provides a gain of the accuracy of the approximation at least of
a factor $\Oscr(\epsilon)$, both in the \emph{internal} and in the
\emph{transversal} dynamics. This is in agreement with the increase of
the order of the remainder of the normal form, which yields to the
estimate~\eqref{e.small.Ups}.

\begin{figure}[!ht]
\centering
\includegraphics[width=0.4\textwidth]{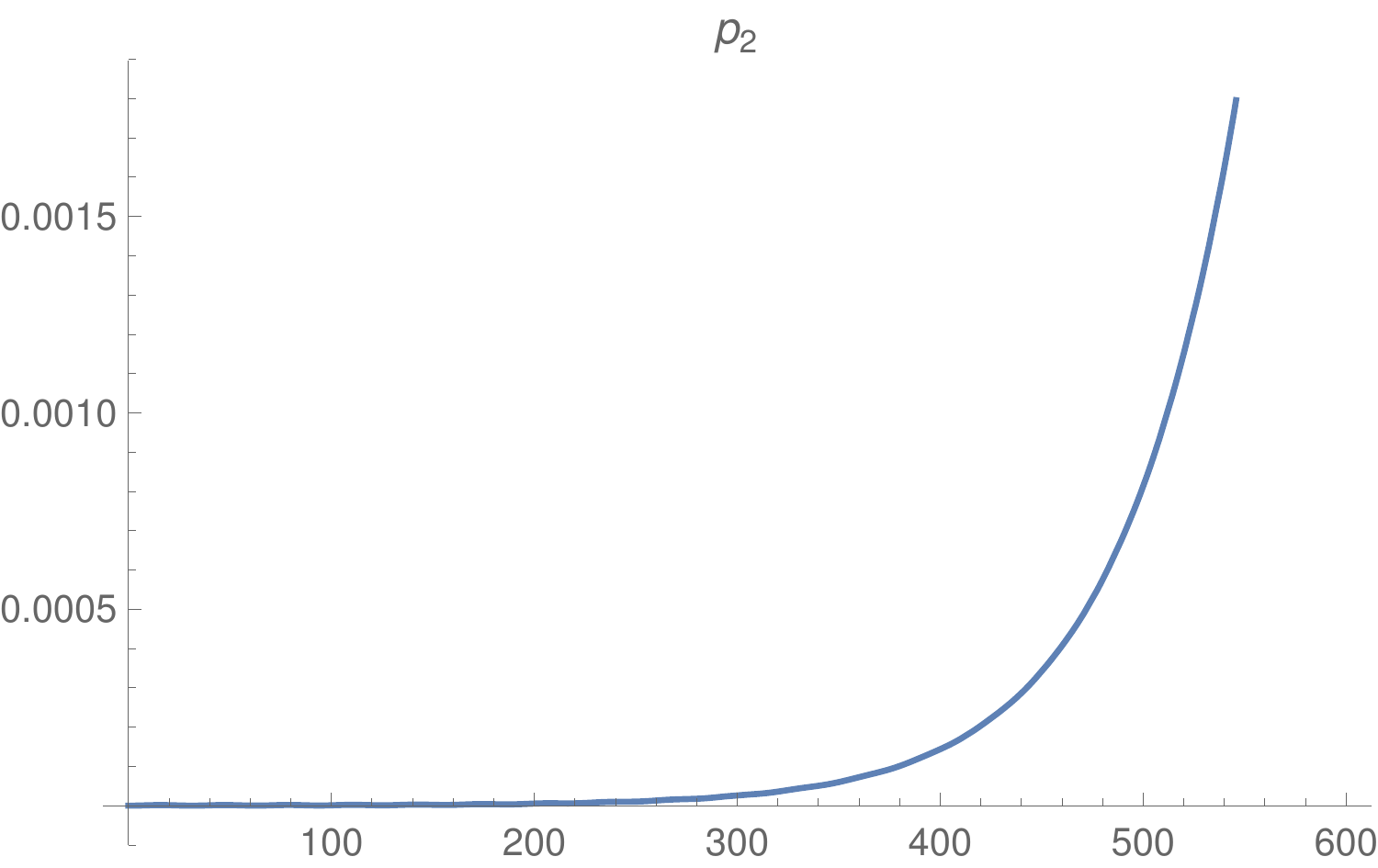}
\includegraphics[width=0.4\textwidth]{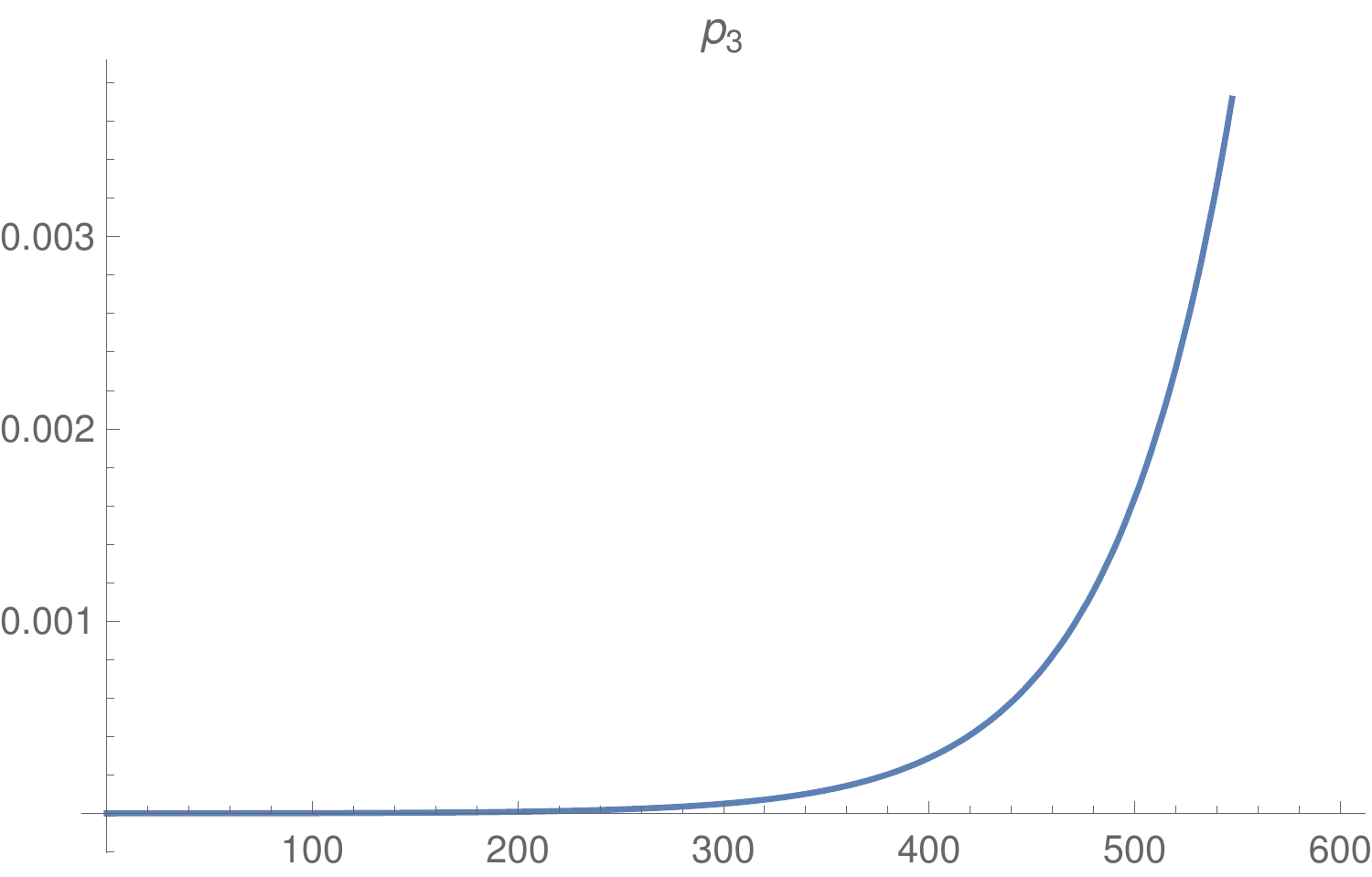}
\includegraphics[width=0.4\textwidth]{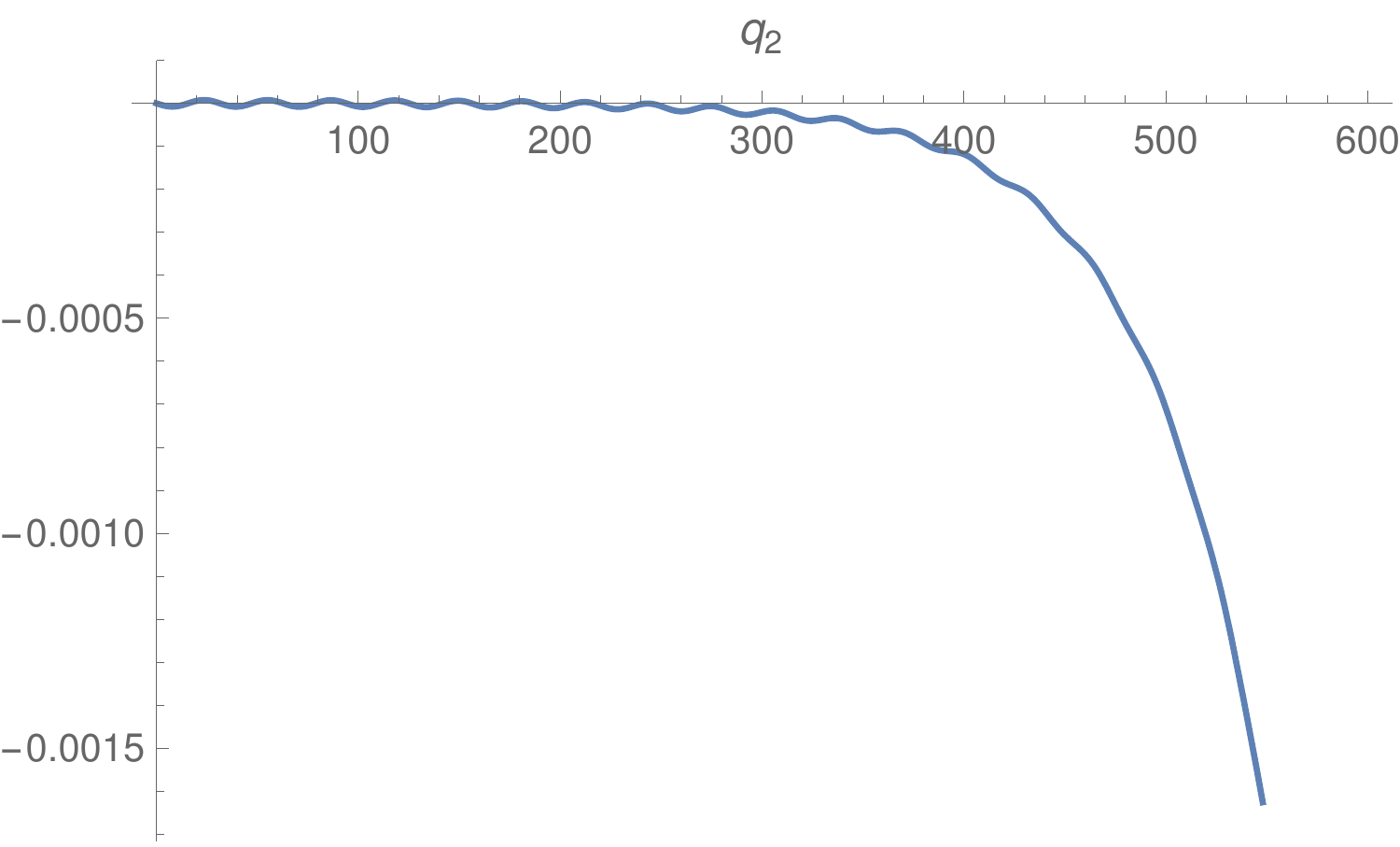}
\includegraphics[width=0.4\textwidth]{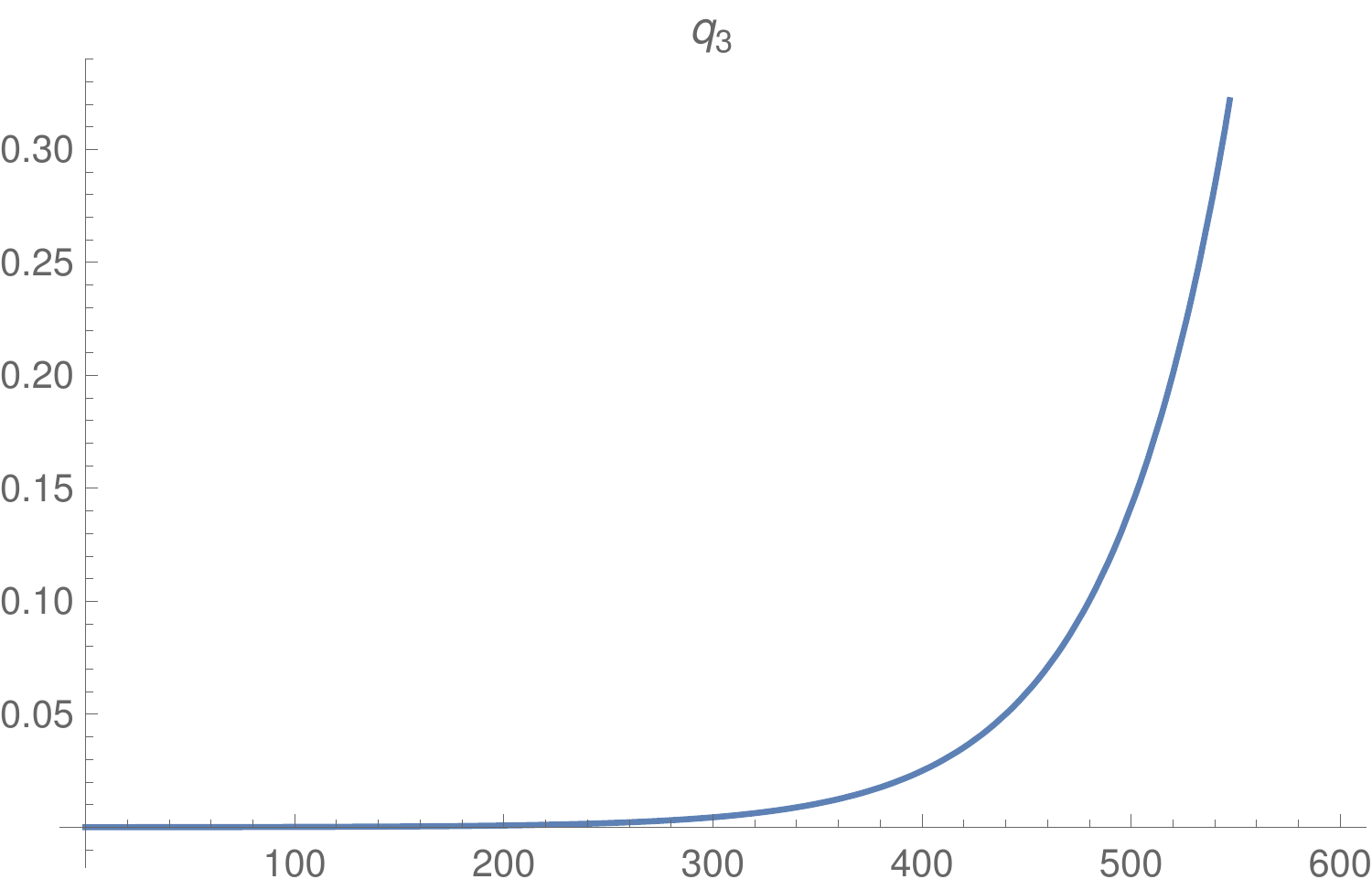}
\caption{Time evolution over the interval $[0,\Oscr(\epsilon^{-1})]$
  of the variables $p_{2}(t)$, $p_{3}(t)$ (top) and $q_{2}(t)$,
  $q_{3}(t)$ (bottom) for $q^*=(0,0)$.}
\label{f.00}
\end{figure}

\begin{figure}[!ht]
\centering
\includegraphics[width=0.4\textwidth]{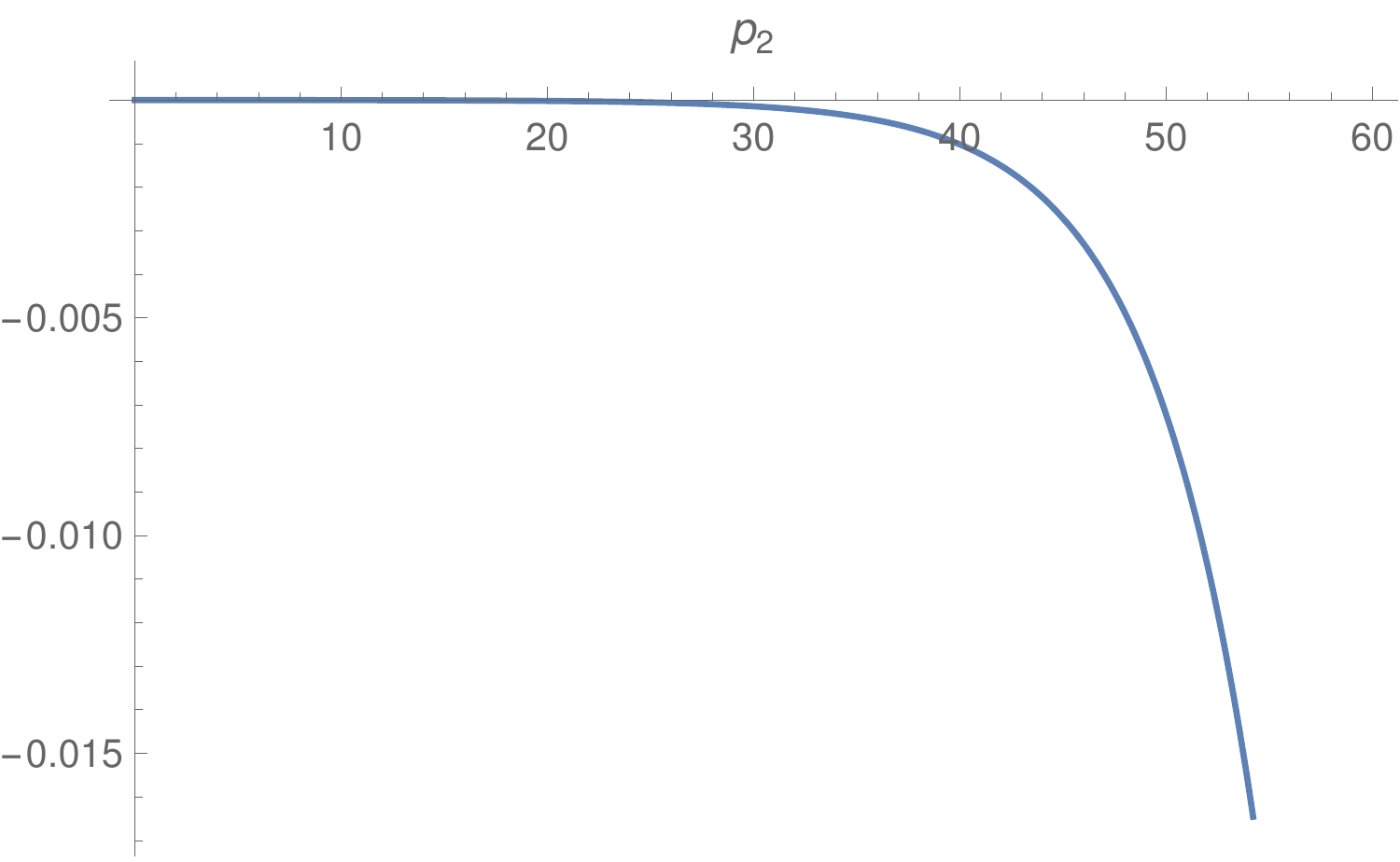}
\includegraphics[width=0.4\textwidth]{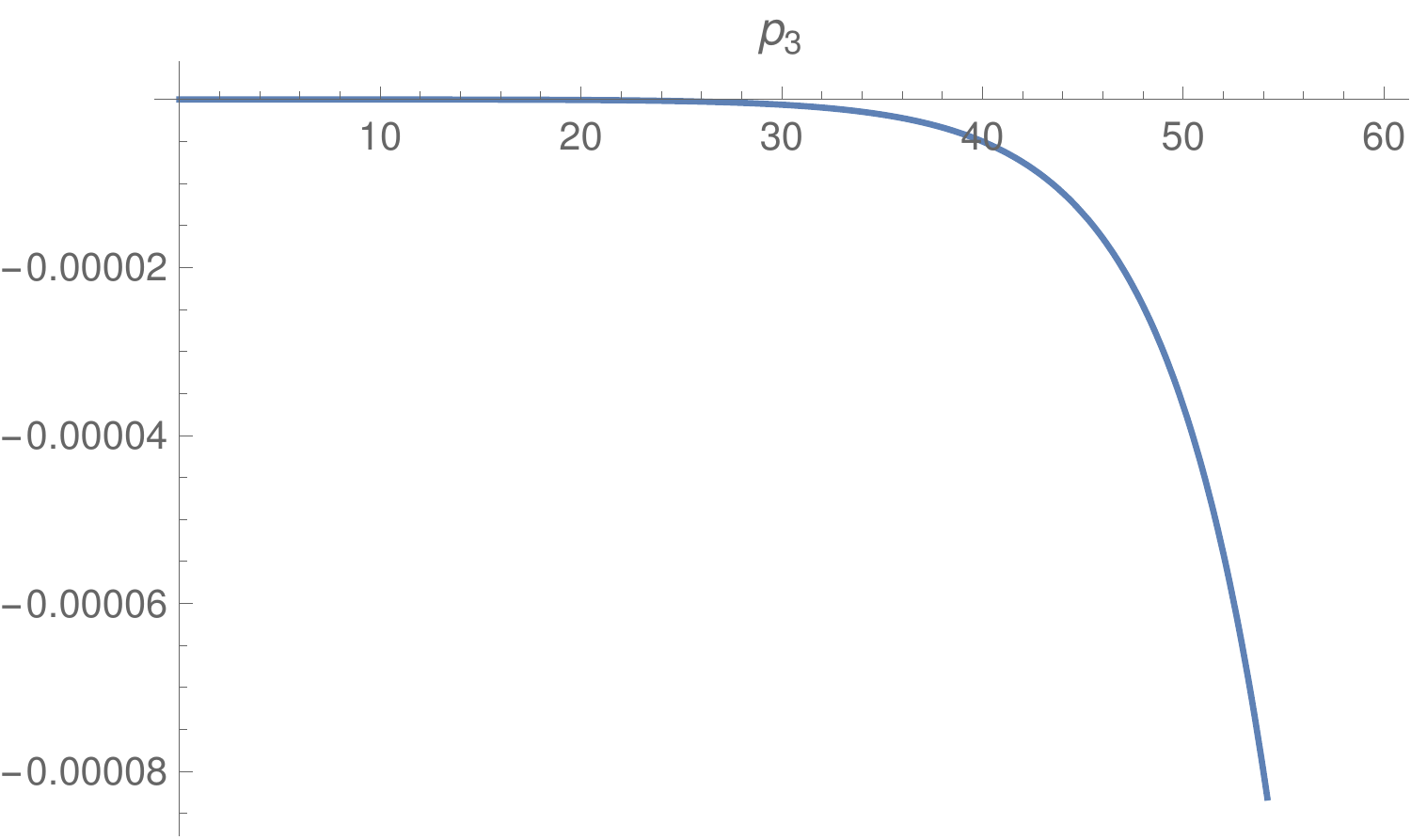}
\includegraphics[width=0.4\textwidth]{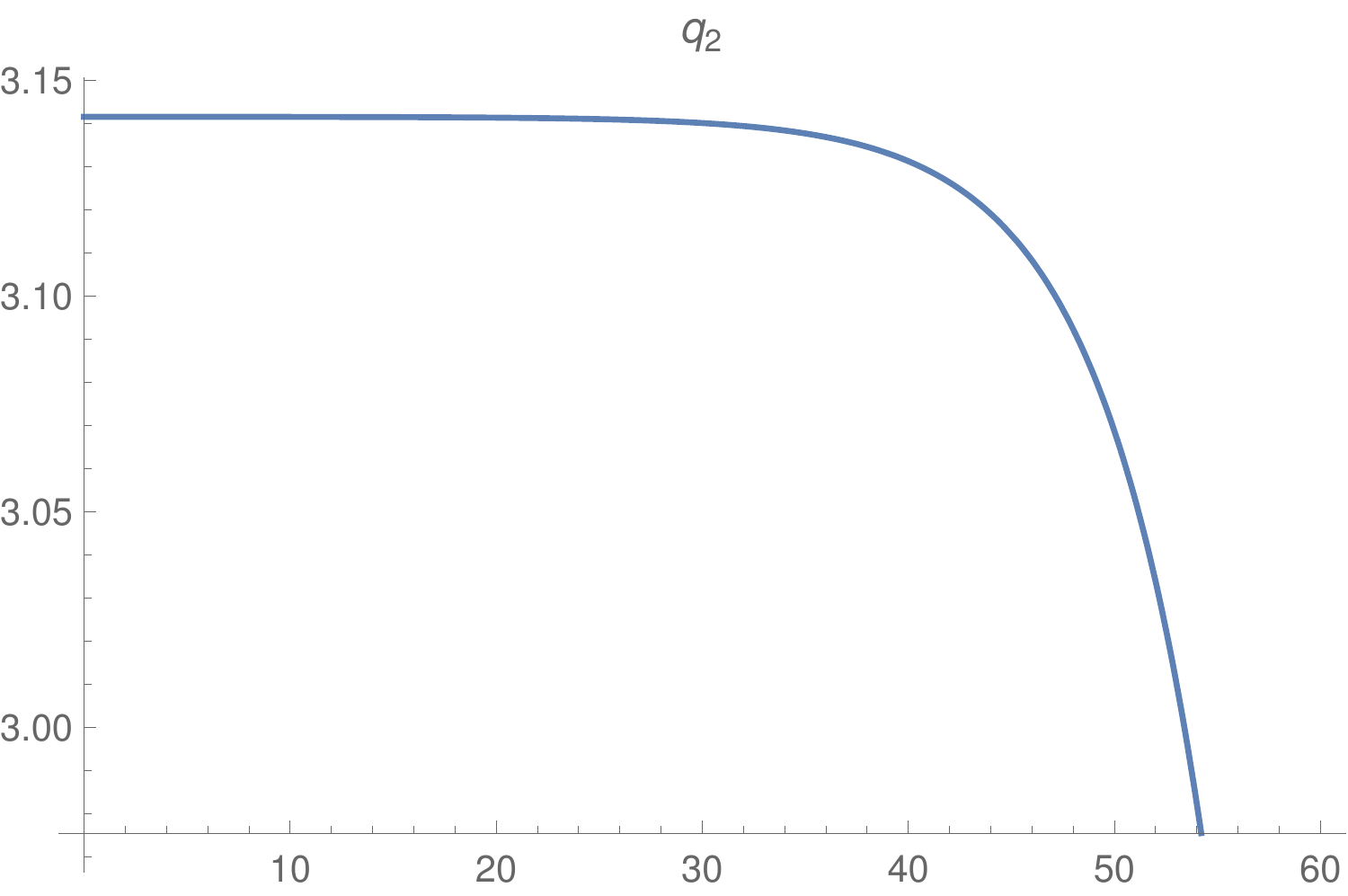}
\includegraphics[width=0.4\textwidth]{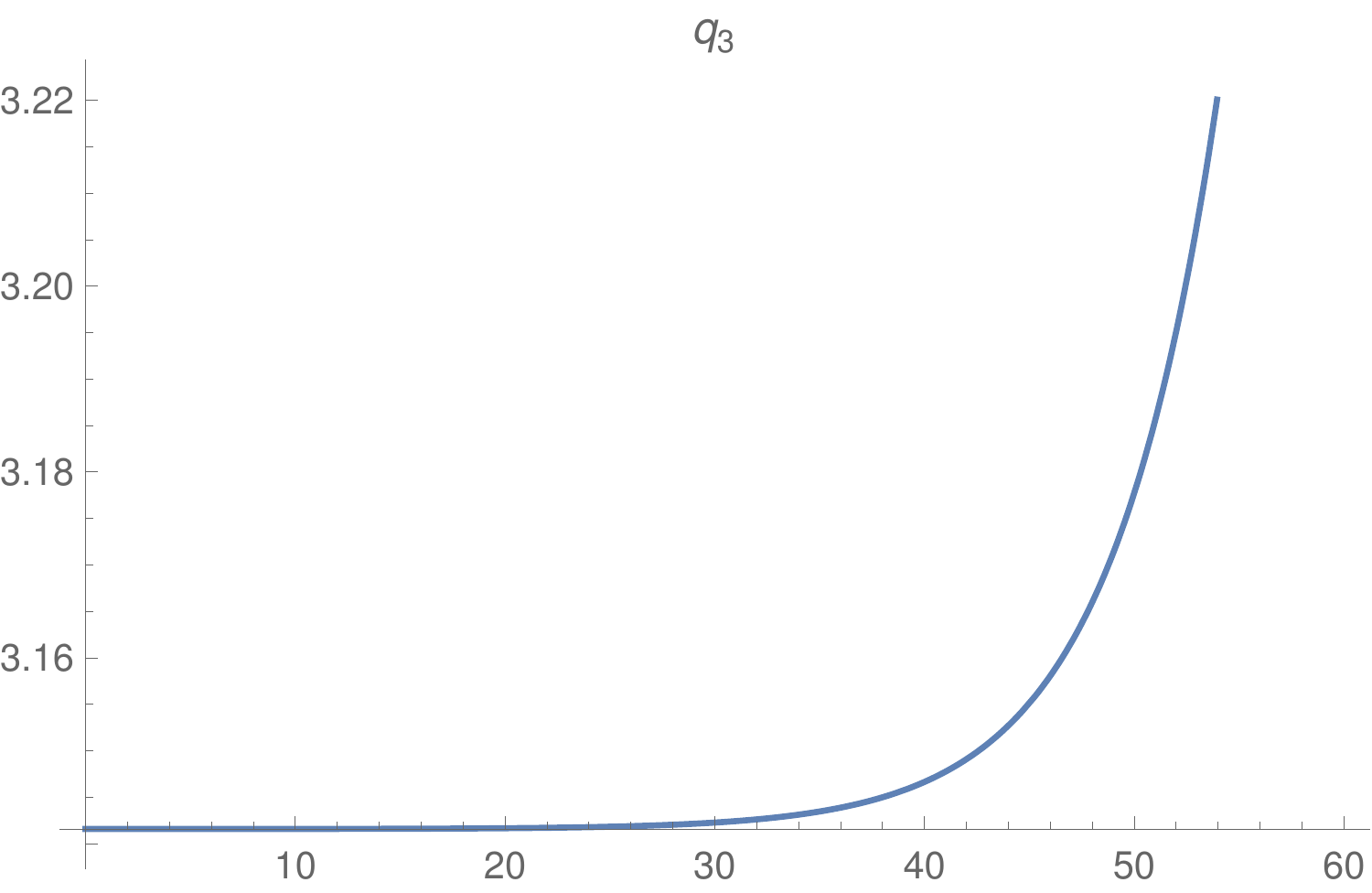}
\caption{Time evolution over the interval $[0,\Oscr(\epsilon^{-\frac12})]$
  of the variables $p_{2}(t)$, $p_{3}(t)$ (top) and $q_{2}(t)$, $q_{3}(t)$ (bottom) for
  $q^*=(\pi,\pi)$.}
\label{f.PiPi}
\end{figure}

\subsection{Linear stability of the approximate periodic orbit}

Looking at Figures \ref{f.0Pi.comp_p} and \ref{f.0Pi.comp_q},
independently of the normalisation order considered, it is evident
both the linear stability of the orbit, over the time interval
according to the slowest frequency, and the effect of the two
frequencies, which have a different scaling in $\epsilon$, as
illustrated in the previous section. In particular, the variables
$p_3$ and $q_3$ clearly show the periodic effect mainly due to the
frequency $\sqrt{3}\epsilon$, while the variables $p_2$ and $q_2$
clearly have a faster oscillation of the order of the frequency
$2\sqrt\epsilon$ (actually $p_2$ have a quasi-periodic dynamics which
exhibits both the frequencies).  In a similar way, Figures \ref{f.00}
and \ref{f.PiPi} show the effect of instability of the approximate
periodic orbit, over two different time scales. Indeed, in Figure
\ref{f.00} the exponentially fast separation from the approximate
periodic orbit requires time of order $\Oscr(\epsilon^{-1})$,
coherently with the real eigenvalues $\sqrt3\epsilon$ of the
$q^*=(0,0)$ solution. On the other hand, in Figure \ref{f.PiPi} the
departure is much faster and occurs already on a time scale of order
$\Oscr(\epsilon^{-\frac12})$, coherently with the real eigenvalues
$2\sqrt\epsilon$ of the $q^*=(\pi,\pi)$ solution.

\section{Conclusions}
\label{s:4}

In this paper we applied an abstract result on the break down of a
fully resonant torus in nearly integrable Hamiltonian system, to
revisit the existence of time periodic and spatially localised
solutions in dNLS lattices, such as discrete solitons or multi-pulse
solitons. We considered several different dNLS models, starting from
the standard one, moving to coupled dNLS chains (Zigzag or railway
models) up to model with a purely nonlinear interaction.  In all these
cases we showed that for $\epsilon$ small enough, i.e., in the limit
of small coupling, this class of solutions are at leading order
degenerate; hence a first order average is not conclusive. The normal
form scheme developed in~\cite{SanDPP20} and the main Theorems on
existence and linear stability there included, allow to investigate,
with the help of a computer algebra system, different kind of
degenerate configurations, thus confirming the practical applicability
of the abstract algorithm. At the same time, it allows to shed some
light on a wider class of localised periodic solutions, leading to the
(expected) existence of 2-dimensional resonant tori, thanks to the
action of the Gauge symmetry of the dNLS models. These are special
localised solutions that typically do not exist in Klein-Gordon
lattices: indeed, the presence of the full Fourier spectrum of the
unperturbed oscillators provides nondegenerate configurations even for
a resonant modulus $M_\omega$ different from the
$1:\text{\textendash}:1$ (see for example
\cite{KouI02,Kou04,PelS12}). Actually, the possibility to apply the
present approach to multibreathers in weakly coupled chain of
oscillators is limited by the need to explicitly transform the excited
oscillators to action-angle variables; this is a problem which might
be overcome with special choices of the nonlinear potential, like the
Morse potential, or with a preliminary dNLS normal form approximation
of the nonlinear lattice (as in \cite{PalP14,PelPP16}). Another
example of Hamiltonian Lattice where we expect that this approach
might led to interesting results is the FPU model. A recent work
\cite{CarL21} has shown how to deal with the original FPU model in
order to split the variables describing a low dimensional elliptic
invariant torus from the variables describing the transversal
dynamics; the same strategy might be adapted in order to study
completely resonant low dimensional tori and the corresponding
periodic orbits, possibly at the thermodynamic limit. Finally, a
different direction of future development could be to extend the
scheme in order to study the existence of degenerate quasi-periodic
solutions (degenerate KAM-subtori, as in \cite{Tre91}), both from an
abstract point of view and in terms of applications to physical
models.

\bigskip
\noindent
{\bf Acknowledgements}
M.S., T. P. and V. D. have been supported by the GNFM - Progetto
Giovani funding ``Low-dimensional Invariant Tori in FPU-like Lattices
via Normal Forms'' and by the MIUR-PRIN 20178CJA2B ``New Frontiers of
Celestial Mechanics: theory and Applications''.  We all thank Vassilis
Koukouloyannis for his visit to Milan in November 2019, which brought
to interesting discussions on the applications here illustrated.

\def\cprime{$'$} \def\i{\ii}\def\cprime{$'$} \def\cprime{$'$}

\end{document}